\documentclass{memo-l}

\usepackage{amssymb,amsfonts,amsmath}
\usepackage{enumerate}
\usepackage{psfrag}
\usepackage{comment}
\usepackage[dvips]{graphicx}

\newtheorem{theo}{Theorem}[chapter]
\newtheorem{prop}[theo]{Proposition} 
\newtheorem{lemma}[theo]{Lemma}
\newtheorem{coro}[theo]{Corollary}
\newtheorem{claim}[theo]{Claim} 

\newtheorem{defn}[theo]{Definition}
\newtheorem{construct}[theo]{Construction}

\numberwithin{section}{chapter}
\numberwithin{equation}{chapter}

\newcommand{\sm}{\setminus}
\newcommand{\eps}{\varepsilon} 
\newcommand{\T}{{\mathcal T}} 
\newcommand{\R}{\mathbb{R}} 
\newcommand{\Q}{{\mathbb Q}} 
\newcommand{\Z}{{\mathbb Z}} 
\newcommand{\C}{{\mathcal C}}

\newcommand{\Prob}{{\mathbb P}}
\newcommand{\Permute}{\textrm{Sym}}

\newcommand{\Sart}{{\mathcal S}}
\newcommand{\Part}{{\mathcal P}}
\newcommand{\Qart}{{\mathcal Q}}
\newcommand{\ab}{{\bf a}}
\newcommand{\ib}{{\bf i}}
\newcommand{\vb}{{\bf v}} 
\newcommand{\ub}{{\bf u}}
\newcommand{\xb}{{\bf x}}
\newcommand{\zb}{{\bf z}}
 
\newcommand{\cb}{{\bf c}} 

\newcommand{\0}{{\bf 0}}
\newcommand{\1}{{\bf 1}}

\newcommand{\mc}[1]{\mathcal{#1}}
\newcommand{\mb}[1]{\mathbb{#1}}
\newcommand{\nib}[1]{\noindent {\bf #1}}

\newcommand{\brac}[1]{\left( #1 \right)}

\newcommand{\bfl}[1]{\left\lfloor #1 \right\rfloor}

\newcommand{\bgen}[1]{\left\langle #1 \right\rangle}
\newcommand{\sub}{\subseteq}

\newcommand{\ov}{\overline}

\newcommand{\es}{\emptyset}

\newcommand{\GG}{\Gamma} 
\newcommand{\N}{{\mathbb N}}
\newcommand{\yb}{{\bf y}}
\newcommand{\unit}{{\bf u}}
\newcommand{\pl}{\partial}

\begin{document}

\pagenumbering{gobble} 

\chapter*{A geometric theory for hypergraph matching}

\hspace{2cm} Peter Keevash\footnote{
Mathematical Institute, University of Oxford, Oxford, UK. 
Email keevash@maths.ox.ac.uk.
Research supported in part by ERC grant 239696 and EPSRC grant EP/G056730/1.}
\hspace{2cm}
Richard Mycroft\footnote{School of Mathematics,  
University of Birmingham, 
Birmingham, UK. Email: r.mycroft@bham.ac.uk.}

\vspace{1cm}

Accepted for publication in Memoirs of the American Mathematical Society

\chapter*{Abstract}

We develop a theory for the existence of perfect matchings in hypergraphs under quite general conditions.  Informally speaking, the obstructions to perfect matchings are geometric, and are of two distinct types: `space barriers' from convex geometry, and `divisibility barriers' from arithmetic lattice-based constructions.  To formulate precise results, we introduce the setting of simplicial complexes with minimum degree sequences, which is a generalisation of the usual minimum degree condition. We determine the essentially best possible minimum degree sequence for finding an almost perfect matching. Furthermore, our main result establishes the stability property: under the same degree assumption, if there is no perfect matching then there must be a space or divisibility barrier. This allows the use of the stability method in proving exact results.  Besides recovering previous results, we apply our theory to the solution of two open problems on hypergraph packings: the minimum degree threshold for packing tetrahedra in $3$-graphs, and Fischer's conjecture on a multipartite form of the Hajnal-Szemer\'edi Theorem. Here we prove the exact result for tetrahedra and the asymptotic result for Fischer's conjecture; since the exact result for the latter is technical we defer it to a subsequent paper. 

\tableofcontents

\mainmatter
\chapter{Introduction} \label{sec:intro}

Hypergraph matchings\footnote{A \emph{hypergraph} $G$ consists of a vertex set $V$ and an edge set $E$, where each edge $e \in E$ is a subset of $V$. We say $G$ is a \emph{$k$-graph} if every edge has size $k$. A \emph{matching} $M$ is a set of vertex disjoint edges in $G$. We call $M$ \emph{perfect} if it covers all of $V$.} provide a general framework for many important Combinatorial problems. Two classical open problems of this nature are the question of whether there exist designs of arbitrary strength, and Ryser's conjecture that every Latin square of odd order has a transversal; these are both equivalent to showing that some particular hypergraphs have perfect matchings. Furthermore, matchings are also an important tool for many practical questions, such as the `Santa Claus' allocation problem (see~\cite{AFS}). However, while Edmonds' algorithm~\cite{E} provides an efficient means to determine whether a graph has a perfect matching, the decision problem is NP-complete in $k$-graphs for $k \ge 3$ (it is one of Karp's original 21 NP-complete problems~\cite{Ka}). Thus we do not expect a nice characterisation, so we concentrate on natural sufficient conditions for finding a perfect matching. 

One natural hypergraph parameter that is widely considered in the literature is the \emph{minimum degree} $\delta(G)$, which is the largest number $m$ such that every set of $k-1$ vertices is contained in at least $m$ edges of $G$. What is the minimum degree threshold for finding a perfect matching? (We assume that $k|n$, where $n=|V|$.)  For graphs ($k=2$) a straightforward greedy argument shows that the threshold is $n/2$ (or one can deduce it from Dirac's theorem \cite{D}, which states that the same threshold even gives a Hamilton cycle). The general case was a long-standing open problem, finally resolved by R\"odl, Ruci\'nski and Szemer\'edi~\cite{RRS}, who determined the threshold precisely for large $n$: it is $n/2-k+C$, where $C \in \{1.5,2,2.5,3\}$ depends on arithmetic properties of $k$ and $n$. There is a large literature on minimum degree problems for hypergraphs, see e.g.~\cite{AFHRRS, Cs, CK, DH, D, F, HPS, KKMO, KM, Kh2, Kh, KO, KO1, KOT, LM, LM2, MM, MR2, MS, P, RRS3, RRS2, RRS, TZ} and the survey by R\"odl and Ruci\'nski~\cite{RR} for details.

\section{Space barriers and divisibility barriers.}
To motivate the results of this paper it is instructive to consider the extremal examples for the minimum degree problems. Consider a graph $G_1$ on $n$ vertices whose edges are all pairs incident to some set $S$ of size $n/2-1$. Then $\delta(G_1)=n/2-1$, and $G_1$ has no perfect matching, as each edge of a matching $M$ uses a vertex in $S$, so $|M| \le |S|$; we say that $G_1$ has a {\em space barrier}. Now suppose $n/2$ is odd and consider a graph $G_2$ on $n$ vertices consisting of two disjoint complete graphs of size $n/2$. Then $\delta(G_2)=n/2-1$, and $G_2$ has no perfect matching, but for the different reason that edges have even size; we say that $G_2$ has a {\em divisibility barrier}. 

While these two examples are equally good for graphs, for general $k$-graphs a separation occurs. For $G_1$ we take a $k$-graph whose edges are all $k$-tuples incident to some set $S$ of size $n/k-1$; this satisfies $\delta(G_1)=n/k-1$. For $G_2$ we take a $k$-graph whose edges are all $k$-tuples that have an even size intersection with some set $S$ such that $n/2-1 \le |S| \le (n+1)/2$ and $|S|$ is odd; this satisfies $\delta(G_1)=n/2-k+C$, being one of the extremal constructions in the result of \cite{RRS} mentioned above. One should also note that space is a robust obstruction to matchings, in that the size of a maximum matching in $G_1$ may be decreased by decreasing $|S|$, whereas divisibility is not robust, in the sense that any construction similar to $G_2$ has an almost perfect matching. In fact, R\"odl, Ruci\'nski and Szemer\'edi~\cite{RRS2} showed that the minimum degree threshold for a matching of size $n/k-t$ with $t \ge k-2$ is $n/k-t$.  Thus we see a sharp contrast between the thresholds for perfect matching and almost perfect matching.

The main message of this paper is that space and divisibility are the determining factors for perfect matchings in hypergraphs under quite general conditions, and that these factors are inherently geometric. The first part of the geometric theory was anticipated by a result on fractional perfect matchings in~\cite{RRS2} that we generalise here. The key point is that fractional perfect matchings correspond to representing a constant vector as a convex combination of edge vectors, so non-existence of fractional perfect matchings can be exploited in geometric form via separating hyperplanes. Furthermore, the fractional problem has bearing on the original problem through standard `regularity machinery', which converts a fractional solution into a matching covering all but $o(n)$ vertices for large $n$. 

The second part of the theory is to understand when the number of uncovered vertices can be reduced to a constant independent of $n$. The idea here can only be properly explained once we have described the regularity embedding strategy, but a brief summary is as follows. Firstly, the $o(n)$ uncovered vertices arise from imbalances created when converting from the fractional to the exact solution. Secondly, the possible `transferrals' of imbalances can be expressed geometrically by defining an appropriate polyhedron and testing for a ball around the origin of some small constant radius; this can also be understood in terms of separating hyperplanes.  Thus the first two parts of the theory are problems of convex geometry, which correspond to space barriers. 

The third part of the theory concerns divisibility barriers, which determine when the number of uncovered vertices can be reduced from a constant to zero. We will see that in the absence of space barriers, perfect matchings exist except in hypergraphs that are structurally close to one of a certain class of arithmetic constructions defined in terms of lattices in $\mb{Z}^d$ for some $d \le k$.  Furthermore, since the constructions with space or divisibility barriers do not have perfect matchings, in a vague sense we have `the correct theory', although this is not a precise statement because of additional assumptions. Our theory is underpinned by the `strong' hypergraph regularity theory independently developed by Gowers~\cite{G1} and R\"odl et al. \cite{FR,NR,RS,RSk}, and the recent hypergraph blowup lemma of Keevash~\cite{K}.  Fortunately, this part of the argument is mostly covered by existing machinery, so the majority of this paper is devoted to the geometric theory outlined above.

To formulate precise results, we introduce the setting of simplicial complexes with minimum degree sequences, which is a generalisation of the minimum degree condition previously considered. In this setting our main theorems (stated in Section~\ref{sec:pm}) give minimum degree sequences that guarantee a perfect matching for hypergraphs that are not close to a lattice construction. These minimum degree sequences are best possible, and furthermore have the `stability' property that, unless the hypergraph is structurally close to a space barrier construction, one can find a perfect matching even with a slightly smaller degree sequence. We defer the statements until we have given the necessary definitions in the next chapter. For now we want to emphasise the power of this framework by describing its application to the solutions of two open problems on packings. Suppose $G$ is a $k$-graph on $n$ vertices and $H$ is a $k$-graph on $h$ vertices, where $h|n$ (we think of $h$ as fixed and $n$ as large). An \emph{$H$-packing} in $G$ is a set of vertex-disjoint copies of $H$ inside $G$; it is \emph{perfect} if there are $n/h$ such copies, so that every vertex of $G$ is covered. In the case when $H$ is a single edge we recover the notion of (perfect) matchings. 

As for matchings we have the natural question: what is the minimum degree threshold for finding a perfect $H$-packing? Even for graphs, this is a difficult question with a long history. One famous result is the Hajnal-Szemer\'edi theorem~\cite{HS}, which determines the threshold for the complete graph $K_r$ on $r$ vertices: if $r|n$ and $\delta(G) \geq (r-1)n/r$ then $G$ has a $K_r$-packing, and this bound is best possible.  It is interesting to note that this is essentially the same threshold as for the Tur\'an problem~\cite{T} for $K_{r+1}$, i.e.\ finding a single copy of the complete graph with one more vertex. The perfect packing problem for general graphs $H$ was essentially solved by K\"uhn and Osthus~\cite{KO1}, who determined the threshold for large $n$ up to an additive constant $C(H)$. The precise statement would take us too far afield here, but we remark that the threshold is determined by either space or divisibility barriers, and that the dependence on the chromatic number of $H$ continues a partial analogy with Tur\'an problems. We refer the reader to their survey~\cite{KO2} for further results of this type.

\section{Tetrahedron packings.}
For hypergraph packing problems, the natural starting point is the tetrahedron $K^3_4$, i.e.\ the complete $3$-graph on $4$ vertices. Here even the asymptotic existence threshold is a long-standing open problem; this is an important test case for general hypergraph Tur\'an problems, which is a difficult area with very few known results (see the survey by Keevash~\cite{K}). In light of this, it is perhaps surprising that we are able here to determine the tetrahedron packing threshold for large $n$, not only asymptotically but precisely.  One should note that the two problems are not unrelated; indeed Tur\'an-type problems for the tetrahedron are required when showing that there are no divisibility barriers (but fortunately they are more tractable than the original problem!) The extremal example for the perfect tetrahedron packing problem is by no means obvious, and it was several years after the problem was posed by Abbasi (reported by Czygrinow and Nagle \cite{CN}) that Pikhurko~\cite{P} provided the optimal construction (we describe it in Chapter~\ref{sec:tetra}). Until recently, the best published upper bounds, also due to Pikhurko~\cite{P}, were $0.8603\dots n$ for the perfect packing threshold and $3n/4$ for the almost perfect packing threshold.  More recent upper bounds for the perfect packing threshold are $4n/5$ by Keevash and Zhao (unpublished) and $(3/4+o(1))n$ independently by Keevash and Mycroft (earlier manuscripts of this paper) and by Lo and Markstr\"om \cite{LM} (posted online very recently). It is instructive to contrast the `absorbing technique' used in the proof of \cite{LM} with the approach here; we will make some remarks later to indicate why our methods seem more general and are able to give the exact result, which is as follows.

\medskip

\begin{theo}~\label{TETRAPACK}
There exists $n_0$ such that if $G$ is a 3-graph on $n \geq n_0$ vertices such that $4 \mid n$ and
$$\delta(G) \geq
\begin{cases}
3n/4 - 2 & \textrm{ if } 8 \mid n\\
3n/4 - 1 & \textrm{ otherwise,}
\end{cases}$$
then $G$ contains a perfect $K^3_4$-packing. This minimum degree bound is best possible.
\end{theo}

The minimum degree bound of Theorem~\ref{TETRAPACK} is best possible. Indeed, consider a $3$-graph $G$ whose vertex set $V$ is the disjoint union of sets $A$, $B$, $C$ and $D$ whose sizes are as equal as possible with $|A|$ odd. The edges of $G$ are all $3$ tuples except those
\begin{enumerate}[(i)]
\item with all vertices in $A$,
\item with one vertex in $A$ and the remaining two vertices in the same vertex class, or
\item with one vertex in each of $B$, $C$ and $D$.
\end{enumerate}
There is then no perfect $K^3_4$-packing in $G$ (see Proposition~\ref{tetrahedrapackextrex} for details), but $\delta(G)$ is equal to $3n/4 - 3$ if $8 \mid n$, and $3n/4 - 2$ otherwise.

\section{A multipartite Hajnal-Szemer\'edi theorem.}

Our second application is to a conjecture of Fischer~\cite{F} on a multipartite form of the Hajnal-Szemer\'edi Theorem.
Suppose $V_1, \dots, V_k$ are disjoint sets of $n$ vertices each, and $G$ is a $k$-partite graph on vertex classes
$V_1, \dots, V_k$ (that is, $G$ is a graph on $V_1 \cup \dots \cup V_k$ such that no edge of $G$ has both vertices in
the same $V_j$). Then we define the \emph{partite minimum degree} of $G$, denoted $\delta^*(G)$, to be the largest $m$
such that every vertex has at least $m$ neighbours in each part other than its own, i.e.\
$$\delta^*(G) = \min_{i \in [k]} \min_{v \in V_i} \min_{j \in [k]\sm\{i\}} |N(v) \cap V_j|,$$
where $N(v)$ denotes the neighbourhood of $v$. Fischer conjectured that if $\delta^*(G) \ge (k-1)n/k + 1$ then $G$ has a perfect $K_k$-packing (actually his original conjecture did not include the $+1$, but this stronger conjecture is known to be false). The case $k=2$ of the conjecture is an immediate corollary of Hall's Theorem, whilst the cases $k=3$ and $k=4$ were proved by Magyar and Martin~\cite{MM} and Martin and Szemer\'edi~\cite{MS} respectively. Also, Csaba and Mydlarz~\cite{Cs} proved a weaker version of the conjecture in which the minimum degree condition has an error term depending on $k$. The following theorem, an almost immediate corollary of our results on hypergraph matchings, gives an asymptotic version of this conjecture for any $k$.

\begin{theo} \label{partitehajnalszem}
For any $k$ and $c > 0$ there exists $n_0$ such that if $G$ is a $k$-partite graph with parts $V_1,\dots,V_k$ of size $n \geq n_0$ and $\delta^*(G) \ge (k-1)n/k + cn$, then $G$ contains a perfect $K_k$-packing.
\end{theo}

This asymptotic result was also proven independently and simultaneously by Lo and Markstr\"om~\cite{LM2} in another application of the `absorbing technique'. As with Theorem~\ref{TETRAPACK}, by considering the  near-extremal cases of the conjecture using the `stability property' of our main theorem, we are able to prove an exact result in~\cite{KM}, namely that Fischer's conjecture holds for any sufficiently large $n$. However, this stability analysis is 
lengthy and technical, so we prefer to divest this application from the 
theory developed in this paper.

As mentioned above, we will deduce both Theorems~\ref{TETRAPACK} and~\ref{partitehajnalszem} from a general framework of matchings in simplicial complexes. These will be formally defined in the next chapter, but we briefly indicate the connection here. For Theorem~\ref{TETRAPACK} we consider the `clique $4$-complex', with the tetrahedra in $G$ as $4$-sets, $G$ as $3$-sets, and all smaller sets; for Theorem~\ref{partitehajnalszem} we consider the `clique $k$-complex', with the $j$-cliques of $G$ as $j$-sets for $j \le k$. In both cases, the required perfect packing is equivalent to a perfect matching using the highest level sets of the clique complex.

\section{Algorithmic aspects of hypergraph matchings.}

As described earlier, the decision problem of whether a $k$-graph $H$ has a perfect matching is NP-complete for $k \geq 3$, motivating our consideration of the minimum degree which guarantees that $H$ contains a perfect matching. Another natural question to ask is for the minimum-degree condition which renders the decision problem tractable. That is, let PM$(k, \delta)$ denote the problem of deciding whether a $k$-graph $H$ on $n$ vertices (where $k \mid n$) with $\delta(H) \geq \delta n$ contains a perfect matching. The result of Karp~\cite{Ka} mentioned earlier shows that PM$(k,0)$ is NP-complete. On the other hand, the theorem of R\"odl, Ruci\'nski and Szemer\'edi~\cite{RRS} described earlier shows that any sufficiently large $k$-graph on $n$ vertices with $\delta(H) \geq n/2$ contains a perfect matching, so PM$(k, \delta)$ can be solved in constant time for $\delta \geq 1/2$ by an algorithm which simply says `yes' if $n$ is sufficiently large, and checks all possible matchings by brute force otherwise.

This question was further studied by Szyma\'nska~\cite{S}, who proved that for $\delta < 1/k$ the problem PM$(k, 0)$ admits a polynomial-time reduction to PM$(k, \delta)$, and so PM$(k, \delta)$ is NP-complete for such $\delta$. In the other direction, Karpi\'nski, Ruci\'nski and Szyma\'nska~\cite{KRS} proved that there exists a constant $\eps > 0$ such that PM$(k, 1/2 - \eps)$ is in P. This leaves a hardness gap for PM$(k, \delta)$ when $\delta \in [1/k, 1/2 - \eps]$.

The connection with the work of this paper is that we can check in polynomial time whether or not the edges of a $k$-graph $H$ satisfy arithmetic conditions of the types which define our notion of a divisibility barrier. We will see that for $\delta>1/k$, any sufficiently large $k$-graph $H$ with $\delta(H) \ge 
\delta |V(H)|$ cannot be close to a space barrier, so our main theorem 
will imply that either $H$ contains a perfect matching or $H$ is close to 
a divisibility barrier.
So to decide PM$(k, \delta)$ for $\delta > 1/k$ it suffices to decide the existence in a perfect matching when $H$ is close to a divisibility barrier. Unfortunately, in our main theorem our notion of `close' means an edit-distance of $o(n^k)$, which is too large to be checked by a brute-force approach. However, for $k$-graphs $H$ of codegree close to $n/2$ we are able to refine our main theorem to prove the following result, which states that either $H$ admits a perfect matching or \emph{every} edge of $H$ satisfies a divisibility condition of a type which defines a divisibility barrier.

\begin{theo} \label{pmstrongstab}
For any $k \geq 3$ there exists $c > 0$ and $n_0$ such that for any $n \geq n_0$ with $k \mid n$ and any $k$-graph $H$ on $n$ vertices with $\delta(H) \geq (1/2 - c)n$ the following statement holds. $H$ does not contain a perfect matching if and only if there exists a partition of $V(H)$ into parts $V_1,V_2$ of size at least $\delta(G)$ and $a \in \{0,1\}$ so that $|V_1| \neq an/k$ mod $2$ and $|e \cap V_1|=a$ mod $2$ for all edges $e$ of $G$.
\end{theo}

It is not hard to construct an algorithm with polynomial running time which checks the existence of such a partition, and so we recover the aforementioned result of Karpi\'nski, Ruci\'nski and Szyma\'nska that PM$(k, 1/2 - c)$ is in P. Together with Knox~\cite{KKM}, we were able to refine the methods of this paper to prove stronger results for $k$-graphs of large minimum degree to almost eliminate the hardness gap referred to above. Indeed, we show that we may replace the condition $\delta(H) \geq (1/2 - c)n$ in Theorem~\ref{pmstrongstab} by the condition $\delta(H) \geq n/3 + o(n)$, which shows that PM$(k, \delta)$ is in P for any $\delta > 1/3$. Furthermore, we are able to then extend Theorem~\ref{pmstrongstab} by further reducing the bound on $\delta(H)$; although this requires a significantly more complicated statement, we obtain a polynomial-time algorithm which decides PM$(k, \delta)$ for $\delta > 1/k$. Together with the work of Szyma\'nska described above, this settles the complexity status of PM$(k, \delta)$ for any $\delta \neq 1/k$. However, these refinements of our results only apply to $k$-graphs of large minimum codegree, whilst the results of this paper are much more general. 

\section{Notation.} The following notation is used throughout the paper: $[k]=\{1, \dots, k\}$; if $X$ is a set then $\binom{X}{k}$ is the set of subsets of $X$ of size $k$, and $\binom{X}{\le k} = \bigcup_{i \le k} \binom{X}{i}$ is the set of subsets of $X$ of size at most $k$; $o(1)$ denotes a term which tends to zero as~$n$ tends to infinity; $x \ll y$ means that for every $y > 0$ there exists some $x_0 > 0$ such that the subsequent statement holds for any $x < x_0$ (such statements with more variables are defined similarly). We write $x = y \pm z$ to mean $y-z \leq x \leq y+z$. Also we denote all vectors in bold font, and their coordinates by subscripts of the same character in standard font, e.g. $\ab = (a_1, \dots, a_n)$.

\chapter{Results and examples} \label{sec:results}

In this chapter we state our main theorems on perfect matchings in simplicial complexes.
For almost perfect matchings it requires no additional work to obtain more general results
that dispense with the `downward closure' assumption. However, for perfect matchings it is
more convenient to assume downward closure, which seems to hold in any natural application
of our results, so we will stick to simplicial complexes, and make some remarks later on
how the approach may be generalised. We also discuss several examples
that illustrate various aspects of the theory: space barriers and tightness of the degree
condition, lattice-based constructions, and generalisations of previous results.

\section{Almost perfect matchings.} \label{sec:almostpm}
We start with some definitions. We identify a hypergraph $H$ with its edge set,
writing $e \in H$ for $e \in E(H)$, and $|H|$ for $|E(H)|$.
A \emph{$k$-system} is a hypergraph $J$ in which every edge of $J$ has size at most $k$
and $\es \in J$. We refer to the edges of size $r$ in $J$
as \emph{$r$-edges of~$J$}, and write $J_r$ to denote the $r$-graph on $V(J)$ formed by these edges. It may be helpful to think of the $r$-graphs $J_r$ as different `levels' of the $k$-system $J$.
A \emph{$k$-complex} $J$ is a $k$-system whose edge set is closed under inclusion,
i.e.\ if $e \in H$ and $e' \sub e$ then $e' \in H$. That is, each level of $J$ is supported by the levels below it.
For any non-empty $k$-graph $G$, we may generate a $k$-complex $G^\le$
whose edges are any $e \sub V(G)$ such that $e \sub e'$ for some edge $e' \in G$.

We introduce the following notion of degree in a $k$-system $J$. For any edge $e$ of $J$,
the \emph{degree} $d(e)$ of $e$ is the number of $(|e|+1)$-edges $e'$ of $J$
which contain $e$ as a subset. (Note that this is {\em not} the standard notion of degree
used in $k$-graphs, in which the degree of a set is the number of edges containing it.)
The {\em minimum $r$-degree} of $J$, denoted $\delta_r(J)$,
is the minimum of $d(e)$ taken over all $r$-edges $e \in J$.
So every $r$-edge of $J$ is contained in at least $\delta_r(J)$ of the $(r+1)$-edges of $J$.
Note that if $J$ is a $k$-complex then $\delta_r(J) \leq \delta_{r-1}(J)$ for each $r \in [k-1]$.
The \emph{degree sequence of $J$} is 
$$\delta(J) = (\delta_0(J), \delta_1(J), \dots, \delta_{k-1}(J)).$$

Our minimum degree assumptions will always be of the form $\delta(J) \ge \ab$ pointwise for some vector
$\ab = (a_0,\dots,a_{k-1})$, i.e.\ $\delta_i(J) \ge a_i$ for $0 \le i \le k-1$. It is helpful to interpret
this `dynamically' as follows: when constructing an edge of $J_k$ by greedily choosing one vertex at
a time, there are at least $a_i$ choices for the $(i+1)$st vertex (this is the reason for the requirement that $\emptyset \in J$, which we need for the first choice in the process). To see that this generalises
the setting of minimum degree in hypergraphs, consider any $k$-graph $G$ and let $J$ be the $k$-system on $V(G)$
with $J_k = G$ and complete lower levels $J_i = \binom{V(G)}{i}$, $0 \le i \le k-1$.
Then $\delta(J) = (n,n-1,\dots,n-k+2,\delta(G))$ is a degree sequence in which the first $k-1$ coordinates
are as large as possible, so any minimum degree condition for this $k$-system $J$ reduces
to a minimum degree condition for $G$. More generally, for any $i \ge 0$, define an {\em $i$-clique} of $G$
to be a set $I \sub V(G)$ of size $i$ such that every $K \sub I$ of size $k$ is an edge of $G$ (this last
condition is vacuous if $i<k$). We can then naturally define a \emph{clique $r$-complex} of $G$, whose edges are the cliques in $G$.

\begin{defn}[Clique $r$-complex] 
The \emph{clique $r$-complex $J_r(G)$} of a $k$-graph $G$ is defined by taking $J_r(G)_i$ to consist of the $i$-cliques of $G$ for $0 \le i \le r$.
\end{defn}

To see the power in the extra generality of the degree sequence setting, consider the $4$-graph of tetrahedra in a $3$-graph $G$. This does not satisfy any minimum degree condition
as a $4$-graph, as non-edges of $G$ are not contained in any $4$-edges; on the other hand, we will see later that
a minimum degree condition on $G$ implies a useful minimum degree condition for the clique $4$-complex.

Our first theorem on hypergraph matchings gives a sufficient minimum degree sequence for a $k$-system to contain a matching covering almost all of its vertices.

\begin{theo} \label{almostpacking}
Suppose that $1/n \ll \alpha \ll 1/k$, and that $J$ is a $k$-system on $n$ vertices with
$$\delta(J) \geq \left(n, \left(\frac{k-1}{k} - \alpha\right)n,
\left(\frac{k-2}{k} - \alpha\right) n, \dots, \left(\frac{1}{k} - \alpha\right) n\right).$$
Then $J_k$ contains a matching which covers all but at most $9k^2\alpha n$ vertices of $J$.
\end{theo}

Next we describe a family of examples showing that Theorem~\ref{almostpacking} is best possible up to the constant $9k^2$ (which we do not attempt to optimise), in the sense that there may not be any almost perfect matching if {\em any} coordinate of the degree sequence is substantially reduced.

\begin{construct}\label{spacebar} (Space barriers)
Suppose $V$ is a set of $n$ vertices, $j \in [k-1]$ and $S \sub V$.  Let $J = J(S,j)$ be the $k$-complex in which $J_i$ (for $0 \le i \le k$) consists of all $i$-sets in $V$ that contain at most $j$ vertices of $S$.  
\end{construct}

\begin{figure}[t] 
\centering
\psfrag{A}{$V\sm S$}
\psfrag{B}{$S$}
\includegraphics[width=5cm]{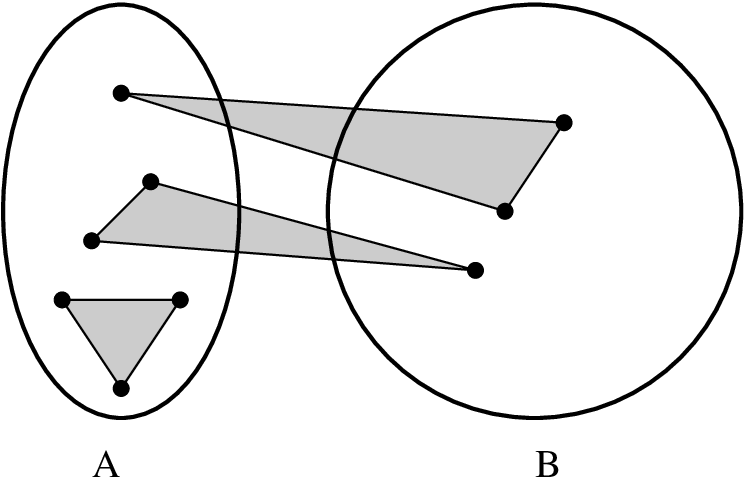} \hspace{1.5cm}
\includegraphics[width=5cm]{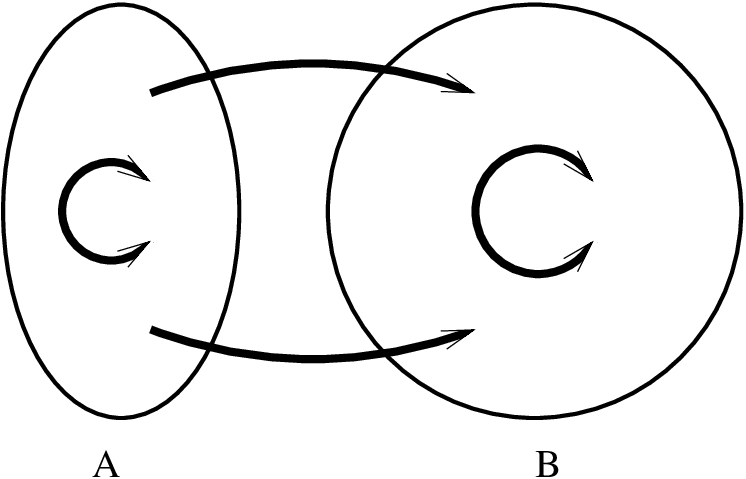} 
\caption{The left hand diagram shows an example of a space barrier $H$: the edges are all $3$-tuples which have at most two vertices in $S$. If $|S| > \frac{2}{3}|V(H)|$ then $H$ has no perfect matching. If instead $|S| = \frac{2}{3}|V(H)|$ then a perfect matching $M$ in $H$ must consist only of edges with exactly two vertices in $S$, whereupon for any $\ell$ the transferral digraph $D_\ell(H, M)$ is as shown in the right hand diagram; it contains all edges except those directed from $S$ to $V \sm S$.}
\label{fig:spacebar} 
\end{figure}

Figure~\ref{fig:spacebar} illustrates this construction in the case $k=3$ and $j=2$, along with the corresponding \emph{transferral digraph}, which will be defined formally later. Observe that for any $i \leq j$ the $i$-graph $J_i = \binom{V}{i}$ is complete, so we have $\delta_i(J) = n-i$ for $0 \le i \le j-1$ and $\delta_i(J)=n-|S|-(i-j)$ for $i \geq j$. Choosing $|S| = \bfl{\brac{\frac{j}{k}+\alpha}n}-k$, we see that the minimum degree condition of Theorem~\ref{almostpacking} is satisfied. However, every edge of $J_k$ has at least $k-j$ vertices in $V \sm S$, so the maximum matching in $J_k$ has size $\bfl{\frac{|V \sm S|}{k-j}}$, which leaves at least $\alpha n - k$ uncovered vertices. 

In these examples we have a space barrier for matchings, which is robust in the sense that the number of uncovered vertices varies smoothly with the minimum degree sequence. Our preliminary stability result is that space barriers are the only obstructions to finding a matching that covers all but a constant number of vertices. Note that we cannot avoid leaving a constant number of uncovered vertices because of divisibility barriers (discussed below). Our structures permit some small imperfections, defined as follows. Suppose $G$ and $H$ are $k$-graphs on the same set $V$ of $n$ vertices and $0 < \beta < 1$. We say that $G$ is {\em $\beta$-contained} in $H$ if all but at most $\beta n^k$ edges of $G$ are edges of $H$.

\begin{theo} \label{prestability}
Suppose that $1/n \ll 1/\ell \ll \alpha \ll \beta \ll 1/k$.
Let $J$ be a $k$-complex on $n$ vertices such that
$$\delta(J) \geq \left(n, \left(\frac{k-1}{k} - \alpha\right)n, \left(\frac{k-2}{k} - \alpha\right)
  n, \dots, \left(\frac{1}{k} - \alpha\right) n\right).$$
Then $J$ has at least one of the following properties:
\begin{description}
\item[1 (Matching)] $J_k$ contains a matching that covers all but at most $\ell$ vertices.
\item[2 (Space barrier)] $J_k$ is $\beta$-contained in $J(S,j)_k$ for some $j \in [k-1]$ and $S \sub V(J)$ with $|S| = \lfloor jn/k \rfloor$.
\end{description}
\end{theo}

\section{Partite systems.} \label{sec:partite}
We will also require `partite' analogues of our hypergraph matching theorems.
A \emph{partition} $\Part$ of a set $V$ is a sequence of disjoint sets $V_1, \dots, V_k$ whose union is $V$; we refer to the
sets $V_i$ as the \emph{parts} of $\Part$, and write $U \in \Part$ to mean that $U$ is a part of $\Part$. Note that we consider the partition $\Part$ to describe not just the contents
of each part but also the order of the parts. We say that the partition $\Part$ is
\emph{balanced} if each part has the same size.

Let $H$ be a hypergraph, and let $\Part$ be a partition of $V(H)$. Then we say a set $S \subseteq V(H)$ is \emph{$\Part$-partite} if if has at most one vertex in any part of $\Part$, and that $H$ is \emph{$\Part$-partite} if
every edge of $H$ is $\Part$-partite. We say that $H$ is \emph{$r$-partite} if there
exists some partition $\Part$ of $V(H)$ into $r$ parts such that $H$ is
$\Part$-partite. For $r$-partite $k$-systems we introduce the following alternative notion of degree.
Let $V$ be a set of vertices, let $\Part$ be a partition of $V$ into $r$ parts $V_1, \dots, V_r$,
and let $J$ be a $\Part$-partite $k$-system on $V$. For each $0 \leq j \leq k-1$ we define the
\emph{partite minimum $j$-degree} $\delta^*_j(J)$ as the largest $m$ such that any $j$-edge $e$ has
at least $m$ extensions to a $(j+1)$-edge in any part not used by $e$, i.e.\
$$\delta^*_j(J) := \min_{e \in J_j} \min_{i : e \cap V_i = \es} |\{v \in V_i : e \cup \{v\} \in J\}|.$$
The \emph{partite degree sequence} is $\delta^*(J) = (\delta_0^*(J), \dots, \delta_{k-1}^*(J))$.
Note that we suppress the dependence on $\Part$ in our notation: this will be clear from the context.

Our next theorem is an analogue of Theorem~\ref{almostpacking} for $r$-partite $k$-systems. Here we may add an additional condition on our matching, for which we need the following definition. Suppose that $\Part$ is a partition of a set of vertices $V$ into vertex classes $V_1, \dots, V_r$, and $H$ is a $\Part$-partite $k$-graph on $V$. Then the \emph{index set} of an edge $e \in H$ is $\{i : |e \cap V_i| = 1\} \in \binom{[r]}{k}$. For a matching $M$ in $H$ and a set $A \in \binom{[r]}{k}$ let $N_A(M)$ denote the number of edges $e \in M$ whose index set is $A$. We say that $M$ is \emph{balanced} if $N_A(M)$ is constant over all $A \in \binom{[r]}{k}$, that is, each index set is represented by equally many edges. The following theorem shows that we may insist that the matching obtained is balanced.

\begin{theo} \label{almostpackingpartite}
Suppose that $1/n \ll \alpha \ll 1/r \le 1/k$, and that $J$ is a $r$-partite $k$-system on vertex classes each of $n$ vertices with
$$\delta^*(J) \geq \left(n, \left(\frac{k-1}{k} - \alpha\right)n, \left(\frac{k-2}{k} - \alpha\right) n, \dots, \left(\frac{1}{k} - \alpha\right) n\right).$$
Then $J_k$ contains a balanced matching which covers all but at most $9k^2r\alpha n$ vertices of $J$.
\end{theo}

To see that this is best possible in the same sense as for Theorem~\ref{almostpacking} we use the natural partite version of Construction~\ref{spacebar}, which shows that there may not be any almost perfect matching if any coordinate of the partite degree sequence is substantially reduced.

\begin{construct}\label{spacebarpartite} (Partite space barriers)
Suppose $\Part$ partitions a set $V$ into $r$ parts $V_1,\dots,V_r$ of size $n$. Let $j \in [k-1]$, $S \sub V$ and $J = J_r(S,j)$ be the $k$-complex in which $J_i$ (for $0 \le i \le k$) consists of all $\Part$-partite $i$-sets in $V$ that contain at most $j$ vertices of $S$. \end{construct}

We choose $S$ to have $s = \bfl{\brac{\frac{j}{k}+\alpha}n}$ vertices in each part. Then $\delta^*_i(J) = n$ for $0 \le i \le j-1$ and $\delta^*_i(J)=n-s$ for $j \le i \le k-1$, so the minimum partite degree condition of Theorem~\ref{almostpackingpartite} is satisfied. However, the maximum matching in $J_k$ has size $\bfl{\frac{|V \sm S|}{k-j}}$, which leaves at least $r(\alpha n - k)$ uncovered vertices.

We also have the following stability result analogous to Theorem~\ref{prestability}.

\begin{theo} \label{prestabilitypartite}
Suppose that $1/n \ll 1/\ell \ll \alpha \ll \beta \ll 1/r \le 1/k$. Let~$\Part$ partition a set $V$ into parts $V_1, \dots, V_r$ each of size~$n$, and let~$J$ be a $\Part$-partite $k$-complex on $V$ with $$\delta^*(J) \geq \left(n, \left( \frac{k-1}{k} - \alpha\right)n, \left(\frac{k-2}{k} - \alpha\right) n, \dots, \left(\frac{1}{k} - \alpha\right) n \right).$$
Then $J$ has at least one of the following properties:
\begin{description}
\item[1 (Matching)] $J_k$ contains a matching that covers all but at most $\ell$ vertices.
\item[2 (Space barrier)] $J_k$ is $\beta$-contained in $J_r(S,j)_k$ for some $j \in [k-1]$ and $S \sub V$ with $\lfloor jn/k \rfloor$ vertices in each $V_i$, $i \in [r]$.
\end{description}
\end{theo}

\section{Lattice-based constructions.} \label{sec:lattice}
Having described the general form of space barriers, we now turn our attention to divisibility barriers. For this we need the notion of \emph{index vectors}, which will play a substantial role in this paper.
Let $V$ be a set of vertices, and let $\Part$ be a partition of $V$ into $d$ parts $V_1, \dots, V_d$.
Then for any $S \sub V$, the \emph{index vector of $S$ with respect to $\Part$}
is the vector 
$$\ib_\Part(S) := (|S \cap V_1|, \dots, |S \cap V_d|) \in \Z^d.$$
When $\Part$ is clear from the context, we write simply $\ib(S)$ for $\ib_\Part(S)$.
So $\ib(S)$ records how many vertices of $S$ are in each part of $\Part$.
Note that $\ib(S)$ is well-defined as we consider the partition $\Part$
to define the order of its parts.

\begin{construct} \label{divbar} (Divisibility barriers)
Suppose $L$ is a lattice in $\Z^d$ (i.e.\ an additive subgroup) with $\ib(V) \notin L$,
fix any $k \ge 2$, and let $G$ be the $k$-graph on $V$ whose edges are all $k$-tuples $e$ with $\ib(e) \in L$.
\end{construct}

For any matching $M$ in $G$ with vertex set $S = \bigcup_{e \in M} e$ we have $\ib(S)  = \sum_{e \in M} \ib(e) \in L$.
Since we assumed that $\ib(V) \notin L$ it follows that $G$ does not have a perfect matching.

\begin{figure}[t] 
\centering
\psfrag{4}{$V_1$}
\psfrag{2}{$V_2$}
\psfrag{3}{$V_3$}
\mbox{
\includegraphics[width=5cm]{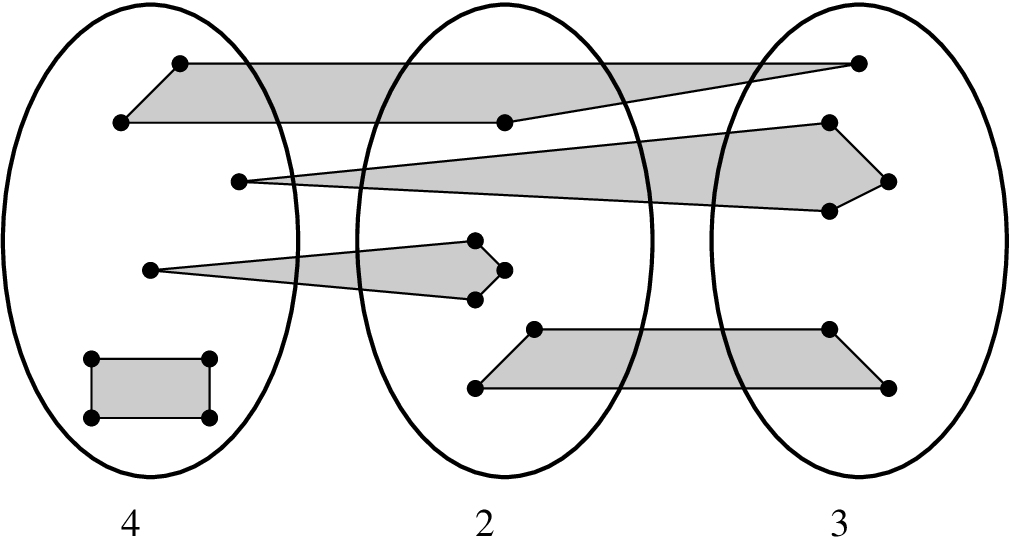} \hspace{1cm}
\includegraphics[width=5cm]{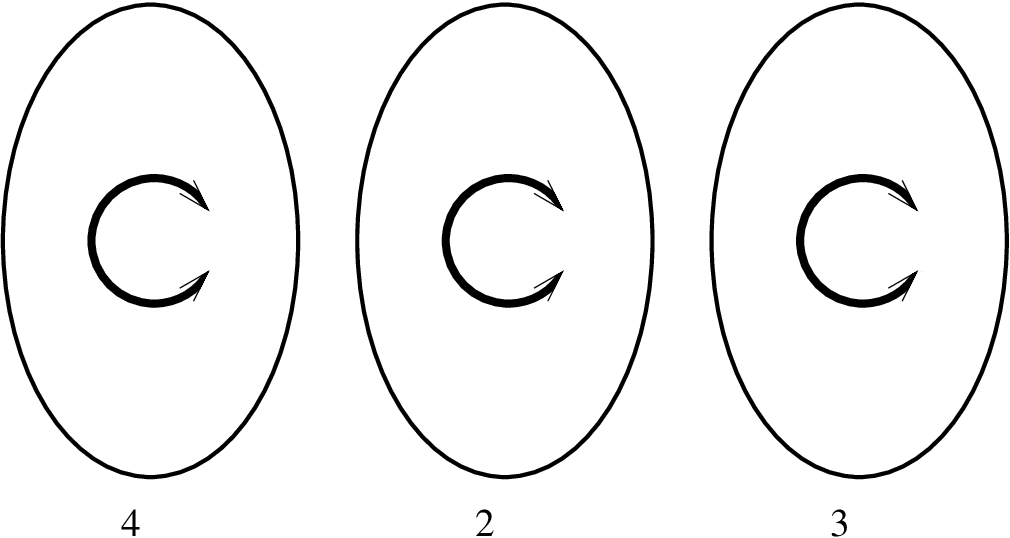}}
\caption{The left hand diagram shows an example of a divisibility barrier $H$: the edges are all $4$-tuples $e$ with $|e \cap V_2| = |e \cap V_3| \mod 3$. So if $|V_2| \neq |V_3| \mod 3$ then there is no perfect matching in $H$. If instead $|V_2| = |V_3|$ and $H$ contains a perfect matching, then for any $\ell$ the transferral digraph $D_\ell(H, M)$ is as shown in the right hand diagram; it consists of three disjoint cliques induced by $V_1, V_2$ and $V_3$.}
\label{fig:divbar} 
\end{figure}

We will now give some concrete illustrations of this construction.
\begin{enumerate}
\item Suppose $d=2$ and $L = \bgen{(-2,2),(0,1)}$. Note that $(x,y) \in L$ precisely when $x$ is even.
Then by definition, $|V_1|$ is odd, and the edges of $G$ are all $k$-tuples $e \sub V$ such that $|e \cap V_1|$
is even. If $|V| = n$ and $|V_1| \sim n/2$, then we recover the extremal example mentioned earlier for the minimum degree perfect matching problem.
The generated $k$-complex $G^\le$ has $\delta_i(G^\le) \sim n$ for $0 \le i \le k-2$ and $\delta_{k-1}(G^\le) \sim n/2$.
\item Suppose $d=3$ and $L = \bgen{(-2,1,1),(1,-2,1),(1,0,0)}$. Note that $(x,y,z) \in L$ precisely
when $y = z$ mod $3$ (This construction is illustrated in Figure~\ref{fig:divbar}, again with the accompanying transferral digraph, which will be defined formally later.). Thus $|V_2| \ne |V_3|$ mod $3$ and the edges of $G$ are  all $k$-tuples $e \sub V$ such that
$|e \cap V_2| = |e \cap V_3|$ mod $3$. Note that for any $(k-1)$-tuple $e' \sub V$ there is a unique $j \in [3]$
such that such that $\ib(e')+\unit_j \in L$, where $\unit_j$ is the $j$th standard basis vector.
(We have $j=1$ if $\ib(e')_2=\ib(e')_3$ mod $3$, $j=2$ if  $\ib(e')_2=\ib(e')_3-1$ mod $3$,
or $j=3$ if $\ib(e')_2=\ib(e')_3+1$ mod $3$.) If $|V| = n$ and $|V_1|, |V_2|, |V_3| \sim n/3$ then $\delta_i(G^\le) \sim n$ for $0 \le i \le k-2$
and $\delta_{k-1}(G^\le) \sim n/3$. Mycroft~\cite{M} showed that this construction is asymptotically extremal for a range of hypergraph packing problems. For example, he showed that any $4$-graph $G$ on $n$ vertices with $7 \mid n$ and $\delta(G) \geq n/3 + o(n)$ contains a perfect $K^4_{4,1,1,1}$-packing, where $K^4_{4,1,1,1}$ denotes the complete $4$-partite $4$-graph with vertex classes of size $4$, $1$, $1$ and $1$ respectively; this construction demonstrates that this bound is best possible up to the $o(n)$ error term. 
\end{enumerate}

For simplicity we only gave approximate formulae for the degree sequences in these examples, but it is not hard to calculate them exactly. Note that the index vectors of edges in a $k$-graph belong to the hyperplane $\{\xb \in \R^d: \sum_{i \in [d]} x_i = k\}$. The intersection of this hyperplane with the lattice $L$ is a coset of the intersection of $L$ with the hyperplane $\Pi^d = \{\xb \in \R^d: \sum_{i \in [d]} x_i = 0\}$. Thus it is most convenient to choose a basis of $L$ that includes a basis of $L \cap \Pi^d$, as in the above two examples. In order for the construction to exist it must be possible to satisfy $\ib(V) \notin L$, so $L$ should not contain the lattice $M^d_k = \{\xb \in \Z^d: k \mid \sum_{i\in[d]} x_i\}$. We say that $L$ is {\em complete} if $M^d_k \sub L$, otherwise $L$ is {\em incomplete}. Assuming that $L$ contains some $\vb$ with $\sum v_i=k$, an equivalent formulation of completeness is that $L \cap \Pi^d = \Z^d \cap \Pi^d$, since this is equivalent to
$L \cap (\Pi^d+\vb) = \Z^d \cap (\Pi^d+\vb) = \{\xb \in \Z^d: \sum_{i \in [d]} x_i = k\}$, and so to $M^d_k \sub L$.

There is a natural notion of minimality for lattice-based constructions. We say that $L$ is {\em transferral-free} if it does not contain any difference of standard basis vectors, i.e.\ a vector of the form $\unit_i-\unit_j$ with $i \ne j$. If $L$ is not transferral-free, then the construction for $L$ can be reformulated using a lattice of smaller dimension. Without loss of generality we may consider the case that $\unit_{d-1}-\unit_d \in L$. Then we replace $\Part$ by the partition of $V$ into $d-1$ parts obtained by combining $V_{d-1}$ and $V_d$, and replace $L$ by the lattice $L'$ consisting of all $\xb \in \Z^{d-1}$ such that $(x_1,\dots,x_{d-1},0) \in L$.  We remark that a transferral-free lattice in $\Z^d$ has index at least $d$ as a subgroup of $\Z^d$.  To see this, note that the cosets $L + (\unit_1-\unit_i)$, $1 \le i \le d$ must be distinct, otherwise there is some $\xb \in Z^d$ and $i,j \in [d]$ for which $\xb = \vb + \unit_1-\unit_i = \vb' + \unit_1 - \unit_j$ with $\vb,\vb' \in L$, which gives $\unit_i - \unit_j = \vb-\vb' \in L$, contradicting the fact that $L$ is transferral-free.  

\section{Perfect matchings.} \label{sec:pm}

For our main theorems we specialise to the setting of simplicial complexes, which will simplify some arguments and is sufficiently general for our applications.  Under the same minimum degree assumption that is tight for almost perfect matchings, we refine the results above, by showing that either we have a perfect matching, or we have a structural description for the complex: it either has a space barrier or a divisibility barrier. It is remarkable that this rigidity emerges from the purely combinatorial degree sequence assumption, and that it exhibits these two very different phenomena, one tied to convex geometry and the other to integer lattices.

Our structures permit some small imperfections, defined as follows. Recall that we say that a $k$-graph $G$ is {\em $\beta$-contained} in a $k$-graph $H$ if all but at most $\beta n^k$ edges of $G$ are edges of $H$. Also, given a partition $\Part$ of $V$ into $d$ parts, we define the {\em $\mu$-robust edge lattice} $L^\mu_\Part(G) \sub \Z^d$ to be the lattice generated by all vectors $\vb \in \Z^d$ such that there are at least $\mu n^k$ edges $e \in G$ with $\ib_\Part(e) = \vb$. Recall that we call a lattice $L \sub \Z^d$ complete if $L \cap \Pi^d = \Z^d \cap \Pi^d$, where $\Pi^d = \{\xb \in \R^d: \sum_{i \in [d]} x_i = 0\}$. Recall also that the space barrier constructions $J(S,j)$ were defined in Section~\ref{sec:almostpm}. Now we can state our first main theorem.

\begin{theo} \label{newmain}
Suppose that $1/n \ll \alpha \ll \mu, \beta \ll 1/k$ and that $k\mid n$.
Let $J$ be a $k$-complex on $n$ vertices such that
$$\delta(J) \geq \left(n, \left(\frac{k-1}{k} - \alpha\right)n, \left(\frac{k-2}{k} - \alpha\right)
  n, \dots, \left(\frac{1}{k} - \alpha\right) n\right).$$
Then $J$ has at least one of the following properties:
\begin{description}
\item[1 (Matching)] $J_k$ contains a perfect matching.
\item[2 (Space barrier)] $J_k$ is $\beta$-contained in $J(S,j)_k$ for some $j \in [k-1]$
and $S \sub V(J)$ with $|S| = jn/k$.
\item[3 (Divisibility barrier)] There is some partition $\Part$ of $V(J)$ into $d \le k$
parts of size at least $\delta_{k-1}(J) - \mu n$ such that $L^\mu_\Part(J_k)$ is incomplete and transferral-free.
\end{description}
\end{theo}

Our second main theorem is a partite version of the previous result.  Recall that the space barrier constructions in this setting are described in Construction~\ref{spacebarpartite}.  We also need to account for the original partition when classifying edge lattices as follows. Suppose $\Part$ is a partition of $V$ into $d$ parts $(V_1, \dots, V_d)$ that refines a partition $\Part'$ of $V$. Let $L_{\Part\Part'} \sub \Z^d \cap \Pi^d$ be the lattice generated by all differences of basis vectors $\unit_i-\unit_j$ for which $V_i,V_j$ are contained in the same part of $\Part'$. We say that a lattice $L \sub \Z^d$ is {\em complete with respect to $\Part'$} if $L_{\Part\Part'} \sub L \cap \Pi^d$, otherwise we say that $L$ is  {\em incomplete with respect to $\Part'$}. 

Similarly to Theorem~\ref{newmain}, this theorem tells us that an $r$-partite $k$-complex satisfying the given minimum degree condition either contains a perfect matching or has a space barrier or divisibility barrier. However, in applications (for example the multipartite Hajnal-Szemer\'edi theorem proved in~\cite{KM}) will may need to know that the perfect matching obtained has roughly the same number of edges of each index. Recall that the index set of an edge $e \in H$ is $\{i : |e \cap V_i| = 1\} \in \binom{[r]}{k}$, that for a matching $M$ in $H$ we write $N_A(M)$ to denote the number of edges $e \in M$ whose index set is $A$, and that $M$ is balanced if $N_A(M)$ is constant over all $A \in \binom{[r]}{k}$. Unfortunately, Construction~\ref{notbalanced} will show that we cannot insist on a balanced matching in our partite analogue of Theorem~\ref{newmain}. Instead we require a weaker property: we say that $M$ is \emph{$\gamma$-balanced} if $N_A(M) \geq (1 - \gamma) N_B(M)$ for any $A, B \in \binom{[r]}{k}$, meaning that $M$ is close to being balanced.

\begin{theo} \label{newmainpartite}
Suppose that $1/n \ll \gamma, \alpha \ll \mu, \beta \ll 1/r \le 1/k$. Let~$\Part'$ partition a set $V$ into parts $V_1, \dots, V_r$ each of size~$n$, where $k \mid rn$. Suppose that $J$ is a $\Part'$-partite $k$-complex with $$\delta^*(J) \geq \left(n, \left(\frac{k-1}{k} - \alpha\right)n, \left(\frac{k-2}{k} - \alpha\right) n, \dots, \left(\frac{1}{k} - \alpha\right) n\right).$$
Then $J$ has at least one of the following properties:
\begin{description}
\item[1 (Matching)] $J_k$ contains a $\gamma$-balanced perfect matching.
\item[2 (Space barrier)] $J_k$ is $\beta$-contained in $J_r(S,j)_k$ for some $j \in [k-1]$ and $S \sub V$ with $\lfloor jn/k \rfloor$ vertices in each $V_i$, $i \in [r]$.
\item[3 (Divisibility barrier)] There is some partition $\Part$ of $V(J)$ into $d \le kr$
parts of size at least  $\delta^*_{k-1}(J) - \mu n$ such that $\Part$ refines $\Part'$
and $L^\mu_\Part(J_k)$ is incomplete with respect to $\Part'$ and transferral-free.
\end{description}
\end{theo} 

A necessary condition for the existence of a balanced perfect matching in an $r$-partite $k$-graph whose vertex classes have size $n$ is that $rn/k$, the number of edges, is divisible by $\binom{r}{k}$, the number of possible index sets of edges. However, even under this additional assumption we cannot replace `$\gamma$-balanced' with `balanced' in option~1 of this theorem; some small imbalance may be inevitable, as shown by the following construction.

\begin{construct} \label{notbalanced}
Choose integers $r \geq 3$ and $n$ so that $2n/(r-1)$ and $rn$ are integers of different parity (e.g. $r=5$ and $n \equiv 2 \mod 4$). 
Let $G$ be an $r$-partite graph whose vertex classes each have size two, say $\{x_i, y_i\}$ for $i \in [r]$, and whose edges are $x_1y_2$, $x_2y_1$, and $x_ix_j$ and $y_iy_j$ for any pair $\{i,j\}$ other than $\{1, 2\}$. Form the `blow-up' $G^*$ by replacing every vertex of $G$ with $n$ vertices (for some even integer $n$), and adding edges between any pair of vertices which replace adjacent vertices in $G$.
\end{construct}

Observe that the graph $G^*$ constructed in Construction~\ref{notbalanced} is an $r$-partite graph whose vertex classes each have size $2n$ and which satisfies $\delta^*(G^*) = n$. So the 2-complex $J = (G^*)^\leq$ satisfies the conditions of Theorem~\ref{newmainpartite}. Furthermore, it is not hard to check that $J$ does not satisfy options 2 or 3 of Theorem~\ref{newmainpartite} for small $\mu, \beta$. So $J_2 = G^*$ must contain a perfect matching which is close to being balanced, and indeed it is not hard to verify this. However, there is no balanced perfect matching in $G^*$. Indeed, let $Y$ be the set of vertices of $G^*$ which replaced one of the vertices $y_i$ of $G$, so $|Y| = rn$. In a balanced perfect matching $M$ each index set would be represented by $rn/\binom{r}{2} = \frac{2n}{r-1}$ edges of $M$. This means that the number of vertices of $Y$ covered by edges of $M$ of index $\{1, 2\}$ is $\frac{2n}{r-1}$, and so $|Y| - \frac{2n}{r-1} = rn -\frac{2n}{r-1}$ must be even, a contradiction.

Theorems~\ref{newmain} and~\ref{newmainpartite} will each be deduced from a more general statement, Theorem~\ref{perfectmatchingF}, in Chapter~\ref{sec:proofs}. However, given a slightly stronger degree sequence, we can rule out the possibility of a space barrier in each of these theorems. This will frequently be the case in applications, so for ease of use we now give the corresponding theorems.

\begin{theo} \label{weakmain}
Suppose that $1/n \ll \mu \ll \alpha, 1/k$ and that $k \mid n$.
Let $J$ be a $k$-complex on $n$ vertices such that
$$\delta(J) \geq \left(n, \left(\frac{k-1}{k} + \alpha\right)n, \left(\frac{k-2}{k} + \alpha \right)
  n, \dots, \left(\frac{1}{k} + \alpha\right) n\right).$$
Then $J$ has at least one of the following properties:
\begin{description}
\item[1 (Matching)] $J_k$ contains a perfect matching.
\item[2 (Divisibility barrier)] There is some partition $\Part$ of $V(J)$ into $d \le k$
parts of size at least $\delta_{k-1}(J) - \mu n$ such that $L^\mu_\Part(J_k)$ is incomplete and transferral-free.
\end{description}
\end{theo}

\proof
Introduce a new constant $\beta$ with $1/n \ll \beta \ll \alpha, 1/k$. Then it suffices to show that option 2 of Theorem~\ref{newmain} is impossible. So suppose for a contradiction that $J_k$ is $\beta$-contained in $J(S,j)_k$ for some $j \in [k-1]$
and $S \sub V(J)$ with $|S| = jn/k$. We now form an ordered $k$-tuple $(x_1, \dots, x_k)$ of vertices of $J$ such that $\{x_1, \dots x_s\} \in J$ for any $s \in [k]$ and $x_1, \dots, x_{j+1} \in S$. Indeed, the minimum degree sequence ensures that when greedily choosing the vertices one at a time there are at least $n-|S| \geq n/k$ choices for the first vertex, at least $\delta_i(J) - (n-|S|) \geq (j-i)n/k+\alpha n$ choices for the $(i+1)$st vertex for any $i \leq j$, and at least $\delta_i(J) \geq (k-i)n/k + \alpha n$ for the $(i+1)$st vertex for any $i > j$. In total this gives at least $\alpha (n/k)^k$ ordered $k$-tuples, each of which is an edge of $J_k$ with at least $j+1$ vertices in $S$, so is not contained in $J(S, j)_k$. This is a contradiction for $k!\beta < \alpha/k^k$. 
\endproof

\begin{theo} \label{weakmainpartite}
Suppose that $1/n \ll \gamma, \mu \ll \alpha, 1/r, 1/k$ and $r \geq k$. Let~$\Part'$ partition a set $V$ into parts $V_1, \dots, V_r$ each of size~$n$, where $k \mid rn$. Suppose that $J$ is a $\Part'$-partite $k$-complex with $$\delta^*(J) \geq \left(n, \left(\frac{k-1}{k} + \alpha\right)n, \left(\frac{k-2}{k} + \alpha\right) n, \dots, \left(\frac{1}{k} + \alpha\right) n\right).$$
Then $J$ has at least one of the following properties:
\begin{description}
\item[1 (Matching)] $J_k$ contains a $\gamma$-balanced perfect matching.
\item[2 (Divisibility barrier)] There is some partition $\Part$ of $V(J)$ into $d \le kr$
parts of size at least  $\delta^*_{k-1}(J) - \mu n$ such that $\Part$ refines $\Part'$
and $L^\mu_\Part(J_k)$ is incomplete with respect to $\Part'$ and transferral-free.
\end{description}
\end{theo}

This theorem follows from Theorem~\ref{newmainpartite} exactly as Theorem~\ref{weakmain} followed from Theorem~\ref{newmain}; we omit the details.

\section{Further results.}
Theorems~\ref{newmain} and~\ref{newmainpartite} can be applied to a variety of matching and packing problems in graphs and hypergraphs. Indeed, in this chapter we give a short deduction of Theorem~\ref{partitehajnalszem} from Theorem~\ref{weakmainpartite} (which was a consequence of Theorem~\ref{newmainpartite}), whilst in Chapter~\ref{sec:tetra} we use Theorem~\ref{newmain} to prove Theorem~\ref{TETRAPACK}. We can also recover and find new variants of existing results. For example, consider the result of R\"odl, Ruci\'nski and Szemer\'edi~\cite{RRS} on the minimum degree threshold for a perfect matching in a $k$-graph. Their proof proceeds by a stability argument, giving a direct argument when the $k$-graph is close to an extremal configuration, and otherwise showing that even a slightly lower minimum degree is sufficient for a perfect matching. Our results give a new proof of stability under a much weaker degree assumption. In  the following result, we only assume that the minimum degree of $G$ is a bit more than $n/3$, and show that if there is no perfect matching then $G$ is almost contained in an extremal example. As described earlier, with a bit more work (and a stronger degree assumption) we can show that $G$ is contained in an extremal example (see Theorem~\ref{pmstrongstab}). This requires some technical preliminaries, so we postpone the proof for now.

\medskip 

\begin{theo} \label{pmstab}
Suppose $1/n \ll b,c \ll 1/k$, $k \ge 3$, $k$ divides $n$ and $G$ is a $k$-graph on $n$ vertices with $\delta(G) \ge (1/3+c)n$
and no perfect matching. Then there is a partition of $V(G)$ into parts $V_1,V_2$ of size at least
$\delta(G)$ and $a \in \{0,1\}$ so that all but at most $bn^k$ edges $e$ of $G$ have $|e \cap V_1|=a$ mod $2$.
\end{theo}

\proof Introduce a constant $\mu$ with $1/n \ll \mu \ll b,c$. Let $J$ be the (clique) $k$-complex with $J_k=G$ and $J_i$ complete for $i<k$. Suppose that $J_k$ does not have a perfect matching. Then option 2 must hold in Theorem~\ref{weakmain}, that is, there is some partition $\Part$ of $V(J)$ into parts of size at least $\delta(G) - \mu n$ such that $L^\mu_\Part(J_k)$ is incomplete. Since $\delta(G) - \mu n > n/3$ there must be $2$ such parts, and $L^\mu_\Part(J_k) \cap \Pi^2$ is generated by $(-t,t)$ for some $t \ge 2$. We cannot have $t>2$, as then neither $(k-2, 2)$ nor $(k-1, 1)$ lie in $L^\mu_\Part(J_k)$, and so there are at most $2 \mu n^k$ edges of $J_k = G$ with one of these two indices. Since there are at least $(n/3)^{k-1}/(k-1)!$ $(k-1)$-tuples of index $(k-2, 1)$, by averaging some such $(k-1)$-tuple must be contained in at most $3^kk!\mu n^k$ edges of $G$, contradicting the minimum degree assumption. So we must have $t=2$, and the result follows.  \endproof

We now present the deduction of Theorem~\ref{partitehajnalszem} from Theorem~\ref{weakmainpartite}.

\medskip \noindent {\bf Proof of Theorem~\ref{partitehajnalszem}.}
Introduce a constant $\mu$ with $1/n \ll \mu \ll c$. Let $J$ be the clique $k$-complex of $G$. Then $\delta_i(J) \ge (k-i)n/k + icn$ for $0 \le i \le k-1$.  Thus we can apply Theorem~\ref{weakmainpartite}. Suppose for a contradiction that option 2 of this theorem holds. Then there is some partition $\Part$ of $V(J)$ into parts of size at least $\delta^*_{k-1}(J) - \mu n$ such that $\Part$ refines the partition $\Part'$ of $V$ into $V_1,\cdots,V_k$ and $L^\mu_\Part(J_k)$ is incomplete with respect to $\Part'$.  We may assume that $L^\mu_\Part(J_k)$ is transferral-free; recall that this means that it does not contain any difference of standard basis vectors $\unit_i-\unit_j$ with $i \ne j$.  Consider a part of $\Part'$ that is refined in $\Part$, without loss of generality it is $V_1$, and vertices $x_1,x'_1\in V_1$ in different parts $U_1$, $U'_1$ of $\Part$.  We can greedily construct many sequences $x_i \in V_i$, $2 \le i \le k$ such that $x_1 x_2 \dots x_k$ and $x'_1 x_2 \dots x_k$ are both $k$-cliques: since $\delta^*(G) \ge (k-1)n/k + cn$ there are at least $(k-i)n/k + icn$ choices for $x_i$, so at least $c(n/k)^{k-1}$ such sequences.  We can repeat this for all choices of $x_1 \in U_1$ and $x'_1 \in U_1'$; there are at least $\delta^*_{k-1}(J) - \mu n > n/k$ choices for each.  Since $c \gg \mu$, there is some choice of parts $U_i \sub V_i$, $2 \le i \le k$ of $\Part$ such that we obtain at least $\mu n^k$ cliques intersecting $U_i$, $i \in [k]$, and at least $\mu n^k$ cliques intersecting $U'_1$ and $U_i$, $2 \le i \le k$.  However, this contradicts the fact that $L^\mu_\Part(J_k)$ is transferral-free. Thus Theorem~\ref{weakmainpartite} implies that $J_k$ has a perfect matching, as required. \endproof

\section{Outline of the proofs.}\label{sec:outline}

The ideas of our arguments can be roughly organised into the following three groups: regularity, transferrals, applications. Most of the new ideas in this paper pertain to transferrals, but to set the scene we start with regularity. The Szemer\'edi Regularity Lemma \cite{Sz} has long been a powerful tool in graph theory. In combination with the blowup lemma of  Koml\'os, S\'ark\"ozy and Szemer\'edi \cite{KSS} it has seen many applications to embeddings of spanning subgraphs (see \cite{KO2}). Recent developments in hypergraph regularity theory have opened the way towards obtaining analogous results for hypergraphs: the decomposition theory (among other things) was developed independently by Gowers \cite{G1} and by R\"odl et al. \cite{FR,NR,RSk,RS}, and the blowup lemma by Keevash \cite{K}. Roughly speaking, the decomposition theory allows us to approximate a $k$-system $J$ on $n$ vertices by a `reduced' $k$-system $R$ on $m$ vertices, where $m$ depends on the accuracy of the approximation, but is independent of $n$. The vertices of $R$ correspond to the parts in some partition of $V(J)$ into `clusters' of equal size, after moving a small number of vertices to an exceptional set $V_0$. The edges of $R$ correspond to groups of clusters for which the restriction of the appropriate level of $J$ is well-approximated by a `dense regular' hypergraph. Furthermore, $R$ inherits from $J$ approximately the same proportional minimum degree sequence. As mentioned earlier, this part of the machinery allows us to reduce the almost perfect matching problem to finding a fractional solution. If $J$ has an extra $o(n)$ in its minimum degree sequence then it is a relatively simple problem in convex geometry to show that $R_k$ has a fractional perfect matching $M$; moreover, in the absence of a space barrier we can find $M$ even without this extra $o(n)$. Then we partition the clusters of $J$ in proportion to the weights of edges in $M$, so that each non-zero edge weight in $M$ is associated to a dense regular $k$-partite $k$-graph with parts of equal size (adding a small number of vertices to $V_0$). It is then straightforward to find almost perfect matchings in each of these $k$-partite $k$-graphs, which together constitute an almost perfect matching in $J_k$.

To find perfect matchings, we start by taking a regularity decomposition and applying the almost perfect matching result in the reduced system. We remove any uncovered clusters, adding their vertices to the exceptional set $V_0$, so that we have a new reduced system $R$ with a perfect matching $M$. We also transfer a small number of vertices from each cluster to $V_0$ so that the edges of $M$ now correspond to dense regular $k$-partite $k$-graphs that have perfect matchings (rather than almost perfect matchings). In fact, these $k$-graphs have the stronger property of `robust universality': even after deleting a small number of additional vertices, we can embed any bounded degree subhypergraph. (Up to some technicalities not mentioned here, this is the form in which we apply the hypergraph blowup lemma.) Next we greedily find a partial matching that covers $V_0$, and also removes some vertices from the clusters covered by $M$, where we take care not to destroy robust universality. Now the edges of $M$ correspond to $k$-partite $k$-graphs that are robustly universal, but have slightly differing part sizes. To find a perfect matching, we will remove another partial matching that balances the part sizes and does not destroy robust universality. Then the edges of $M$ will now correspond to $k$-partite $k$-graphs that have perfect matchings, and together with the two partial matchings these constitute a perfect matching in $J_k$.

\begin{figure}[t] 
\centering 
\includegraphics[width=10cm]{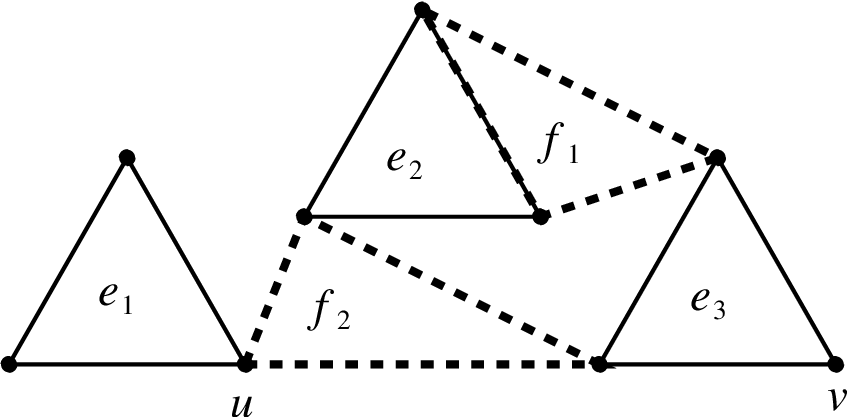} 
\caption{An example of a $1$-fold $(u, v)$-transferral in a matched $3$-graph $(J, M)$. Here $e_1$, $e_2$ and $e_3$ are edges of the matching $M$ in $J$, whilst the dashed edges $f_1$ and $f_2$ are edges of $J$ alone. Taking $T = \{f_1, f_2\}$ and $T' = \{e_2, e_3\}$, we have $\chi(T) - \chi(T') = \chi(\{u\}) - \chi(\{v\})$.}
\label{fig:transferral} 
\end{figure}

Transferrals come into play in the step of removing a partial matching to balance the part sizes. Given clusters $U$ and $U'$ and $b \in \N$, a $b$-fold $(U,U')$-transferral in $(R,M)$ consists of a pair $(T,T')$ of multisets of edges, with $T \sub R$, $T' \sub M$ and $|T|=|T'|$, such that every cluster is covered equally often by $T$ and $T'$, except that $U$ is covered $b$ times more in $T$ than in $T'$, and $U'$ is covered $b$ times more in $T'$ than in $T$. An example is illustrated in Figure~\ref{fig:transferral}. Given such a pair $(T,T')$ for which $U$ is too large and $U'$ is too small for our perfect matching strategy, we can reduce the imbalance by $b$; to achieve this we choose a partial matching in $J$ with edges corresponding to each edge of $T$ (disjoint from all edges chosen so far), and we note for future reference that the perfect matchings corresponding to edges of $M$ chosen in the final step will use one fewer edge corresponding to each edge of $T'$. We say that $(R,M)$ is $(B,C)$-irreducible if for any $U$, $U'$ there is a $b$-fold $(U,U')$-transferral of size $c$ for some $b\le B$ and $c \le C$ (here we are sticking to the non-partite case for simplicity). 

Irreducibility allows us to reduce the cluster imbalances to a constant independent of $n$, and so find a matching in $J_k$ covering all but a constant number of vertices. Furthermore, it can be expressed in terms of the following geometric condition, alluded to in the introduction. Let 
$$X=X(R_k,M) = \{\chi(e)-\chi(e'): e \in R_k,e' \in M\},$$ 
where $\chi(S)$ denotes the characteristic vector of a set $S \sub V(R)$. The required condition is that the convex hull of $X$ should contain a ball of some small constant radius centred at the origin. Thus irreducibility becomes a question of convex geometry. It is also the second point at which space barriers come into play, as our minimum degree assumption implies the existence of such a ball about the origin in the absence of a space barrier. This leads to our first stability result: under the same minimum degree sequence assumption on $J$ needed for a matching in $J_k$ that is almost perfect (i.e. covers all but $o(n)$ vertices), we can in fact find a matching in $J_k$ that covers all but a constant number of vertices, unless $J$ is structurally close to a space barrier construction.

For perfect matchings, we need more precise balancing operations, namely simple transferrals, by which we mean $1$-fold transferrals. We can still use general transferrals for most of the balancing operations, but we require simple transferrals to make the cluster sizes precisely equal. We introduce transferral digraphs $D_\ell(R,M)$, where there is an edge from $U$ to $U'$ if and only if $(R,M)$ contains a simple $(U,U')$-transferral of size at most $\ell$. If we now assume for simplicity that $J$ is a $k$-complex (rather than just a $k$-system), then it is immediate that every vertex in $D_1(R,M)$ has outdegree at least $\delta^*_{k-1}(R)$. We will prove a structure result for digraphs with linear minimum outdegree, which when combined with irreducibility gives a partition $\Part$ of $V(R)$ into a constant number of parts, such that there are simple transferrals between any two vertices in the same part. This is the point at which divisibility barriers come into play. If the robust edge lattice of $R_k$ with respect to $\Part$ is incomplete then $J_k$ is structurally close to a divisibility barrier construction. Otherwise, for any pair of parts $P_i,P_j \in \Part$ we have some simple $(U_i,U_j)$-transferrals with $U_i \in P_i$ and $U_j \in P_j$, and the robustness of the edge lattice gives enough of these transferrals to make the cluster sizes precisely equal. Thus we obtain the full stability theorem:  under the same minimum degree sequence assumption on $J$ needed for a matching in $J_k$ that is almost perfect, we can find a perfect matching in $J_k$, unless $J$ is structurally close to a space or divisibility barrier construction.

\chapter{Geometric Motifs} \label{sec:prelims}

In this chapter we prove or cite various lemmas that underpin the convex geometry in our theory,
and also demonstrate the connection with fractional perfect matchings in hypergraphs.
For the classical results cited we refer the reader to the book by Schrijver \cite{S}.
We begin with Carath\'eodory's theorem, in a slightly unusual form.
Given points $\xb_1, \dots, \xb_s \in \mathbb{R}^d$, we define their \emph{convex hull} as
$$CH(\xb_1, \dots, \xb_s) := \left\{\sum_{j \in [s]} \lambda_j \xb_j : \lambda_j \in [0,1], \sum_{j \in [s]} \lambda_j = 1\right\}.$$

\begin{theo}[Carath\'eodory's Theorem] \label{caratheodory}
Suppose $X \sub \R^d$ and $\xb \in CH(X)$.
Then there are $\lambda_1, \dots, \lambda_{r} \geq 0$ and $\xb_1, \dots, \xb_{r} \in X$ such that
\begin{itemize}
\item[(a)] $\sum_{j \in [r]} \lambda_j = 1$,
\item[(b)] the vectors $\xb_j - \xb_r$ are linearly independent for each $j \in [r-1]$, and
\item[(c)] $\xb = \sum_{j \in [r]} \lambda_j \xb_j$.
\end{itemize}
\end{theo}

Note that condition~(b) implies that $r \leq d+1$. Indeed, a more standard statement of Theorem~\ref{caratheodory} has the condition $r \leq d+1$ in place of (b); this is commonly proved by first proving our formulation of the theorem.
This means that we may write any point in~$CH(X)$ as a positive linear combination of a small number of members of~$X$.
The following proposition gives conditions under which we can arrange that all of the coefficients~$\lambda_j$ are rationals of small denominator. We say that a real number~$x$ is \emph{$q$-rational} if we can write $x=a/b$ for integers $a,b$
with $1  \le b \le q$. We let $\Q^d_q$ denote the set of points~$\xb \in \mathbb{R}^d$ such that every coordinate of~$\xb$ is $q$-rational. Also, given $\xb \in \R^d$ and $r \ge 0$ we let 
$$B^d(\xb,r) = \{\zb: \|\zb-\xb\| \le r\}$$ 
denote the ball of radius $r$ centred at $\xb$ (we sometimes drop the dimension superscript when this is clear from the context).

\begin{prop} \label{coeffs_are_rational}
Suppose that $1/q' \ll 1/q, 1/d, 1/k$, and let~$X$ be a subset of $\Q^d_q \cap B^d(\0, 2k)$. Then for any $\xb \in CH(X) \cap \Q_q^d$ there exist $\lambda_1, \dots, \lambda_r \geq 0$ and $\xb_1, \dots, \xb_r \in X$ such that $r \leq d+1$, each $\lambda_j$ is $q'$-rational, $\sum_{j \in [r]} \lambda_{j} =1$, and $\xb = \sum_{j \in [r]} \lambda_j \xb_j$.
\end{prop}

\proof By Theorem~\ref{caratheodory} we may choose $z_1, \dots, z_{r} \geq 0$ and $\xb_1, \dots, \xb_{r} \in X$ such that $\sum z_j =1$, the vectors $\xb_j - \xb_r$ are linearly independent for each $j \in [r-1]$ (and hence $r \leq d+1$) and $\xb = \sum_{j \in [r]} z_j \xb_j$. We can write $\xb = A\zb + \xb_r$, where $A$ is the $d$ by $r-1$ matrix whose columns are the vectors $\xb_i-\xb_r$, $i \in [r-1]$, and $\zb = (z_1, \dots, z_{r-1})^\intercal$.
Note that $A$ has rank $r-1$, since its columns are linearly independent. Let $\vb_1, \dots, \vb_d$ be the row vectors of $A$.
Then we can choose $S \sub [d]$ of size $r-1$ such that the row vectors $\vb_j$, $j \in S$ are linearly independent.
Let~$B$ be the $r-1$ by $r-1$ square matrix with rows $\vb_j$ for $j \in S$. Then $\det(B) \ne 0$. Also, since every entry of~$A$ has absolute value at most~$4k$, we have $\det(B) \leq (r-1)!(4k)^{r-1}$. (A better bound is available from Hadamard's inequality, but it suffices to use this crude bound which follows by estimating each term in the expansion of the determinant.)
Now we write $\xb' - \xb'_r = B\zb$, where~$\xb'$ is the restriction of $\xb$ to the coordinates $j \in S$, and~$\xb'_r$ is defined similarly. Then $\zb = B^{-1}(\xb' - \xb'_r)$, and so every coordinate of~$\zb$ can be expressed as a fraction with denominator at most $q^2\det(B) \leq q'$, as required.
\endproof

Next we need the classical theorem on the equivalence of vertex and half-plane representations of convex polytopes.
This is commonly known as the Weyl-Minkowski theorem; it is also implied by results of Farkas.
Given points $\vb_1, \dots, \vb_r \in \R^d$, we define their {\em positive cone} as
$$PC(\{\vb_1, \dots, \vb_r\}) := \{\sum_{j \in [r]} \lambda_{j} \vb_j : \lambda_1, \dots, \lambda_r \geq 0\}.$$
The {\em Minkowski sum} of two sets $A,B \sub \R^d$ is $A+B = \{a+b: a \in A, b \in B\}$.

\begin{theo}[Weyl-Minkowski Theorem] \label{weyl-minkowski}
Let $P \sub \R^d$. Then the following statements are equivalent.
\begin{itemize}
\item[(i)] $P = \{\xb \in \R^d : \ab_j \cdot \xb \geq b_j \textrm{ for all } j \in [s]\}$ for some $\ab_1, \dots, \ab_s \in \R^d$ and $b_1, \dots, b_s \in \R$.
\item[(ii)] $P = CH(X) + PC(Y)$ for some finite sets $X, Y \sub \R^d$.
\end{itemize}
\end{theo}

An important case of Theorem~\ref{weyl-minkowski} is when $P =  PC(Y)$ for some finite set $Y \sub \R^d$.
Then we can write $P = \{\xb \in \R^d : \ab_j \cdot \xb \geq 0 \textrm{ for all } j \in S\}$ for some $\ab_1, \dots, \ab_s \in \R^d$, since we have $\0 \in P$, and if $\xb \in P$ then $2\xb \in P$. The following result of Farkas follows.

\begin{lemma}[Farkas' Lemma] \label{farkas}
Suppose $\vb \in \R^d \sm PC(Y)$ for some finite set $Y \sub \R^d$.
Then there is some $\ab \in \R^d$ such that $\ab \cdot \yb \ge 0$ for every $\yb \in Y$ and $\ab \cdot \vb < 0$.
\end{lemma} 

Our next result exploits the discrete nature of bounded integer polytopes.
First we need a convenient description for the faces of a polytope.
A {\em face} of a polytope $P$ is the intersection of $P$ with a set $\mc{H}$ of hyperplanes,
such that for each $H \in \mc{H}$, $P$ is contained in one of the closed halfspaces defined by $H$.
More concretely, consider a finite set $X \sub \R^d$ with $\0 \in X$. By Theorem~\ref{weyl-minkowski}
we can write $CH(X) = \{\xb : \ab_j \cdot \xb \geq b_j \textrm{ for all } j \in [s]\}$
for some $\ab_1,\dots,\ab_s \in \R^d$ and $b_1,\dots,b_s \in \R$. Let $S = \{j \in [s]: b_j = 0\}$,
\[\Pi^X_0 = \{\xb :  \ab_j \cdot \xb = 0 \textrm{ for all } j \in S\}
\textrm{ and }  F^X_0 = CH(X) \cap \Pi^X_0.\]
Then $F^X_0$ is the (unique) minimum face of $CH(X)$ containing $\0$.
Note that in the extreme cases, $F^X_0$ could be all of $CH(X)$, or just the single point $\0$.
The following result gives a lower bound for the distance of $\0$ from the boundary of $F^X_0$.

\begin{lemma} \label{separating_plane_through_zero}
Suppose that $0 < \delta \ll 1/k, 1/d$ and $\0 \in X \sub \Z^d \cap B^d(\0, 2k)$.
Then $\Pi^X_0 \cap B^d(\0, \delta) \sub CH(X)$.
\end{lemma}

\proof
For each $X$ with $\0 \in X \sub \Z^d \cap B^d(\0, 2k)$ we fix a representation $CH(X) = \{\xb : \ab_j^X \cdot \xb \geq b_j^X \textrm{ for all } j \in [s^X]\}$ such that $b_j^X \in \R$ and $\ab_j^X \in \R^d$ with $\|\ab_j^X\|=1$ for $j \in [s^X]$.
Since $\0 \in X$ we must have $b_j^X \leq 0$ for every~$j$ and~$X$. Note that there are only finitely many possible choices of~$X$, and each~$s^X$ is finite and depends only on~$X$. Choosing $\delta$ sufficiently small, we may assume that for any~$j$ and~$X$ either $b_j^X = 0$ or $b_j^X < - \delta$. Now fix any such set~$X$, and write $F^X_0 = CH(X) \cap \Pi^X_0$, where $S^X = \{j \in [s^X]: b_j^X = 0\}$ and $\Pi^X_0 = \{\xb :  \ab^X_j \cdot \xb = 0 \textrm{ for all } j \in S^X\}$. Consider any $\xb \in \Pi^X_0 \sm CH(X)$. Since $\xb \in \Pi^X_0$, $\xb$ satisfies all constraints $\ab_j^X \cdot \xb \geq b_j^X$ for $CH(X)$ with $b_j^X=0$. Since $\xb \notin CH(X)$, $\xb$ must fail some constraint  $\ab_j^X \cdot \xb \geq b_j^X$ for $CH(X)$ with $b_j^X \ne 0$. It follows that $\ab_j^X \cdot \xb < b_j^X < - \delta$. Applying the Cauchy-Schwartz inequality, we have
$\delta < |\ab_j^X \cdot \xb| \le  \|\ab_j^X \| \|\xb\| = \|\xb\|$, as required.
\endproof

Now we will demonstrate the connection between convex geometry and fractional perfect matchings, and how these can be obtained under our minimum degree sequence assumption. A \emph{fractional perfect matching} in a hypergraph $G$ is an assignment of non-negative weights to the edges of $G$ such that for any vertex $v$,
the sum of the weights of all edges incident to~$v$ is equal to $1$. We also require the following notation, which will be used throughout the paper. For any $S \sub [n]$, the \emph{characteristic vector}~$\chi(S)$ of~$S$ is the vector in~$\R^n$ given by
$$ \chi(S)_i = \begin{cases}
1 & i \in S\\
0 & i \notin S.
\end{cases}$$
If~$G$ has~$n$ vertices, then by identifying~$V(G)$ with~$[n]$ we may refer to the characteristic vector~$\chi(S)$ of any $S \sub V(G)$. Whilst~$\chi(S)$ is then dependent on the chosen identification of~$V(G)$ with~$[n]$, the effect of changing this identification is simply to permute the vectors~$\chi(S)$ for each $S \sub V(G)$ by the same permutation, and all properties we shall consider will be invariant under this isomorphism. So we will often speak of~$\chi(S)$ for $S \sub V(G)$ without specifying the identification of~$V(G)$ and~$[n]$. Then a fractional perfect matching is an assignment of weights $w_e \ge 0$ to each $e \in G$ such that $\sum_{e \in G} w_e \chi(e) = \1$; throughout this paper $\1$ denotes a vector of the appropriate dimension in which every coordinate is $1$. Thus $G$ has a fractional perfect matching precisely when $\1$ belongs to the positive cone $PC(\chi(e): e \in G)$. The next lemma shows that this holds under a similar minimum degree sequence assumption to that considered earlier. We include the proof for the purpose of exposition, as a similar argument later in Lemma~\ref{fractionalmatchingFpartite} will simultaneously generalise both this statement and a multipartite version of it. The method of proof used in both cases adapts the separating hyperplane argument of R\"odl, Ruci\'nski and Szemer\'edi~\cite[Proposition~3.1]{RRS}, which was used to prove the existence of a fractional perfect matching in a $k$-graph in which all but a few sets $S$ of $k-1$ vertices satisfy $d(S) \geq n/k$.

\begin{lemma} \label{fractionalmatching}
Suppose that $k \mid n$, and that~$J$ is a $k$-system on~$n$ vertices with
$$\delta(J) \geq \left(n, (k-1)n/k, (k-2)n/k, \dots, n/k \right).$$
Then~$J_k$ admits a fractional perfect matching.
\end{lemma}
 
\proof
Suppose for a contradiction that $J_k$ does not admit a fractional perfect matching. As noted above, this means that
$\1 \notin PC(\chi(e): e \in J_k)$. Then by Farkas' Lemma, there is some $\ab \in \R^n$ such that $\ab \cdot \1 < 0$ and $\ab \cdot \chi(e) \geq 0$ for every $e \in J_k$. Let $v_1, \dots, v_n$ be the vertices of $J$, and let $a_1, \dots, a_n$ be the corresponding coordinates of $\ab$, with the labels chosen so that $a_1 \leq a_2 \leq \dots \leq a_n$. For any sets~$S$ and~$S'$ of~$k$ vertices of~$J$, we say that~$S$ \emph{dominates}~$S'$ if we may write $S = \{v_{i_1}, \dots, v_{i_k}\}$ and $S' = \{v_{j_1}, \dots, v_{j_k}\}$ so that $j_\ell \leq i_\ell$ for each $\ell \in [k]$. Note that if~$S$ dominates~$S'$ then $\ab \cdot \chi(S') \leq \ab \cdot \chi(S)$; this follows from the fact that the coordinates of~$\ab$ are increasing.

For each $j \in [n/k]$, let~$S_j$ be the set $\{v_j, v_{j+n/k}, v_{j+2n/k}, \dots, v_{j + (k-1)n/k}\}$. Then the sets~$S_j$ partition~$V(J)$, so $\sum_{j \in  [n/k]} \ab \cdot \chi(S_j) = \ab \cdot \1$. We claim that there is some edge~$e \in J_k$ which is dominated by every~$S_j$. To see this, we let $d_1 = 1$, then apply the minimum degree sequence condition on $J$ to choose $d_2, \dots, d_k$ greedily so that for each $j \in [k]$ we have $\{v_{d_1}, \dots v_{d_j}\} \in J$ and $d_j \leq (j-1)n/k +1$. Then $e := \{v_{d_1}, \dots, v_{d_k}\}$ is dominated by $S_j$ for each $j \in [n/k]$. We therefore have
$$0 \leq \sum_{j \in [n/k]} \ab \cdot \chi(e) \leq \sum_{j \in [n/k]} \ab \cdot \chi(S_j) = \ab \cdot \1 < 0, $$
which is a contradiction.
\endproof 

\chapter{Transferrals} \label{sec:transferrals}

To motivate the results of the next two chapters, we start by recalling the proof strategy discussed in Section~\ref{sec:outline}. Hypergraph regularity theory (presented in Chapter~\ref{sec:regularity}) will enable us to approximate our original $k$-system by a reduced $k$-system. To avoid confusion over notation, we should emphasise that all of our transferral results will be applied when $J$ is equal to the reduced system, rather than the original system. We will also have a perfect matching $M$ in $J_k$, whose edges represent $k$-tuples of clusters in which we would be able to find a perfect matching, if we were able to make the cluster sizes equal. The role of transferrals is to achieve this by removing a small partial matching. While they are motivated by the proof strategy, the definition and analysis of transferrals does not require any regularity theory (apart from one technical lemma in the next chapter). Suppose that $J$ is a $k$-graph on $[n]$ and $M$ is a perfect matching in $J$; for convenience we call the pair $(J,M)$ a \emph{matched $k$-graph}. Recall that $\chi(e) \in \R^n$ denotes the characteristic vector of $e \sub [n]$. If $T$ is a multiset of subsets of $[n]$ we write $\chi(T) = \sum_{e \in T} \chi(e)$, thus identifying $T$ with the multiset in $[n]$ in which the multiplicity of $i \in [n]$ is the number of sets in $T$ containing $i$, counting with repetition. Also, we say, e.g.~`a multiset $T$ in $J$' to mean a multiset $T$ of members of $J$.

\section{Irreducibility.}
Now we make some important definitions. Given $b \in \N$, $u, v \in V(J)$ and multisets $T$ in $J$ and $T'$ in $M$, we say that $(T, T')$ is a \emph{$b$-fold $(u,v)$-transferral in $(J, M)$} if $$\chi(T) - \chi(T') = b(\chi(\{u\}) - \chi(\{v\})).$$
That is to say that every vertex of $J$ appears equally many times in $T$ as in $T'$, with the exception of $u$, which appears $b$ times more in $T$ than in $T'$, and $v$, which appears $b$ times more in $T'$ than in $T$. An example is shown in Figure~\ref{fig:transferral}. Note that if $(T, T')$ is a $b$-fold $(u,v)$-transferral in $(J, M)$ then we must have $|T| = |T'|$; we refer to this common size as the $\emph{size}$ of the transferral. We say that $(J, M)$ is \emph{$(B,C)$-irreducible} if for any $u,v \in V(J)$ there exist $b \leq B$ and $c \leq C$ such that $(J, M)$ contains a $b$-fold $(u,v)$-transferral of size $c$. We also make the following partite version of this definition, in which we only require transferrals within parts. Given a partition $\Part$ of $V(J)$ into parts $V_1, \dots, V_r$, then we say that $(J, M)$ is \emph{$(B,C)$-irreducible with respect to $\Part$} if for any $i \in [r]$ and any $u,v \in V_i$ there exist $b \leq B$ and $c \leq C$ such that $(J, M)$ contains a $b$-fold $(u,v)$-transferral of size $c$. As described in Section~\ref{sec:outline}, irreducibility will allow us to reduce the cluster imbalances to a constant, and so find a matching in original system covering all but a constant number of vertices.

It is instructive to consider how irreducibility fails in the space barrier constructions (Construction~\ref{spacebar}, for example as shown in Figure~\ref{fig:spacebar}). Consider for simplicity the $3$-complex $J = J(S,2)$ when $|S|=2n/3$ (suppose that $3 \mid n$). By definition we have $|e \cap S| \leq 2$ for every $e \in J_3$. So for any perfect matching $M$ in $J_3$, each $e' \in M$ has $|e' \cap S|= 2$. Thus for any multisets $T, T'$ in $J$, $M$ with $|T|=|T'|$, the vector $\xb = \chi(T) - \chi(T')$ satisfies $\sum_{i \in S} x_i \leq 0$. It follows that there are no $b$-fold $(u,v)$-transferrals for any $b \in \N$, $u \in S$ and $v \notin S$. Conversely, we will see in the next chapter that if there is no space barrier, then the minimum degree sequence assumption implies irreducibility (with respect to the partition in the partite case). In fact, it implies a geometric condition which is the fundamental property behind irreducibility. To formulate this, we define
\[ X=X(J,M) = \{\chi(e)-\chi(e'): e \in J,e' \in M\}.\]
Also, if $\Part$ is a partition of $V$ we define
\[ \Pi_\Part = \{\xb \in \Z^{n} : \xb \cdot \chi(U) = 0 \textrm{ for all } U \in \Part\}.\]
The required geometric condition is that the convex hull of $X$ should contain a ball about the origin within $\Pi_\Part$. The following lemma shows that this implies irreducibility. (For now we gloss over the fact that $B$ and $C$ depend on $n$, which is not permissible in our proof strategy; this dependence will be removed by a random reduction in Lemma~\ref{RANDOMEDGESELECTION}.)

\begin{lemma} \label{ballgivesgentransferrals}
Suppose that $1/B, 1/C \ll \delta, 1/n, 1/k$. Let $V$ be a set of $n$ vertices and let $\Part$ partition $V$. Suppose that $(J, M)$ is a matched $k$-graph on~$V$ such that $X=X(J,M)$ satisfies $B(\0, \delta) \cap \Pi_\Part \sub CH(X)$. Then $(J, M)$ is $(B, C)$-irreducible with respect to $\Part$.
\end{lemma}

\proof Let $q$ be such that $1/B, 1/C \ll 1/q \ll \delta, 1/n, 1/k$. Fix any $U \in \Part$ and $u, v \in U$. Let $t = \lceil 2/\delta \rceil$, and let $\xb \in \Q^n$ have coordinates $x_u=1/t$, $x_v = -1/t$ and $x_w = 0$ for $w \ne u,v$. Then $\|\xb\|<\delta$ and $\xb \in \Pi_\Part$, so our hypothesis implies that $\xb \in CH(X)$. Applying Proposition~\ref{coeffs_are_rational}, we may choose $\lambda_1, \dots, \lambda_s \geq 0$ and $\xb_1, \dots, \xb_s \in X$ such that $s \leq n+1$, each $\lambda_j$ is $q$-rational, $\sum_{j \in [s]} \lambda_{j} =1$, and $\xb = \sum_{j \in [s]} \lambda_j \xb_j$. By definition of $X$, we can choose $e_j \in J$ and $e'_j \in M$ such that $\chi(e_j) - \chi(e_j') = \xb_j$ for $j \in [s]$. Now let the multiset $T$ in $J$ consist of $tq! \lambda_j$ copies of $e_j$ for each $j \in [s]$, and similarly let the multiset $T'$ of edges of $M$ consist of $tq! \lambda_j$ copies of $e_j'$ for each $j \in [s]$. Then $$\chi(T) - \chi(T') = tq! \sum_{j \in s} \lambda_j \xb_j = tq! \xb.$$
So the paired multiset $(T, T')$ is a $q!$-fold $(u,v)$-transferral in $(J, M)$ with size at most $tq! \leq C$. Since $U \in \Part$ and $u,v \in U$ were arbitrary, we deduce that $(J, M)$ is $(B, C)$-irreducible with respect to $\Part$, as required.
\endproof

\section{Transferral digraphs.}
To obtain perfect matchings, we need more precise balancing operations, namely \emph{simple transferrals}, by which we mean $1$-fold transferrals. We represent these using digraphs (directed graphs), for which we make the following standard definitions. A \emph{digraph} $D$ consists of a vertex set $V(D)$ and an edge set $E(D)$, where the edges are each ordered pairs of vertices; we allow loops $(u,u)$, and $2$-cycles $\{(u,v),(v,u)\}$, but do not allow multiple copies of a given ordered pair. We think of an edge $(u, v) \in E(D)$ as being \emph{directed} from $u$ to $v$; if $(u, v) \in E(D)$ then we say that $v$ is an \emph{outneighbour} of $u$ and that $u$ is an \emph{inneighbour} of $v$. For any vertex $v \in V(D)$, the \emph{outneighbourhood} $N^+(v)$ is the set of outneighbours of $v$, and the \emph{inneighbourhood} $N^-(v)$ is the set of inneighbours of $v$. We also define the \emph{outdegree} $d^+(v) = |N^+(v)|$ and the \emph{indegree} $d^-(v) = |N^-(v)|$, and the \emph{minimum outdegree} $\delta^+(D)$ of $D$, which is the minimum of $d^+(v)$ over all $v \in D$.

We represent the simple transferrals in $(J, M)$ by the \emph{$\ell$-transferral digraph of $(J,M$)}, denoted $D_\ell(J, M)$, where $\ell \in \N$. This is defined to be the digraph on $V(J)$ in which there is an edge from $u$ to $v$ if and only if $(J, M)$ contains a simple $(u,v)$-transferral of size at most~$\ell$. For future reference we record here some basic properties of transferral digraphs:
\begin{enumerate}[(i)]
\item $D_\ell(J,M)$ contains all loops $(u,u)$, as a $(u,u)$-transferral is trivially achieved by the empty multisets.
\item If $\ell \leq \ell'$ then $D_\ell(J, M)$ is a subgraph of $D_{\ell'}(J, M)$.
\item If $(u,v) \in D_{\ell}(J,M)$ and $(v,w) \in D_{\ell'}(J,M)$ then $(v,w) \in D_{\ell+\ell'}(J,M)$. Indeed, if $(S,S')$ is a simple $(u,v)$-transferral of size at most $\ell$ and $(T,T')$ is a simple $(v,w)$-transferral of size at most $\ell'$, then $(S+T,S'+T')$ is a simple $(u,w)$-transferral of size at most $\ell+\ell'$.
\end{enumerate}
In the last property we used the notation $A+B$ for the multiset union of multisets $A$ and $B$. We will also use the notation $\sum_{i \in [t]} A_i$ for the multiset union of multisets $A_i$, $i \in [t]$, and the notation $pA$ for the multiset union of $p$ copies of $A$.
 
Again, it is instructive to consider the examples of the space barrier and divisibility barrier constructions (Constructions~\ref{spacebar} and~\ref{divbar}, for example as shown in Figures~\ref{fig:spacebar} and~\ref{fig:divbar}). As described in the last section, if $J$ is the 3-complex $J(S, 2)$ from Construction~\ref{spacebar}, where $|S| = 2n/3$, and $M$ is a perfect matching in $J$, then the matched 3-complex $(J, M)$ does not contain a $b$-fold $(u, v)$-transferral for any $b \in \mathbb{N}, u \in S$ and $v \notin S$. However, it is not hard to see that $(J, M)$ does contain a simple $(u, v)$-transferral of size 1 for any other $u$ and $v$. So for any $\ell \in \mathbb{N}$ the sets $S$ and $V \sm S$ both induce complete subgraphs in the transferral digraph $D_\ell(J, M)$, and $D_\ell(J, M)$ also contains all edges directed from $V \sm S$ to $S$. Now instead consider the $k$-graph $G$ on vertex set $V = V_1 \cup V_2 \cup V_3$ whose edges are all $k$-tuples $e$ with $|e \cap V_2| = |e \cap V_3| \mod 3$ (this is example (2) following Construction~\ref{divbar}), and let $M$ be a perfect matching in $G$. Then the matched $k$-graph $(G, M)$ contains a simple $(u, v)$-transferral of size 1 for any $u, v$ with $u, v \in V_j$ for some $j \in [3]$, but there are no simple $(u, v)$-transferrals (of any size) with $u \in V_i$ and $v \in V_j$ for $i \neq j$. So for any $\ell \in \mathbb{N}$ the transferral digraph $D_\ell(G, M)$ is the disjoint union of complete digraphs on $V_1$, $V_2$ and $V_3$.


As described in the proof outline earlier, our strategy for finding a perfect matching in a $k$-complex $J$ uses simple transferrals in the matched $k$-system $(R, M)$ to `rebalance' the sizes of clusters of the reduced $k$-system $R$. For this we need the transferral digraph $D_\ell(R, M)$ to be complete for some $\ell$ (or complete on each part of $\Part$ if $R$ is $\Part$-partite). In this way the transferral digraphs of the two constructions considered above capture the essence of why these constructions do not contain a perfect matching.

Our goal in the remainder of this chapter is to describe the conditions under which, given a matched $k$-graph $(J,M)$ on $V$ and a partition $\Part$ of $V$, there is some transferral digraph that is complete on each part of $\Part$, i.e.\ there is some $\ell$ so that $(u,v) \in D_\ell(J,M)$ for every $u,v \in U \in \Part$. This is a convenient context that in fact gives more general results than the two cases that are needed for our main theorems, namely (i) when $\Part$ consists of a single part $V$, and (ii) when $\Part$ consists of at least $k$ parts and $(J,M)$ is $\Part$-partite. In this section, we prove a structure theorem for general digraphs under a minimum outdegree condition. This introduces a partition $\Part'$ that refines $\Part$, and we show that irreducibility implies that there is some $C'$ for which $D_{C'}(J,M)$ is complete on the parts of $\Part'$. In the final section we show that if there is no divisibility barrier then there is some $C$ for which $D_C(J,M)$ is complete on the parts of $\Part$, as required.

We require the following lemma, which states that in a digraph with linear minimum outdegree, one can choose a constant number of vertices so that from any fixed vertex there is a short directed path to one of these chosen vertices.

\begin{lemma} \label{digraphpartslemma}
Suppose that $1/\ell \ll \alpha$, and let $D$ be a digraph on $n$ vertices. Then we may choose vertices $v_1, \dots, v_t \in V(D)$ and sets $V_1, \dots, V_t \sub V(D)$ such that
\begin{itemize}
\item[(a)] for any $j$ and any $u \in V_j$ there exists a directed path from $u$ to $v_j$ of length at most~$\ell$,
\item[(b)] the sets $V_1, \dots, V_t$ partition $V(D)$, and
\item[(c)] $|V_j| \geq \delta^+(D) - \alpha n$ for each $j$.
\end{itemize}
\end{lemma}

\proof
Let $\gamma'$ satisfy $1/(\ell-1) \ll \gamma' \ll \alpha$. Clearly we may assume that $\delta^+(D) \geq \alpha n$. For any $S \sub V(D)$, we let $L(S)$ be the number of edges of $D$ that {\em leave} $S$, i.e.\ are of the form $(u,v)$ with $u \in S$ and $v \in \ov{S} = V(D) \sm S$. We require the following claim.

\medskip

\begin{claim}\label{digraphpartsclaim}
Let $1/(\ell-1) \leq \gamma \leq \gamma'$, and let $S \sub V(D)$ be such that $|S| \geq \alpha n$ and $L(S) \le \gamma n^2$. Then we may choose $v \in S$ and $S' \sub S$ such that $|S'| > \delta^+(D) - \alpha n$, $L(S \sm S') \le 2\gamma n^2$, and for any $u \in S'$ there is a path from $u$ to $v$ in $D$ of length at most $\ell - 1$.
\end{claim}

To prove the claim, we start by showing that there must be some $v \in S$ with $d^-_S(v) := |N^-(v) \cap S| \geq \delta^+(D) - \alpha n$. To see this, we note that $\sum_{v \in S} d^-_S(v) \ge \sum_{v \in S} d^+(v) - L(S) \ge |S|\delta^+(D)-\gamma n^2$. Thus by averaging we can choose $v \in S$ with $d^-_S(v) \ge \delta^+(D)-\gamma n^2/|S| \ge \delta^+(D) - \alpha n$. Now consider the `iterated inneighbourhood' $N^-_j$ of $v$ for $j \ge 1$, defined as the set of vertices $u \in S$ such that there exists a path from $u$ to $v$ in $D$ of length at most $j$. Note that $N_1^-\sub N_2^- \sub \dots$, so we can choose some $j \leq 1/\gamma$ for which $|N^-_{j+1}| \leq |N^-_{j}| + \gamma n$. Now we set $S' = N^-_j$. By construction this satisfies the required property of paths (since $j \leq 1/\gamma \leq \ell-1$). We also have $|S'| \geq |N^-_1| \geq \delta^+(D) - \alpha n$, so it remains to estimate $L(S \sm S')$. Any edge leaving $S \sm S'$ is either an edge leaving $S$ or an edge from $S \sm S'$ to $S'$. By assumption there are at most $\gamma n^2$ edges of the first type. Also, for any edge $(u,v)$ of the second type we have $u \in N^-_{j+1} \sm N^-_j$, so there are at most $|N^-_{j+1} \sm N^-_j|n \le \gamma n^2$ such edges. This completes the proof of Claim~\ref{digraphpartsclaim}.

\medskip

Returning to the proof of the lemma, we apply the claim to $S=V(D)$ with $\gamma = 1/(\ell-1)$, obtaining $v \in S$ and $S' \sub S$ with the stated properties. We set $v_1=v$ and $V_1=S'$. Then we repeatedly apply the claim for each $j \ge 1$ to $S = V(D) \sm (V_1 \cup \dots \cup V_j)$ with $\gamma = 2^j/(\ell-1)$; in each application we obtain $v \in S$ and $S' \sub S$ for the current set $S$ and we set $v_{j+1}=v$ and $V_{j+1}=S'$. Since $\delta^+(D) \geq \alpha n$, we can repeat this process until we reach some $t \ge 1$ for which $S = V(D) \sm (V_1 \cup \dots \cup V_t)$ has $|S| < \delta^+(D)$. Then any $u \in S$ must have an outneighbour $u' \in V(D) \sm S$, i.e. $u' \in V_j$ for some $j \in [t]$. By choice of $V_j$ there is a path of length at most $\ell-1$ from $u'$ to $v_j$, so we have a path of length at most $\ell$ from $u$ to $v_j$. For each such~$u$, choose any such $j$ arbitrarily, and add $u$ to $V_j$. Then the vertices $v_1, \dots, v_t$ and the new sets $V_1, \dots, V_t$ satisfy the properties required for the lemma. \endproof

Now we define the structure we will find in our digraphs. Suppose $D$ is a digraph. We say that $S \sub V(D)$ is a \emph{dominated set in $D$} if $(u,v)$ is an edge of $D$ for every $u \in V(D)$ and $v \in S$. Now let $\Part$ be a partition of $V(D)$ into parts $V_1, \dots, V_r$. We say that $\Sart = (S_1, \dots, S_r)$ with $S_1 \sub V_1, \dots, S_r \sub V_r$ is a \emph{dominated $\Part$-tuple} if each $S_j$ is a dominated set in $D[V_j]$. If such an $\Sart$ exists, we say that $\Part$ is a \emph{dominated partition} of $D$. The \emph{$\ell$th power} of $D$ is the digraph $D^\ell$ with vertex set $V(D)$ where $(u, v) \in D^\ell$ if and only if there exists a path from $u$ to $v$ in $D$ of length at most $\ell$. Now we can state our structural result, which is that some power of $D$ admits a dominated partition.

\begin{coro} \label{digraphparts}
Suppose that $1/\ell \ll \alpha$, and let $D$ be a digraph on $n$ vertices. Then $D^\ell$ admits a dominated partition $\Part$ in which each part has size at least $\delta^+(D) - \alpha n$.
\end{coro}

\proof Apply Lemma~\ref{digraphpartslemma} to obtain sets $V_1, \dots, V_t \sub V(D)$ and $v_1, \dots, v_t \in V(D)$ for some $t$. Let $\Part$ be the partition of $V(D)$ into parts $V_1, \dots, V_t$. So each part of $\Part$ has size at least  $\delta^+(D) - \alpha n$. Now, for any $j \in [t]$ and any $u \in V_j$ there is a path in $D$ from $u$ to $v_j$ of length at most $\ell$, and therefore an edge in $D^\ell$ from $u$ to $v_j$. Then $\Sart := (\{v_1\}, \dots, \{v_t\})$ is a dominated $\Part$-tuple in $D^\ell$, as required.
\endproof

Now we show that in combination with irreducibility, a dominated partition $\Part$ has the property that some transferral digraph is complete on every part of $\Part$. We need the following lemma.

\begin{lemma} \label{edges_bidirectional}
Let $(J, M)$ be a matched $k$-graph, $b,c,c' \in \N$, and vertices $u, v \in V(J)$ be such that $(J, M)$ contains a $b$-fold $(u,v)$-transferral of size $c$, and $(J, M)$ contains a simple $(v,u)$-transferral of size $c'$.
Then $(J, M)$ contains a simple $(u,v)$-transferral of size $(b-1)c'+c$.
\end{lemma}

\proof By definition of transferrals, there are multisets $T_1,T_2$ in $J$ and $T'_1,T'_2$ in $M$ such that $|T_1| = |T'_1| = c$, $|T_2| = |T'_2| = c'$, $$\chi(T_1) - \chi(T'_1) = b(\chi(\{u\}) - \chi(\{v\})), \textrm{ and }
\chi(T_2) - \chi(T'_2) = \chi(\{v\}) - \chi(\{u\}).$$
Let $T = T_1 + (b-1)T_2$ and $T' = T'_1 + (b-1)T'_2$. Then $|T| = |T'| = (b-1)c' + c$, and
$$\chi(T) - \chi(T') = \chi(T_1) - \chi(T'_1) + (b-1)(\chi(T_2)-\chi(T'_2)) =  \chi(\{u\}) - \chi(\{v\}),$$
so $(T,T')$ is the desired transferral.
\endproof

Finally, we show the required completeness property of the transferral digraph.

\begin{lemma} \label{complete_on_receiving}
Suppose that $1/\ell' \ll 1/B,1/C, 1/\ell \ll \alpha, \gamma$. Let $\Part'$ be a partition of a set $V$ of size $n$ and $(J, M)$ be a matched $k$-graph on $V$ which is $(B, C)$-irreducible with respect to $\Part'$. Suppose also that $\delta^+(D_\ell(J, M)[U]) \geq \gamma n$ for every $U \in \Part'$. Then $D_{\ell^2}(J,M)$ admits a partition $\Part$ refining $\Part'$ in which every part has size at least $(\gamma-\alpha)n$, and $D_{\ell'}(J,M)[U]$ is complete for every $U \in \Part$.
\end{lemma}

\proof
Applying Corollary~\ref{digraphparts} to $D_\ell(J,M)[U]$ for each $U \in \Part'$, we see that $D_\ell(J,M)^\ell$ admits a dominated partition $\Part$ refining $\Part'$ in which every part has size at least $(\gamma-\alpha)n$. We have $D_\ell(J,M)^\ell \sub D_{\ell^2}(J,M)$ by property (iii) of transferral digraphs, so $\Part$ is also a dominated partition for $D_{\ell^2}(J,M)$. Now suppose $U \in \Part$ and $S \sub U$ is the dominated set for $D_{\ell^2}(J,M)[U]$. Consider any $u,v \in U$. We need to show that there is a simple $(u,v)$-transferral of size at most $\ell'$. Fix any $w \in S$. By choice of $\Part$, we have a simple $(u,w)$-transferral and a simple $(v,w)$-transferral each of size at most $\ell^2$. By irreducibility, we also have a $b$-fold $(w,v)$-transferral of size $c$ for some $b \le B$ and $c \le C$. Applying Lemma~\ref{edges_bidirectional}, we have a simple $(w,v)$-transferral of size $(b-1)\ell^2 + c$. Combining this with the simple $(u,w)$-transferral, we obtain a simple $(u,v)$-transferral of size at most $b\ell^2+c \le \ell'$.
\endproof

\section{Completion of the transferral digraph.}

In this final section, we show that if in addition to the previous assumptions there is no divisibility barrier, then we can strengthen the previous structure to obtain a transferral digraph that is complete on every part of the original partition. We start with a useful consequence of irreducibility in the following lemma, which allows us to extend a pair of multisets in $J$ and $M$ to a pair of multisets that cover the same vertices (including multiplicities). This corresponds to the geometric intuition that we can use a ball about the origin to counterbalance any given vector. We need to impose the following property on $(J,M)$ that will be automatically satisfied in our applications: we say that $(J,M)$ is \emph{$\Part$-proper} if for any $e \in J$ there is $e' \in M$ with $\ib_\Part(e')=\ib_\Part(e)$.

\begin{lemma} \label{cancel_transferral_partite}
Suppose that $1/C' \ll 1/B, 1/C, 1/k$, let $\Part$ be a partition of a set of vertices~$V$, and let $(J, M)$ be a matched $k$-graph on $V$ which is $\Part$-proper and $(B,C)$-irreducible with respect to $\Part$. Then for any multisets $A$ in $J$ and $A'$ in $M$ with $|A|,|A'| \le C$ there exist multisets $T$ in $J$ and $T'$ in $M$ such that $A \sub T$, $A' \sub T'$, $|T| = |T'| \leq C'$ and $\chi(T) = \chi(T')$.
\end{lemma}

\proof
First we find multisets $S$ in $J$ and $S'$ in $M$ such that $A \sub S$, $A' \sub S'$ and $\sum_{e \in S} \ib_\Part(e) = \sum_{e' \in S'} \ib_\Part(e')$. To do this, we let $S = A+A'$ and $S' = A^*+A'$, where $A^*$ is formed by taking some $e' \in M$ with $\ib_\Part(e') = \ib_\Part(e)$ for every $e \in A$ (which exists because $(J,M)$ is $\Part$-proper). Since $\sum_{e \in S} \ib_\Part(e) = \sum_{e' \in S'} \ib_\Part(e')$ we may label the vertices of $S$ and $S'$ as $u_1, \dots, u_{k|S|}$ and $v_1, \dots, v_{k|S|}$ respectively so that $u_i$ and $v_i$ lie in the same part of $\Part$ for each $i \in [k|S|]$. By irreducibility, for each $i \in [k|S|]$ there is a $b$-fold $(u_i,v_i)$-transferral in $(J,M)$ of size $c$ for some $b \le B$ and $c \le C$. Combining $B!/b$ copies of this transferral, we obtain a $B!$-fold $(u_i,v_i)$-transferral $(T_i, T'_i)$ of size at most $CB!$. Let $T = B!S + \sum_{i \in [k|S|]} T_i$ and $T' = B!S' + \sum_{i \in [k|S|]} T'_i$. Then $|T| = |T'| \leq CB! + 2Ck \cdot CB! \leq C'$. Also,
\begin{align*}
\chi(T)-\chi(T') &= B!(\chi(S) - \chi(S')) + \sum_{i \in [k|S|]} \chi(T_i)-\chi(T'_i) \\
& = B! \sum_{i \in [k|S|]} \big(\chi(\{u_i\}) - \chi(\{v_i\}) - \chi(\{u_i\}) + \chi(\{v_i\})\big) = 0,
\end{align*}
which completes the proof.
\endproof

Next we need a simple proposition that allows us to efficiently represent vectors in a lattice, in that we have bounds on the number of terms and size of the coefficients in the representation.

\begin{prop} \label{vector_as_small_sum}
Suppose that $1/s \ll 1/k, 1/d$. Let $X \sub \Z^d \cap B^d(\0, 2k)$, and let $L_X$ be the sublattice of $\Z^d$ generated by $X$. Then for any vector $\xb \in L_X \cap B^d(\0, 2k)$ we may choose multisets $S_1$ and $S_2$ of elements of $X$ such that $|S_1|, |S_2| \leq s$ and $\sum_{\vb \in S_1} \vb - \sum_{\vb \in S_2} \vb = \xb$.
\end{prop}

\proof
The number of pairs $(X, \xb)$ as in the statement of the proposition depends only on $k$ and $d$. Furthermore, for any such pair $(X, \xb)$ we may choose multisets $S_1$ and $S_2$ of elements of $X$ such that $\sum_{\vb \in S_1} \vb - \sum_{\vb \in S_2} \vb = \xb$. Let $s_{X, \xb}$ be minimal such that we may do this with $|S_1|, |S_2| \leq s_{X, \xb}$. Since $1/s \ll 1/k, 1/d$, we may assume that $s \geq \max_{(X, \xb)} s_{X, \xb}$. Thus for any $X$ and $\xb$ we may express $\xb$ in the required manner.
\endproof

After these preparations, we are ready to prove the main lemma of this chapter. Suppose that $\Part'$ is a partition of
a set $V$ and $\Part$ is a partition into $d$ parts refining $\Part'$. Given a $k$-graph $G$ on $V$, we define the \emph{edge lattice} $L_\Part(G) \sub \Z^d$ to be the lattice generated by all vectors $\ib_\Part(e)$ with $e \in G$. (Note that this definition is similar to that of the robust edge lattices defined earlier; indeed the purpose of robustness is for the same edge lattice to be inherited by the reduced graph.) Recall that a lattice $L \sub \Z^d$ is complete with respect to $\Part'$ if $L_{\Part\Part'} \sub L \cap \Pi^d$, where  $L_{\Part\Part'} \sub \Z^d \cap \Pi^d$ is the lattice generated by all differences of basis vectors $\unit_i-\unit_j$ for which $V_i,V_j$ are contained in the same part of $\Part'$ and $\Pi^d = \{\xb \in \R^d: \sum_{i \in [d]} x_i = 0\}$.

\begin{lemma} \label{make_rectuple_spanning_partite}
Suppose that $1/C' \ll 1/B, 1/C, 1/d, 1/k$. Let $\Part'$ be a partition of a set of vertices $V$, and let $(J, M)$ be a matched $k$-graph on $V$ which is $\Part'$-proper and $(B,C)$-irreducible with respect to $\Part'$. Suppose that $\Part$ is a partition of $V$ into $d$ parts $V_1, \dots, V_d$ which refines $\Part'$ such that $D_C(J, M)[V_j]$ is complete for each $j \in [d]$ and $L_\Part(J)$ is complete with respect to $\Part'$. Then $(J, M)$ is $(1, C')$-irreducible with respect to $\Part'$, i.e.\ $D_{C'}(J, M)[U]$ is complete for each $U \in \Part'$.
\end{lemma}

\proof
Introduce new constants with $1/C' \ll 1/C_2 \ll 1/C_1 \ll 1/C_0 \ll 1/B, 1/C, 1/d, 1/k$. Fix any $i,j \in [d]$ with $i \neq j$ such that $V_i$ and $V_j$ are both subsets of the same part of $\Part'$ (if no such $i, j$ exist then there is nothing to prove). Then $\xb := \unit_i - \unit_j \in L_{\Part\Part'} \sub L_\Part(J)$ by assumption. By Proposition~\ref{vector_as_small_sum} we may therefore choose multisets $S$ and $T$ in $J$ such that $\ib_\Part(S) - \ib_\Part(T) = \xb$ and $|S|, |T| \leq C_0$. Note that 
$$\sum_{j \in [d]} {i_\Part(S)}_j - \sum_{j \in [d]} {i_\Part(T)}_j = \sum_{j \in [d]} x_j = 0,$$ $\sum_{j \in [d]} {i_\Part(S)}_j = k|S|$ and $\sum_{j \in [d]} {i_\Part(T)}_j = k|T|$, so $S$ and $T$ are of equal size. So we may label the vertices of $S$ as $u_1, \dots, u_r$ and $T$ as $v_1, \dots, v_r$ so that $u_1 \in V_i$, $v_1 \in V_j$, and $u_\ell$ and $v_\ell$ lie in the same part of~$\Part$ for every $2 \leq \ell \leq r$. Since $D_C(J, M)[V_\ell]$ is complete for each $\ell \in [r]$, for each $2 \leq \ell \leq r$ we may choose $(T_\ell, T'_\ell)$ to be a simple $(u_\ell,v_\ell)$-transferral of size at most~$C$. Then
\begin{align*}
\chi(S) - \chi(T) &= \sum_{\ell \in [r]} \chi(\{u_\ell\}) - \chi(\{v_\ell\}) = \chi(\{u_1\}) - \chi(\{v_1\})+ \sum_{2 \leq \ell \leq r} \chi(T_\ell) - \chi(T'_\ell).
\end{align*}
Now, by Lemma~\ref{cancel_transferral_partite} there exist multisets $A_1$ in $J$ and $A'_1$ in $M$ such that $|A_1| = |A'_1| \leq C_1$, $S \sub A_1$ and $\chi(A_1) = \chi(A'_1)$. Let $A_2=A_1-S+T$ be formed by deleting $S$ from $A_1$ and replacing it by $T$. Then $\chi(A_2) - \chi(A'_1) = \chi(A_1) - \chi(S) + \chi(T) - \chi(A'_1) = \chi(T) - \chi(S)$. Finally, let $A_3 = A_2 + \sum_{\ell=2}^r T_\ell$ and $A'_3 = A'_1 + \sum_{\ell=2}^r T'_\ell$. Then
\begin{align*}
\chi(A_3) - \chi(A'_3) &= \chi(A_2) - \chi(A'_1) + \sum_{2 \leq \ell \leq r} \chi(T_\ell) - \chi(T'_\ell) \\ &= \chi(\{v_1\}) - \chi(\{u_1\}).
\end{align*}
So $(A_3, A'_3)$ is a simple $(v_1, u_1)$-transferral of size at most $C_2$. Since $D_C(J, M)[V_i]$ and $D_C(J, M)[V_j]$ are each complete, we deduce that there exists a simple $(v,u)$-transferral of size at most $C_2 + 2C \leq C'$ for each $u \in V_i$ and $v \in V_j$. This holds for any $i,j$ such that $V_i$ and $V_j$ are both subsets of the same part $U$ of $\Part'$, so $D_{C'}(J, M)[U]$ is complete for each $U \in \Part'$.
\endproof

Combining Lemmas~\ref{make_rectuple_spanning_partite} and~\ref{complete_on_receiving} we obtain the following result. 

\begin{lemma} \label{complete_digraph}
Suppose that $1/\ell' \ll 1/B,1/C, 1/\ell \ll \alpha, \gamma$. Let $\Part'$ be a partition of a set $V$ of size $n$ and $(J, M)$ be a matched $k$-graph on $V$ that is $\Part'$-proper and $(B, C)$-irreducible with respect to $\Part'$. Suppose that $\delta^+(D_\ell(J,M)[U]) \geq \gamma n$ for any $U \in \Part'$, and $L_\Part(J)$ is complete with respect to $\Part'$
 for any partition $\Part$ of $V(J)$ into $d$ parts of size at least $(\gamma-\alpha)n$ that refines $\Part'$. Then $D_{\ell'}(J,M)[U]$ is complete for every $U \in \Part'$.
\end{lemma}

\chapter{Transferrals via the minimum degree sequence} \label{sec:mindeg}

In the previous chapter we demonstrated how irreducibility and completeness of the edge lattice imply the existence of all possible simple transferrals. In this chapter we will obtain the same result in the minimum degree sequence settings of our main theorems. In fact we will work in the following setting that simultaneously generalises the partite and non-partite settings. 

Let $J$ be a $k$-system on $V$ and $\Part$ be a partition of $V$ into sets $V_1, \dots, V_r$. An \emph{allocation function} is a function $f:[k] \to [r]$. Intuitively an allocation function should be viewed as a rule for forming an edge of $J$ by choosing vertices one by one; the $i$th vertex should be chosen from part $V_{f(i)}$ of $\Part$. This naturally leads to the notion of the \emph{minimum $f$-degree sequence} $$\delta^f(J) := \left(\delta^f_0(J),\dots,\delta^f_{k-1}(J)\right),$$ where $\delta^f_j(J)$ is the largest $m$ such that for any $\{v_1, \dots, v_j\} \in J$ with $v_i \in V_{f(i)}$ for $i \in [j]$ there are at least $m$ vertices $v_{j+1} \in V_{f(j+1)}$ such that $\{v_1, \dots, v_{j+1}\} \in J$. So $\delta^f_{i-1}(J)$ is a lower bound on the number of choices for the $i$th vertex in the process of forming an edge of $J$ described above. For any allocation function $f$, the index vector $\ib(f)$ of $f$ is the vector 
$$\ib(f) := (|f^{-1}(1)|, \dots, |f^{-1}(r)|) \in \Z^r$$ 
whose $j$th coordinate is the cardinality of the preimage of $j$ under $f$. So any edge $e \in J_k$ formed by the described process will have $\ib(e) = \ib(f)$.

Let $I$ be a multiset of index vectors with respect to $\Part$ such that $\sum_{j \in [r]} {i_j} = k$ for each $\ib = (i_1, \dots, i_r)$ in $I$. Then we may form a multiset $F$ of allocation functions $f : [k] \to [r]$ as follows. For each $\ib \in I$ (with repetition) choose an allocation function $f$ with $\ib(f) = \ib$, and include in $F$ each of the $k!$ permutations $f^\sigma$ for $\sigma \in \Permute_k$ (again including repetitions), where $f^\sigma(i) = f(\sigma(i))$ for $i \in [k]$. Note that the multiset $F$ so obtained does not depend on the choices of allocation function $f$. Also observe that $|F| = k!|I|$, and for any $f \in F$ with $\ib(f) = \ib \in I$, the multiplicity of $f$ in $F$ is $m_\ib \prod_{j \in [r]} i_j!$, where $m_\ib$ is the multiplicity of $\ib$ in $I$. If $F$ can be obtained in this manner, then we say that $F$ is an \emph{allocation}, and we write $I(F)$ for the multiset $I$ from which $F$ was obtained. We say that an allocation $F$ is \emph{$(k,r)$-uniform} if for every $i \in [k]$ and $j \in [r]$ there are $|F|/r$ functions $f \in F$ with $f(i)=j$. We also say that $F$ is \emph{connected} if there is a connected graph $G_F$ on $[r]$ such that for every $jj' \in E(G_F)$ and $i, i' \in [k]$ with $i \ne i'$ there is $f\in F$ with $f(i)=j$ and $f(i')=j'$.

All the allocations $F$ considered in this paper will have the property that $I(F)$ is a set (as opposed to a multiset). In this case we have $I(F) = \{\ib(f) : f \in F\}$, and so $|I| \leq r^k$ and $|F| \leq k!|I| \leq k!r^k$. However, we allow the possibility that $I(F)$ is a multiset as this may be useful for future applications. In this case, we shall usually bound $|F|$ by a function $D_F(r,k)$. Also, the reader should take care with any statement regarding (e.g.) $\ib \in I(F)$ to interpret the statement with multiplicity. We slightly abuse the notation along these lines, writing (e.g.) ``the number of edges of index~$\ib$ is constant over all $\ib \in I(F)$" to mean that the number of edges of index~$\ib$ is a constant multiple of the multiplicity of $\ib$ in $I(F)$.

For an allocation $F$, a $k$-system $J$ on $V$ is \emph{$\Part F$-partite} if for any $j \in [k]$ and $e \in J_j$ there is some $f \in F$ so that $e = \{v_1,\dots,v_j\}$ with $v_i \in V_{f(i)}$ for $i \in [j]$ (so every edge of $J$ can be constructed through the process above for some $f \in F$). The \emph{minimum $F$-degree sequence} of $J$ is then defined to be $$\delta^F(J) := (\delta^F_0(J),\dots,\delta^F_{k-1}(J)),$$ where $\delta^F_j(J) = \min_{f \in F} \delta^f_j(J)$. We note two special cases of this definition that recover our earlier settings. In the case $r=1$ (the non-partite setting), there is a unique function $f:[k] \to [1]$. Let $F$ consist of $k!$ copies of $f$; then $J$ is $\Part F$-partite and $\delta^F(J) = \delta(J)$ is the (usual) minimum degree sequence. In the case when $r \ge k$ and $J$ is $\Part$-partite we can instead choose $F$ to consist of all injective functions from $[k]$ to $[r]$, and then being $\Part$-partite is equivalent to being $\Part F$-partite and  $\delta^F(J)=\delta^*(J)$ is the partite minimum degree sequence. Note that in both cases $F$ is $(k,r)$-uniform and connected.

Recall that we established a geometric criterion for irreducibility in Lemma~\ref{ballgivesgentransferrals}. Our first lemma states that if there is no space barrier then this criterion holds under our generalised minimum degree assumption. Note we have slightly altered the formulation of the space barrier from that used earlier, so that we can state the lemma in the more general context of $k$-systems; it is not hard to demonstrate that in the case of $k$-complexes it is an essentially equivalent condition (up to changing the value of $\beta$). We say that $(J,M)$ is a \emph{matched $k$-system} if $J$ is a $k$-system and $M$ is a perfect matching in $J_k$. In this case we write $X = X(J,M) = \{\chi(e)-\chi(e'): e \in J_k,e' \in M\}$ and $D_\ell(J, M)$ to mean $D_\ell(J_k, M)$.

\begin{lemma} \label{gen_trans_1_partite}
Suppose that $\delta \ll 1/n \ll \alpha \ll \beta, 1/D_F, 1/r, 1/k$, and $k \mid n$. Let $\Part'$ be a partition of a set $V$ into sets $V_1, \dots, V_r$ each of size $n$, $F$ be a $(k,r)$-uniform connected allocation with $|F| \leq D_F$, and $(J, M)$ be a $\Part' F$-partite matched $k$-system on $V$ such that
$$\delta^F(J) \geq \left(n, \left(\frac{k-1}{k} - \alpha \right)n, \left(\frac{k-2}{k} - \alpha \right) n, \dots, \left(\frac{1}{k} - \alpha \right) n\right).$$
Suppose that for any $p \in [k-1]$ and sets $S_i \sub V_i$ such that $|S_i| = pn/k$ for each $i \in [r]$ we have $|J_{p+1}[S]| \geq \beta n^{p+1}$, where $S = S_1 \cup \dots \cup S_r$.
Then $$ B^{rn}(\0, \delta) \cap \Pi_{\Part'} \sub CH(X(J,M)).$$
\end{lemma} 

We first sketch the main ideas of the proof for the non-partite case (where $r = 1$, and $\Part'$ has only one part, namely $V$). Suppose for a contradiction that $B^{n}(\0, \delta) \cap \Pi^k \not\sub CH(X)$, where $X = X(J,M)$. We apply Lemma~\ref{separating_plane_through_zero} to deduce that then there must exist some $\ab \in \R^n$ such that $\ab \cdot \xb \leq 0$ for every $\xb \in X$ and $\ab$ is not constant on $V$ (i.e. $\ab$ is not a constant multiple of $\1$). We use this $\ab$ to partition $V$ into $k$ parts $V_1, \dots, V_k$. Indeed, this partition is chosen so that $|V_1| = n/k - \alpha n$, $|V_2| = |V_3| = \dots = |V_{k-1}| = n/k$, $|V_k| = n/k + \alpha n$, and moreover for any $i < j$ and any vertices $u \in V_i$ and $v \in V_j$ we have $a_u \leq a_v$, where $a_u$ and $a_v$ are the coordinates of $\ab$ corresponding to $u$ and $v$ respectively. Let $u_j$ be the maximum value of $\ab$ on part $V_j$. 

The reason for choosing parts of these sizes is that, due to our assumption on the minimum degree sequence of $J$, given any edge $e$ of $J_{t}$ (for some $t \in [k-1]$) we may greedily extend $e$ to an edge $e \cup \{v_{t+1}, \dots, v_{k}\} \in J_k$ in which $v_{j} \in V_{k-j+1}$ attains at least the maximum value $u_{k-j+1}$ of $\ab$ on $V_{k-j+1}$. We use this fact repeatedly in the proof.


Let $v_n \in V_k$ be a vertex of $V$ with maximum $\ab$-coordinate; then greedily extending $\{v_n\} \in J_1$ to an edge $e \in J_k$ as described above gives $\ab \cdot \chi(e) \geq \sum_{j \in [k]} u_j =: U$. In particular, the definition of $\ab$ then implies that any edge $e' \in M$ must have $\ab \cdot \chi(e') \geq U$. 
However, since $M$ is a perfect matching in $J_k$ we must have $\sum_{e' \in M} \ab' \cdot \chi(e') = \ab' \cdot \1 \approx U|M'|$. Indeed, the final approximation would be an equality if every $V_j$ had the same size. Since $\ab' \cdot \chi(e') \geq \ab' \cdot \chi(e)$ for any $e \in J_k$, we conclude that $\ab' \cdot \chi(e') \approx \ab \cdot \chi(e') \approx U$ for almost every $e' \in M$, and therefore that $\ab'$ is close to $\ab$ (this argument is formalised in Claim~\ref{c5.1a}).

In particular, we find that $\ab'$ is not constant on $V$. So we may choose some $p \in [k-1]$ so that $u_{p+1} - u_p \gg u_k - u_{p+1}$, that is, the gap between the values of $\ab'$ on $V_{p+1}$ and $V_p$ is much greater than the gaps between the values of $\ab'$ on $V_{p+1}, \dots, V_k$. This is the point where we use our assumption that $J$ is far from a space barrier: taking $S = \bigcup_{j > p} V_j$, we can find an edge $e^* \in J_{k-p+1}[S]$. We greedily extend $e^*$ to an edge $e \in J_k$ as described earlier, giving $\ab \cdot \chi(e) \geq (k-p+1) u_{p+1} + \sum_{j \in [p-1]} u_j$, which by choice of $p$ is significantly greater than $U$. For any $e' \in M$ with $\ab \cdot \chi(e') \approx U$ this gives $\ab \cdot \chi(e) > \ab \cdot \chi(e')$, contradicting the choice of $\ab$ for $\xb = \chi(e) - \chi(e') \in X$. 

We now give the full details of the proof.

\medskip \nib{Proof of Lemma~\ref{gen_trans_1_partite}.} 
We start by applying Lemma~\ref{separating_plane_through_zero} to $X = X(J,M) \sub \Z^{rn} \cap B^{rn}(\0, 2k)$. This gives $\Pi^X_0 \cap B^{rn}(\0, \delta) \sub CH(X)$, where $F^X_0 = CH(X) \cap \Pi^X_0$ is the minimum face of $CH(X)$ containing $\0$. We can write $\Pi^X_0 = \{ \xb \in \Z^{rn}: \ab \cdot \xb = 0 \text{ for all } \ab \in A\}$ for some $A \sub \R^{rn}$ such that $\ab \cdot \xb \leq 0$ for every $\xb \in X$ and $\ab \in A$. To prove that $B^{rn}(\0, \delta) \cap \Pi_{\Part'} \sub CH(X)$, it suffices to show that $A$ is contained in the subspace generated by $\{\chi(V_i) : i \in [r]\}$. Suppose for a contradiction that some $\ab \in A$ is not in $\langle \{\chi(V_i) : i \in [r]\}\rangle$. Then there is some $i' \in [r]$ such that $\ab$ is not constant on $V_{i'}$.

We label $V$ so that for each $i \in [r]$, the vertices of $V_i$ are labelled $v_{i, 1}, \dots, v_{i,n}$, the corresponding coordinates of $\ab$ are $a_{i,1}, \dots, a_{i,n}$, and $$a_{i,1} \leq a_{i,2} \leq \dots \leq a_{i,n}.$$
It is convenient to assume that we have a strict inequality $\delta^F_i(J) > \left(\frac{k-i}{k} - \alpha \right)n$ for $i \in [k-1]$. This can be achieved by replacing $\alpha$ with a slightly smaller value; we can also assume that $\alpha n \in \N$.
Next we partition each $V_i$ into $k$ parts $V_{i,1}, \dots, V_{i,k}$ as follows. For $i \in [r]$ we let
\begin{align*}
V_{i, j} =
\begin{cases}
\{v_{i,\ell} : 0 < \ell \leq n/k - \alpha n\} & \textrm{ if } j = 1,\\
\{v_{i,\ell} : (j-1)n/k - \alpha n < \ell \leq jn/k - \alpha n\} & \textrm{ if } 2 \leq j \leq k-1,\\
\{v_{i,\ell} : (k-1)n/k - \alpha n < \ell \leq n\} & \textrm{ if } j = k.
\end{cases}
\end{align*}
Then for each $i \in [r]$ we have
\begin{align}\label{eq:sizesofVi}
|V_{i, j}| =
\begin{cases}
\frac{n}{k} - \alpha n & \textrm{ if } j = 1,\\
\frac{n}{k} & \textrm{ if } 2 \leq j \leq k-1,\\
\frac{n}{k} + \alpha n & \textrm{ if } j =k.
\end{cases}
\end{align}
We use these partitions to define a simpler vector $\ab'$ that will be a useful approximation to $\ab$. For $i \in [r]$ and $j \in [n]$ we let $p(j) \in [k]$ be such that $v_{i,j} \in V_{i,p(j)}$. For $i \in [r]$ and $\ell \in [k]$ we let $u_{i,\ell} = \max\{ a_{i,j} : p(j) = \ell \}$. Then we define $\ab' \in \Z^{rn}$ by
\[ \ab' := (a'_{i,j} : i \in [r], j \in [n]), \textrm{ where } a'_{i,j} = u_{i,p(j)}.\]
Note that $a'_{i,j} \geq a_{i,j}$ for every $i \in [r]$ and $j \in [n]$. We can assume without loss of generality that $u_{i,k}-u_{i,1}$ is maximised over $i \in [r]$ when $i=1$; then we define
\[d = u_{1,k}-u_{1,1}, \quad U = \frac{1}{r}\sum_{i \in [r],j \in [k]} u_{i,j},
\textrm{ and } U_f = \sum_{j \in [k]} u_{f(j),j} \textrm{ for } f \in F.\]
Note that uniformity of $F$ implies that $U$ is the average of $U_f$ over $f \in F$. We will prove the following claim. \medskip

\begin{claim}\label{c5.1a} We have the following properties.
\begin{enumerate}[(i)]
\item For any $e' \in M$ and $f \in F$ we have $\ab' \cdot \chi(e') \ge \ab \cdot \chi(e') \ge U_f$ and $\ab \cdot \chi(e') \ge U$,
\item $d>0$ and $\sum_{e' \in M} (\ab' \cdot \chi(e') - U) \le rd\alpha n$,
\item There is some $e' \in M$ such that $\ab \cdot \chi(e') \leq U_f + D_F kd\alpha$ for every $f \in F$,
\item For any $i, i' \in [r]$ and $j, j' \in [k]$ we have
$|(u_{i, j} - u_{i, j'}) - (u_{i',j} - u_{i',j'})| \leq D_F rkd\alpha$.
\end{enumerate}
\end{claim}

To see property (i), note that the first inequality follows from $a'_{i,j} \geq a_{i,j}$. For the second, we use the minimum $f$-degree sequence of $J$ to greedily form an edge $e =  \{v_{f(k), \ell_k}, \dots, v_{f(1), \ell_1}\} \in J_k$, starting with $\ell_k = n$, then choosing $\ell_j > jn/k - \alpha n$ for each $j = k-1,k-2,\dots,1$. Then $p(\ell_j) \ge j+1$ for $j \in [k-1]$, so $\ab \cdot \chi(e) \geq U_f$. Now the inequality follows from the definition of $A$, which gives $\ab \cdot (\chi(e')-\chi(e)) \ge 0$. We noted above that $U$ can be expressed as an average of some $U_f$'s, so $\ab \cdot \chi(e') \ge U$. This proves property (i).

To see that $d>0$, suppose for a contradiction that $u_{i,k}=u_{i,1}$ for $i \in [r]$. Recall that there is some $i' \in [r]$ such that $\ab$ is not constant on $V_{i'}$. Let $e' \in M$ be the edge of $M$ containing $v_{i',1}$. Since $(J,M)$ is $\Part' F$-partite we can write $e' = \{v_{f(i),n_i}\}_{i \in [k]}$ for some $f \in F$ and $n_i \in [n]$ for $i \in [k]$, where without loss of generality $v_{f(1),n_1}=v_{i',1}$. Then $$\ab \cdot \chi(e') \leq a_{i',1} + \sum_{i=2}^k a_{f(i),n_i} < \sum_{i=1}^k a_{f(i),n}  = U_f.$$ However, this contradicts property (i), so $d>0$. Next note that $\sum_{e' \in M} \chi(e') = \1$ and $|M|=rn/k$, so by (\ref{eq:sizesofVi}) and the definition of~$d$ we have
\begin{align*}
\sum_{e' \in M} (\ab' \cdot \chi(e') - U) = \ab' \cdot \1 - |M|U =& \sum_{i \in [r],j \in [k]} (|V_i| - n/k)u_{i,j} \\=& \sum_{i \in [r]} (u_{i,k}-u_{i,1})\alpha n \leq rd \alpha n,
\end{align*} so we have proved property (ii).

Now we may choose some $e' \in M$ with $0 \le \ab' \cdot \chi(e') - U \le rd\alpha n/|M| = kd\alpha$. Since $U = |F|^{-1}\sum_{f \in F} U_f$, $\ab' \cdot \chi(e') \ge \ab \cdot \chi(e') \ge U_f$ for any $f \in F$ and $|F| \le D_F$, we have $0 \leq \ab' \cdot \chi(e') - U_f \leq D_F kd\alpha$; this implies (iii), using $\ab' \cdot \chi(e') \ge \ab \cdot \chi(e')$. It also implies $|U_f - U_{f'}| \le D_F kd\alpha$ for any $f, f' \in F$. Now, by definition of $G_F$ and invariance of $F$ under $\Permute_k$, for any $i, i' \in [r]$ and $j, j' \in [k]$ such that $ii' \in G_F$, we can choose $f\in F$ with $f(j)=i$ and $f(j')=i'$, and let $f' \in F$ be obtained from $f$ by transposing the values on $j$ and $j'$. Then $|U_f - U_{f'}| = |(u_{i, j} - u_{i, j'}) - (u_{i',j} - u_{i',j'})| \le D_F kd\alpha$. Since $G_F$ is connected we obtain statement (iv), so we have proved Claim~\ref{c5.1a}.

\medskip

To continue the proof of Lemma~\ref{gen_trans_1_partite}, we say that an edge $e' \in M$ is \emph{good} if $\ab' \cdot \chi(e') \leq U + d\sqrt{\alpha}$. We say that a vertex $v \in V(J)$ is \emph{good} if it lies in a good edge of $M$. Note that if $v_{i,j}$ is a good vertex and $v_{i,j} \in e' \in M$, then Claim~\ref{c5.1a} (i) gives $0 \le (\ab'-\ab) \cdot \chi(e') \leq d\sqrt{\alpha}$, so $a'_{i,j} - a_{i,j} < d \sqrt{\alpha}$. Writing $B$ for the set of bad edges of $M$, by Claim~\ref{c5.1a} (i) and (ii) we have $$|B| d\sqrt{\alpha} < \sum_{e' \in B} (\ab' \cdot \chi(e') - U) \le \sum_{e' \in M} (\ab' \cdot \chi(e') - U) \le rd\alpha n,$$ so $|B| < r\sqrt{\alpha} n$. Thus all but at most $rk \sqrt{\alpha}n$ vertices $v \in V(J)$ are good. Next we need another claim.

\begin{claim}\label{c5.1b}
There is some $p \in [k-1]$ and an edge $e^* \in J_{k-p+1}$ in which all vertices are good, such that $\ab' \cdot \chi(e^*) \ge (k+1)d\sqrt{\alpha} + \sum_{p \le j \le k} u_{f(j),j}$, where $f \in F$ is such that $e^* = \{v_p,\dots,v_k\}$ with  $v_j \in V_{f(j)}$ for $p \le j \le k$.
\end{claim}

To prove this claim, we start by choosing $p$ so that the gap $u_{1,p+1} - u_{1,p}$ is considerably larger than $u_{1,k} - u_{1,p+1}$. More precisely, we require \[ u_{1,p+1} - u_{1,p} > k (u_{1,k} - u_{1,p+1}) + (k+2)d \sqrt{\alpha}.\]
Suppose for a contradiction that this is not possible. Then for every $p \in [k-1]$ we have 
$$u_{1,k} - u_{1,p} = (u_{1,k} - u_{1,p+1}) + (u_{1,p+1} - u_{1,p}) \le (k+2)(u_{1,k} - u_{1,p+1} + d \sqrt{\alpha}).$$
Iterating this inequality starting from $u_{1,k} - u_{1,k-1} \le (k+2)d \sqrt{\alpha}$ we obtain the bound $u_{1,k} - u_{1,p} \le (k+3)^{k-p} d \sqrt{\alpha}$. However, for $p=1$ we obtain $d = u_{1,k} - u_{1,1} \le (k+3)^{k-1} d \sqrt{\alpha}$, which is a contradiction for $\alpha \ll 1/k$. Thus we can choose $p$ as required. 

Now consider $S = \bigcup_{i \in [r], p < j \leq k} V_{i,j}$. Note that $|S \cap V_i|  = (k-p)n/k + \alpha n$ for each $i \in [r]$. Thus we can apply the assumption that there is no space barrier, which gives at least $\beta n^{k-p+1}$ edges in $J_{k-p+1}[S]$.  At most $rk \sqrt{\alpha} n^{k-p+1}$ of these edges contain a vertex which is not good, so we may choose an edge $e^* \in J_{k-p+1}$ whose vertices are all good and lie in $S$. Since $J$ is $\Part' F$-partite, we can choose $f \in F$ so that $e^* = \{v_{f(p), \ell_p}, \dots, v_{f(k), \ell_k}\}$ for some $\ell_p, \dots, \ell_k$, and since $e^* \sub S$ we have $\ab' \cdot \chi(e^*) \ge \sum_{p \le j \le k} u_{f(j),p+1}$. Then by Claim~\ref{c5.1a} (iv) we can estimate \[\ab' \cdot \chi(e^*) - \sum_{p \le j \le k} u_{f(j),j} \ge \sum_{p \leq j \leq k} (u_{1,p+1} - u_{1,j}) - D_Frk^2d\alpha.\]
By choice of $p$ we have 
$$\sum_{p \leq j \leq k} (u_{1,p+1} - u_{1,j}) \ge (u_{1,p+1} - u_{1,p}) - (k-p)(u_{1,k} - u_{1,p+1}) > (k+2)d \sqrt{\alpha}.$$
The required bound on $\ab' \cdot \chi(e^*)$ follows, so this proves Claim~\ref{c5.1b}.

\medskip

To finish the proof of Lemma~\ref{gen_trans_1_partite}, using the minimum degree sequence of $J$ as above, we greedily extend $e^*$ to an edge $e =  \{v_{f(k), \ell_k}, \dots, v_{f(1), \ell_1}\} \in J_k$ with $\ell_j > jn/k - \alpha n$ for $j=p-1,p-2,\dots,1$. This gives
$\ab \cdot \chi(e \sm e^*) \geq \sum_{j \in [p-1]} u_{f(j),j}$. Also, since every vertex of $e^*$ is good we have $(\ab' - \ab) \cdot \chi(e^*) \leq kd\sqrt{\alpha}$. It follows that
\begin{align*} \ab \cdot \chi(e) & = \ab' \cdot \chi(e^*) - (\ab' - \ab) \cdot \chi(e^*)  + \ab \cdot \chi(e \sm e^*) \\
& \ge (k+1)d\sqrt{\alpha} + \sum_{p \le j \le k} u_{f(j),j} - kd\sqrt{\alpha} +  \sum_{j \in [p-1]} u_{f(j),j}
=  U_f + d\sqrt{\alpha}. \end{align*}
On the other hand, for any $e' \in M$ we have $\ab \cdot \chi(e) \le \ab \cdot \chi(e')$ by definition of $A$,
and by Claim~\ref{c5.1a} (iii) we can choose $e'$ so that $\ab \cdot \chi(e') \le U_f + D_Fkd\alpha$. Thus we have a contradiction to our original assumption that $A \not\sub \langle \{\chi(V_i) : i \in [r]\}\rangle$, which proves the Lemma.
\endproof

Next we need to address a technical complication alluded to earlier, which is that we cannot satisfy the assumption $\delta \ll 1/n$ in Lemma~\ref{gen_trans_1_partite} by working directly with $J$. Thus we need the next lemma, which will be proved by taking a random selection of edges in $M$, that allows us to reduce to a small matched subsystem $(J', M')$ with similar properties to $(J, M)$. For this we need the following definitions.

\begin{defn}
Suppose that $\Part$ partitions a vertex set $V$ into parts $V_1, \dots, V_r$, that $(J, M)$ is a matched $k$-graph or $k$-system on $V$, and that $F$ is is a $(k,r)$-uniform allocation. We say that \emph{$M$ $\alpha$-represents $F$} if for any $\ib, \ib'$ in $I(F)$, letting $N, N'$ denote the number of edges $e' \in M$ with index $\ib, \ib'$ respectively, we have $N' \ge (1-\alpha)N$. That is, each $\ib \in I(F)$ is represented by approximately the same number of edges of $M$.

In particular, if $M$ $0$-represents $F$ then the number of edges $e' \in M$ with $\ib(e') = \ib$ is the same for all $\ib \in I(F)$; in this case we say that $M$ is \emph{$F$-balanced}.
\end{defn}

 Note that when $(J, M)$ is $\Part F$-partite, the fact that $M$ $\alpha$-represents $F$ implies that $(J,M)$ is $\Part$-proper (recall that this was a condition needed for Lemma~\ref{make_rectuple_spanning_partite}).

\begin{lemma} \label{RANDOMEDGESELECTION}
Suppose that $1/n \ll 1/n' \ll \alpha' \ll \alpha \ll \beta' \ll \beta \ll 1/D_F, 1/r, 1/k$.
Let $\Part$ partition a set $V$ into $r$ parts $V_1, \dots, V_r$ each of size $n$. Suppose $F$ is a $(k,r)$-uniform allocation with $|F| \leq D_F$ and $(J, M)$ is a $\Part F$-partite matched $k$-system on $V$ such that $M$ $\alpha'$-represents $F$. Suppose also that
\begin{enumerate}[(i)]
\item $\delta^F(J) \geq \left(n, \left(\frac{k-1}{k} - \alpha \right)n, \left(\frac{k-2}{k} - \alpha \right) n, \dots, \left(\frac{1}{k} - \alpha \right) n\right)$, and
\item for any $p \in [k-1]$ and sets $S_i \sub V_i$ such that $|S_i| = \lfloor pn/k \rfloor$  for each $i \in [r]$ there are at least $\beta n^{p+1}$ edges in $J_{p+1}[S]$, where $S := \bigcup_{i \in [r]} S_i$.
\end{enumerate}
Then for any $a \in [r]$ and $u,v \in V_a$, there exists a set $M' \sub M$ such that, writing $V' = \bigcup_{e \in M'} e$, $J' = J[V']$, and $\Part'$ for the partition of $V'$ into $V'_i := V_i \cap V'$, $i \in [r]$, we have $u,v \in V'_a$, $|V_1'| = \dots = |V_r'| = n'$,
and $(J', M')$ is a $\Part' F$-partite matched $k$-system on $V'$ such that
\begin{enumerate}[(i)]
\item $\delta^F(J') \geq \left(n', \left(\frac{k-1}{k} - 2\alpha \right)n', \left(\frac{k-2}{k} - 2\alpha \right) n', \dots, \left(\frac{1}{k} - 2\alpha \right) n'\right)$, and
\item for any $p \in [k-1]$ and sets $S'_i \sub V'_i$ such that $|S'_i| = \lfloor pn'/k \rfloor$ for each $i \in [r]$ there are at least $\beta' (n')^{p+1}$ edges in $J'_{p+1}[S']$, where $S' := \bigcup_{i \in [r]} S'_i$.
\end{enumerate}
\end{lemma}

The proof of this lemma requires some regularity theory, so we postpone it to the next chapter. Now we can combine Lemma~\ref{gen_trans_1_partite} and Lemma~\ref{RANDOMEDGESELECTION}, to obtain the following lemma, which shows that our assumptions guarantee irreducibility with the correct dependence of parameters. Note that the condition on $M$ in the final statement is automatically satisfied if $J$ is a complex.

\begin{lemma} \label{gen_trans_2_partite}
Suppose that $1/n \ll 1/\ell \ll 1/B, 1/C \ll \alpha' \ll \alpha \ll \beta \ll 1/D_F, 1/r, 1/k$. Let $V$ be a set of $rn$ vertices, and let $\Part'$ partition $V$ into $r$ parts $V_1, \dots, V_r$ each of size~$n$. Suppose $F$ is a $(k,r)$-uniform connected allocation with $|F| \leq D_F$, and $(J, M)$ is a $\Part' F$-partite matched $k$-system on $V$ such that $M$ $\alpha'$-represents $F$. Suppose also that
\begin{enumerate}[(i)]
\item $\delta^F(J) \geq \left(n, \left(\frac{k-1}{k} - \alpha \right)n, \left(\frac{k-2}{k} - \alpha \right) n, \dots, \left(\frac{1}{k} - \alpha \right) n\right)$, and
\item for any $p \in [k-1]$ and sets $S_i \sub V_i$ such that $|S_i| = \lfloor pn/k \rfloor$ for each $i \in [r]$ there are at least $\beta n^{p+1}$ edges in $J_{p+1}[S]$, where $S := \bigcup_{i \in [r]} S_i$.
\end{enumerate}
Then $(J_k, M)$ is $(B, C)$-irreducible with respect to $\Part'$.

If, in addition, for any $v \in e' \in M$ we have $e' \sm \{v\} \in J$, then there is a partition $\Part$ refining $\Part'$ such that  $|U| \ge n/k - 2\alpha n$ and $D_\ell(J,M)[U]$ is complete for every $U \in \Part$.
\end{lemma}

\proof
Choose $n',\delta$ with $1/B,1/C \ll \delta \ll 1/n' \ll \alpha'$ and $kr \mid n'$, fix any $i \in [r]$,  $u, v \in V_i$, and
let $M'$, $V'$, $J'$ be given by applying Lemma~\ref{RANDOMEDGESELECTION}. Let $\Part''$ be the partition of $V'$ into $V' \cap V_1, \dots, V' \cap V_r$ and $X' = X(J',M') \sub \R^{n'}$. Applying Lemma~\ref{gen_trans_1_partite}, we have $B^{rn'}(\0,\delta) \cap \Pi_{\Part''} \sub CH(X')$. Then by Lemma~\ref{ballgivesgentransferrals} we deduce that $(J'_k, M')$ is $(B, C)$-irreducible with respect to $\Part''$. Since $i \in [r]$ and $u,v \in V_i$ were arbitrary, it follows that $(J_k, M)$ is $(B, C)$-irreducible with respect to $\Part'$. For the final statement we claim that $\delta^+(D_1(J_k, M)[V_q]) \geq \delta^F_{k-1}(J) \geq n/k - \alpha n$ for each $q \in [r]$. Indeed, consider any $v \in V_q$, and let $e'$ be the edge of $M$ containing $v$. Then we can write $e' = \{v_1,\dots,v_k\}$ where $v_k=v$ and $v_i \in V_{f(i)}$ for $i \in [k]$ for some $f \in F$. By assumption $e' \sm \{v\} \in J$, so there are at least $\delta^f_{k-1}(J)$ vertices $u \in V_q$ such that $\{u\} \cup e' \sm \{v\} \in J$. For each such $u$, $(\{\{u\} \cup e' \sm \{v\}\}, \{e'\})$ is a simple $(u,v)$-transferral in $(J_k, M)$ of size one, so this proves the claim. The conclusion then follows from Lemma~\ref{complete_on_receiving}.
\endproof

Now we can formulate our main transferral lemma, which is that the minimum degree sequence assumption, with no divisibility or space barrier, implies the existence of all required simple transferrals. The proof is immediate from our previous lemmas.

\begin{lemma} \label{main-reduction-partiteF}
Suppose that $1/n \ll 1/\ell \ll \alpha' \ll \alpha \ll \beta \ll 1/D_F, 1/r, 1/k$. Let $V$ be a set of $rn$ vertices, and let $\Part'$ be a partition of $V$ into parts $V_1, \dots, V_r$ each of size $n$. Suppose $F$ is a $(k,r)$-uniform connected allocation with $|F| \leq D_F$, $(J, M)$ is a $\Part' F$-partite matched $k$-system on $V$ such that $M$ $\alpha'$-represents $F$, and for any $v \in e' \in M$ we have $e' \sm \{v\} \in J$. Suppose also that 
\begin{enumerate}[(i)]
\item $\delta^F(J) \geq \left(n, \left(\frac{k-1}{k} - \alpha \right)n, \left(\frac{k-2}{k} - \alpha \right) n, \dots, \left(\frac{1}{k} - \alpha \right) n\right)$,
\item for any $p \in [k-1]$ and sets $S_i \sub V_i$ such that $|S_i| = \lfloor pn/k \rfloor$  for each $i \in [r]$ there are at least $\beta n^{p+1}$ edges in $J_{p+1}[S]$, where $S := \bigcup_{i \in [r]} S_i$, and
\item $L_\Part(J_k)$ is complete with respect to $\Part'$ for any partition~$\Part$ of~$V(J)$ which refines~$\Part'$ and whose parts each have size at least $n/k - 2\alpha n$.
\end{enumerate}
Then $D_\ell(J_k, M)[V_j]$ is complete for each $j \in [r]$. That is, for any $j \in [k]$ and $u, v \in V_j$ the matched $k$-graph $(J_k, M)$ contains a simple $(u,v)$-transferral of size at most $\ell$.
\end{lemma}

\proof
Let $\ell',B$ and $C$ satisfy $1/\ell \ll 1/\ell' \ll 1/C,1/B \ll \alpha$. By Lemma~\ref{gen_trans_2_partite}, $(J_k, M)$ is $(B,C)$-irreducible with respect to $\Part'$, and there is a partition $\Part$ refining $\Part'$ such that  $|U| \ge n/k - 2\alpha n$ and $D_{\ell'}(J,M)[U]$ is complete for every $U \in \Part$. Then $L_\Part(J_k)$ is complete with respect to $\Part'$ by assumption, so Lemma~\ref{make_rectuple_spanning_partite} applies to give the required result.
\endproof

\chapter{Hypergraph Regularity Theory} \label{sec:regularity}

This chapter contains the technical background that we need when applying the method of hypergraph regularity. Most of the machinery will be quoted from previous work, although we also give some definitions and short lemmas that are adapted to our applications. We start in the first section by introducing our notation, defining hypergraph regularity, and stating the `regular restriction lemma'. The second section states the Regular Approximation Lemma of R\"odl and Schacht. In the third section we state the Hypergraph Blow-up Lemma, and prove a short accompanying result that is adapted to finding perfect matchings. In the fourth section we define reduced $k$-systems, and develop their theory, including the properties that they inherit degree sequences, and that their edges represent $k$-graphs to which the Hypergraph Blow-up Lemma applies. The final section contains the proof of Lemma~\ref{RANDOMEDGESELECTION}, using the simpler theory of `weak regularity' (as opposed to the `strong regularity' in the rest of the chapter).

\section{Hypergraph regularity.}

Let $X$ be a set of vertices, and let $\Qart$ be a partition of $X$ into $r$ parts $X_1, \dots, X_r$.  The \emph{index} $i_\Qart(S)$ of a set $S \sub X$ is the multiset in $[r]$ where the multiplicity of $j \in [r]$ is $\ib_\Qart(S)_j = |S \cap X_j|$; we generally write $i(S)=i_\Qart(S)$ when $\Qart$ is clear from the context. Recall that a set $S \sub X$ is $\Qart$-partite if it has at most one vertex in each part of $\Qart$; for such sets $i(S)$ is a set. For any $A \sub [r]$ we write $X_A$ for the set $\bigcup_{i \in A} X_i$, and $K_A(X)$ for the complete $\Qart$-partite $|A|$-graph on $X_A$, whose edges are all $\Qart$-partite sets $S \sub X$ with $i(S) = A$. For a $\Qart$-partite set $S \sub X$ and $A \sub i(S)$, we write $S_A = S \cap X_A$. If $H$ is a $\Qart$-partite $k$-graph or $k$-system on $X$, then for any $A \sub [r]$ we define $H_A := H \cap K_A(X)$. Equivalently, $H_A$ consists of all edges of $H$ with index $A$, which we regard as an $|A|$-graph on vertex set $X_A$. If additionally $H$ is a $k$-complex, then we write $H_{A^\leq}$ for the $k$-complex on $X_A$ with edge set $\bigcup_{B \sub A} H_B$, and $H_{A^<}$ for the $k$-complex on $X_A$ with edge set $\bigcup_{B \subset A} H_B$. Similarly, we write $H_{\ib}$ for the set of edges in $H$ with index vector $\ib$. This is a $|\ib|$-graph, where $|\ib| := \sum_{j \in [r]} i_j$.

To understand the definition of hypergraph regularity it is helpful to start with the case of graphs. If $G$ is a bipartite graph on vertex classes $U$ and $V$,  we say that $G$ is \emph{$\eps$-regular} if for any $U' \sub U$ and $V' \sub V$ with $|U'| > \eps |U|$ and $|V'| > \eps |V|$ we have $d(G[U' \cup V']) = d(G) \pm \eps$.  Similarly, for a $k$-complex $G$, an informal statement of regularity is that the restriction of $G$ to any large subcomplex of the `lower levels' of $G$ has similar densities to $G$. Formally, let $\Qart$ partition a set $X$ into parts $X_1, \dots, X_r$, and $G$ be an $\Qart$-partite $k$-complex on $X$. We denote by $G^*_A$ the $|A|$-graph on $X_A$ whose edges are all those $S \in K_A(X)$ such that $S' \in G$ for every strict subset $S' \subset S$. So $G^*_A$ consists of all sets which could be edges of $G_A$, in the sense that they are supported by edges at `lower levels'. The \emph{relative density of $G$ at $A$} is
$$d_A(G) := \frac{|G_A|}{|G^*_A|};$$
this is the proportion of `possible edges' (given $G_{A<}$) that are in fact edges of $G_A$.
(It should not be confused with the \emph{absolute density of $G$ at $A$}, which is $d(G_A) := |G_A|/|K_A(X)|$.)

For any $A \in \binom{[r]}{\leq k}$, we say that $G_A$ is \emph{$\eps$-regular} if for any subcomplex $H\sub G_{A^<}$ with $|H^*_A| \geq \eps |G^*_A|$ we have
$$\frac{|G_A \cap H_A^*|}{|H_A^*|} = d_A(G) \pm \eps.$$
We say $G$ is \emph{$\eps$-regular} if $G_A$ is $\eps$-regular for every $A \in \binom{[r]}{\leq k}$.

The following lemma states that the restriction of any regular and dense $k$-partite $k$-complex to large subsets of its vertex classes is also regular and dense (it is a weakened version of \cite[Theorem~6.18]{K}).

\begin{lemma} (Regular restriction) \label{regularrestriction}
Suppose that $1/n \ll \eps \ll c, 1/k$. Let $\Qart$ partition a set $X$ into $X_1, \dots, X_k$, and $G$ be an $\eps$-regular $\Qart$-partite $k$-complex on $X$ with $d(G) \ge c$. Then for any subsets $X_1' \sub X_1$, \dots, $X_k' \sub X_k$ each of size at least $\eps^{1/2k} n$, the restriction $G' = G[X_1' \cup \dots \cup X_k']$ is $\sqrt{\eps}$-regular with $d(G') \ge d(G)/2$ and $d_{[k]}(G') \ge d_{[k]}(G)/2$.
\end{lemma}

\section{The Regular Approximation Lemma.}

Roughly speaking, hypergraph regularity theory shows that an arbitrary $k$-graph can be split into pieces, each of which forms the `top level' of a regular $k$-complex. To describe the splitting, we require the following definitions. Let $\Qart$ partition a set $X$ into $r$ parts $X_1, \dots, X_r$. A \emph{$\Qart$-partition $k$-system} on $X$ consists of a partition $\C_A$ of $K_A(X)$ for each $A \in \binom{[r]}{\leq k}$. We refer to the parts of each $\C_A$ as \emph{cells}, and to the cells of $\C_{\{i\}}$ for each $i \in [k]$ as the \emph{clusters} of $P$. Observe that the clusters of $P$ form a partition of $X$ which refines $\Qart$. For any $\Qart$-partite set $S \sub X$ with $|S| \leq k$, we write $C_S$ for the set of all edges of $K_{i(S)}(X)$ lying in the same cell of $P$ as~$S$. We write $C_{S^\leq}$ for the $\Qart$-partite $k$-system with vertex set $X$ and edge set $\bigcup_{S' \sub S} C_{S'}$.

Let $P$ be a $\Qart$-partition $k$-system on $X$. We say that $P$ is a \emph{$\Qart$-partition $k$-complex} if $P$ satisfies the additional condition that whenever edges $S, S' \in K_A(X)$ lie in the same cell of $\C_A$, the edges $S_B, S_B'$ of $K_B(X)$ lie in the same cell of $\C_B$ for any $B \sub A$. Note that then $C_{S^\le}$ is a $k$-complex.  We say that $P$ is \emph{vertex-equitable} if every cluster of $P$ has the same size. We say that $P$ is \emph{$a$-bounded} if for every $A \in \binom{[r]}{\leq k}$ the partition $\C_A$ partitions $K_A(X)$ into at most $a$ cells. We say that $P$ is \emph{$\eps$-regular} if $C_{S^\leq}$ is $\eps$-regular for every $\Qart$-partite $S \sub X$ with $|S| \leq k$.

Now instead let $P$ be a $\Qart$-partition $(k-1)$-complex on $X$. Then $P$ naturally induces a $\Qart$-partition $k$-complex $P'$ on $X$. Indeed, for any $A \in \binom{[r]}{k}$ we say that $S, S' \in K_A(X)$ are \emph{weakly equivalent} if for any strict subset $B \subset A$ we have that $S_B$ and $S_B'$ lie in the same cell of $\C_B$; this forms an equivalence relation on $K_A(X)$. Then for each $A \in \binom{[r]}{\leq k-1}$ the partition $\C'_A$ of $P'$ is identical to the partition $\C_A$ of $P$, and for each $A \in \binom{[r]}{k}$ the partition $\C'_A$ of $P'$ has the equivalence classes of the weak equivalence relation as its cells. We refer to $P'$ as the $\Qart$-partition $k$-complex generated from $P$ by weak equivalence. Note that if $P$ is $a$-bounded then $P'$ is $a^k$-bounded, as for each $A \in \binom{[r]}{k}$ we have that $K_A(X)$ is divided into at most $a^k$ cells by weak equivalence. Now let $G$ be a $\Qart$-partite $k$-graph on $X$. We denote by $G[P]$ the $\Qart$-partition $k$-complex formed by using weak equivalence to refine the partition $\{G_A, K_A(X) \sm G_A\}$ of $K_A(X)$ for each $A \in \binom{[r]}{k}$, i.e.\ two edges of $G_A$ are in the same cell if they are weakly equivalent, and similarly for two $k$-sets in $K_A(X) \sm G_A$. Together with $P$, this yields a partition $k$-complex which we denote by $G[P]$. If $G[P]$ is $\eps$-regular then we say that $G$ is \emph{perfectly $\eps$-regular with respect to $P$}.

We use the following form of hypergraph regularity due to R\"odl and Schacht~\cite{RS} (it is a slight reformulation of their result). It states that any given $k$-graph $H$ can be approximated by another $k$-graph $G$ that is regular with respect to some partition $(k-1)$-complex $P$. It is convenient to consider the $k$-graph $G$, but we take care not to use any edges in $G \sm H$, to ensure that every edge we use actually lies in the original $k$-graph $H$. There are various other forms of the regularity lemma for $k$-graphs which give information on $H$ itself, but these do not have the hierarchy of densities necessary for the application of the blow-up lemma (see \cite{K} for discussion of this point). In the setting of the theorem, we say that $G$ and $H$ are \emph{$\nu$-close} (with respect to $\Qart$) if $|G_A \bigtriangleup H_A| < \nu |K_A(V)|$ for every $A \in \binom{[r]}{k}$.

\begin{theo} (Regular Approximation Lemma) \label{eq-partition}
Suppose $1/n \ll \eps \ll 1/a \ll \nu, 1/r, 1/k$ and $a!r|n$. Let $V$ be a set of $n$ vertices, $\Qart$ be an balanced partition of $V$ into $r$ parts, and $H$ be a $\Qart$-partite $k$-graph on $V$. Then there is an $a$-bounded $\eps$-regular vertex-equitable $\Qart$-partition $(k-1)$-complex $P$ on $V$, and an $\Qart$-partite $k$-graph $G$ on $V$, such that $G$ is $\nu$-close to $H$ and perfectly $\eps$-regular with respect to $P$.
\end{theo}

\section{The hypergraph blowup lemma.}

While hypergraph regularity theory is a relatively recent development, still more recent is the hypergraph blow-up lemma due to Keevash~\cite{K}, which makes it possible to apply hypergraph regularity theory to embeddings of spanning subcomplexes. Indeed, it is similar to the blow-up lemma for graphs, insomuch as it states that by deleting a small number of vertices from a regular $k$-complex we may obtain a \emph{super-regular} complex, in which we can find any spanning subcomplex of bounded maximum degree. However, unlike the graph case, the definition of super-regularity for complexes is extremely technical, so it is more convenient to work with a formulation using robustly universal complexes. In essence, a complex $J'$ is robustly universal if even after the deletion of many vertices (with certain conditions), we may find any spanning subcomplex of bounded degree within the complex $J$ that remains. The formal definition is as follows (we have simplified it by removing the option of `restricted positions', which are not required in this paper).

\begin{defn}[Robustly universal]
Suppose that $V'$ is a set of vertices, $\Qart$ is a partition of $V'$ into $k$ parts $V'_1 \cup \dots \cup V'_k$, and $J'$ is a $\Qart$-partite $k$-complex on $V'$ with $J'_{\{i\}} = V'_i$ for each $i \in [k]$. Then we say that $J'$ is \emph{$c$-robustly
$D$-universal} if whenever
\begin{itemize}
\item[(i)] $V_j \sub V'_j$ are sets with $|V_j| \ge c|V'_j|$ for each $j \in [k]$, such that writing $V = \bigcup_{j \in [k]} V_j$ and $J=J'[V]$ we have $|J_k(v)| \ge c|J'_k(v)|$ for any $j \in [k]$ and $v \in V_j$, and
\item[(ii)] $L$ is a $k$-partite $k$-complex of maximum vertex degree at most $D$ whose vertex classes $U_1, \dots, U_k$ satisfy $|U_j|=|V_j|$ for each $j \in [k]$,
\end{itemize}
\noindent then $J$ contains a copy of $L$, in which for each $j \in [k]$ the vertices of $U_j$ correspond to the vertices of $V_j$.
\end{defn}

The following version of the hypergraph blow-up lemma states that we may obtain a robustly universal complex from a regular complex by deleting a small number of vertices (it is a special case of \cite[Theorem 6.32]{K}). After applying Theorem~\ref{eq-partition}, we regard $Z = G \sm H$ as the `forbidden' edges of $G$; so with this choice of $Z$ in Theorem~\ref{blowup}, the output $G' \sm Z'$ is a subgraph of $H$.

\begin{theo} (Blow-up Lemma) \label{blowup}
Suppose $1/n \ll \eps \ll d^* \ll d_a \ll \theta \ll d, c, 1/k, 1/D, 1/C$.
Let $V$ be a set of vertices, $\Qart$ be a partition of $V$ into $k$ parts $V_1, \dots, V_k$ with $n \leq |V_j| \leq Cn$ for each $j \in [k]$, and $G$ be an $\eps$-regular $\Qart$-partite $k$-complex on $V$ such that $|G_{\{j\}}| = |V_j|$ for each $j \in [k]$, $d_{[k]}(G) \geq d$ and $d(G) \geq d_a$. Suppose $Z \sub G_k$ satisfies $|Z| \leq \theta |G_k|$.
Then we can delete at most $2\theta^{1/3} |V_j|$ vertices
from each $V_j$ to obtain $V' = V_1' \cup \dots \cup V_k'$, $G' = G[V']$ and $Z' = Z[V']$ such that
\begin{itemize}
\item[(i)] $d(G') > d^*$ and $|G'(v)_k| > d^*|G'_k|/|V'_i|$ for every $v \in V'_i$, and
\item[(ii)] $G' \sm Z'$ is $c$-robustly $D$-universal.
\end{itemize}
\end{theo}

We will apply the blow-up lemma in conjunction with the following lemma for finding perfect matchings in subcomplexes of robustly universal complexes. The set $X$ forms an `ideal' for the perfect matching property, in that {\em any} extension $W$ of $X$ with parts of equal size contains a perfect matching. We will see later that a random choice of $X$ has this property with high probability.

\begin{lemma} \label{randomsplitkeepsmatching}
Let $\Qart$ partition a set $U$ into $U_1, \dots, U_k$, and $G$ be a $\Qart$-partite $k$-complex on $U$ which is $c$-robustly $2^k$-universal. Suppose we have $X_j \sub U_j$ for each $j \in [k]$ such that
\begin{itemize}
\item[(i)] $|X_j| \geq c|U_j|$ for each $j \in [k]$, and
\item[(ii)] $|G_k[X \cup \{v\}](v)| \geq c|G_k(v)|$ for any $v \in U$,
where $X := X_1 \cup \dots \cup X_k$.
\end{itemize}
Then for any sets $W_j$ with $X_j \sub W_j \sub U_j$ for $j \in [k]$ and $|W_1|= \dots =|W_k|$, writing $W = W_1 \cup \dots \cup W_k$, the $k$-complex $G[W]$ contains a perfect matching,
\end{lemma}

\proof
By (i) we have $|W_j| \geq |X_j| \geq c|U_j|$ for each $j$, and by (ii) we have $$|G_k[W](v)| \geq |G_k[X \cup \{v\}](v)| \geq c|G_k(v)|$$ for every $v \in W$. So by definition of a $c$-robustly $2^k$-universal complex, $G$ contains any $k$-partite $k$-complex on $W$ with maximum vertex degree at most $2^k$. In particular, this includes the complex obtained by the downward closure of a perfect matching in $K_{[k]}[W]$.
\endproof

\section{Reduced $k$-systems.} Now we introduce and develop the theory of reduced $k$-systems, whose role in the $k$-system setting is analogous to that of reduced graphs in the graph setting. Informally speaking, an edge of the reduced system represents a set of clusters, which is dense in the original $k$-system $J$, and sparse in the forbidden $k$-graph $Z$. While this does not contain enough information for embeddings of general $k$-graphs, it is sufficient for matchings, which is our concern here. Note that this definition includes a partition $(k-1)$-complex $P$, which will be obtained from the regularity lemma, but only the first level of this complex (the partition of $X$ into clusters) is used in the definition.

\begin{defn}[Reduced $k$-system]
Let~$\Qart$ partition a set $X$ of size $n$ into parts of equal size, and let $J$ be a $\Qart$-partite $k$-system, $Z$ a $\Qart$-partite $k$-graph, and $P$ a vertex-equitable $\Qart$-partition $(k-1)$-complex on $X$. Suppose $\nu \in \R$ and $\cb = (c_1,\dots,c_k) \in \R^k$. We define the \emph{reduced $k$-system} $R := R^{JZ}_{P\Qart}(\nu,\cb)$ as follows.

Let $V_1, \dots, V_m$ be the clusters of $P$, and let $n_1$ denote their common size. The vertex set of $R$ is $[m]$, where vertex~$i$ corresponds to cluster~$V_i$ of~$P$. Any partition $\Part$ of $X$ which is refined by the partition of $X$ into clusters of $P$ naturally induces a partition of~$[m]$, which we denote by $\Part_R$: $i$ and $j$ lie in the same part of $\Part_R$ if and only if $V_i$ and $V_j$ are subsets of the same part of $\Part$. Note that $\Qart_R$ partitions $[m]$ into parts of equal size.

The edges of $R$ are defined as $\es$, and those $S \in \binom{[m]}{j}$ for each $j \in [k]$ such that $S$ is $\Qart_R$-partite, $|J[\bigcup_{i \in S} V_i]| \geq c_jn_1^j$, and for any $S' \sub S$ of size $j'$, at most $\nu^{2^{-j'}} n_1^{j'}n^{k-j'}$ edges of~$Z$ intersect each member of $\{V_i : i \in S'\}$.
\end{defn}

Ideally, we would like the reduced $k$-system of a $k$-complex $J$ to be itself a $k$-complex. However, it does not appear to be possible to define a reduced $k$-system that both accomplishes this and inherits a minimum degree condition similar to that of $J$. Instead, we have a weaker property set out by the following result, which shows that any subset of an edge of a reduced $k$-system of $J$ is an edge of another reduced $k$-system of $J$, where the latter has weaker density parameters.

\begin{lemma} \label{edgescloseddownwards}
Let~$\Qart$ partition a set $X$ into parts of equal size, $J$ be a $\Qart$-partite $k$-system, $Z$ a $\Qart$-partite $k$-graph, and $P$ a vertex-equitable $\Qart$-partition $(k-1)$-complex on $X$, with clusters of size $n_1$. Suppose $\cb,\cb' \in \R^k$ with $0 \le c'_i \le c_j$ for all $i, j \in [k]$ with $i \leq j$.
Then for any $e \in R := R^{JZ}_{P\Qart}(\nu,\cb)$ and $e' \sub e$ we have $e' \in R' := R^{JZ}_{P\Qart}(\nu,\cb')$.
\end{lemma}

\proof
Since $e \in R$, $e$ is $\Qart_{R}$-partite, and so $e'$ is $\Qart_{R}$-partite, which is identical to being $\Qart_{R'}$-partite. Next, since $e \in R$, for any $e'' \sub e$ of size $j$ at most $\nu^{2^{-j}} n_1^{j}n^{k-j}$ edges of~$Z$ intersect each member of $\{V_i : i \in e''\}$. In particular, this property holds for any $e'' \sub e'$. Finally, since $e \in R$, we have $|J[\bigcup_{i \in e} V_i]| \geq c_{|e|}n_1^{|e|}$. Since $J$ is a $k$-complex and $c'_{|e'|} \le c_{|e|}$, it follows that $|J[\bigcup_{i \in e'} V_i]| \geq c'_{|e'|}n_1^{|e'|}$.
\endproof

The next lemma shows that when the $k$-graph~$Z$ is sparse, we do indeed have the property mentioned above, namely that the minimum degree sequence of the original $k$-system~$J$ is `inherited' by the reduced $k$-system. We work in the more general context of minimum $F$-degree sequences defined in Chapter~\ref{sec:mindeg}.

\begin{lemma} \label{reducedgraphminimumdegree}
Suppose that $1/n \ll 1/h, \nu \ll c_k \ll \dots \ll c_1 \ll \alpha, 1/k, 1/r$. Let~$X$ be a set of~$rn$ vertices, $\Part$ partition $X$ into $r$ parts $X_1,\dots,X_r$ of size $n$, and $\Qart$ refine $\Part$ into $h$ parts of equal size. Let $J$ be a $\Qart$-partite $(k-1)$-system on $X$, $Z$ a $\Qart$-partite $k$-graph on $X$ with $|Z| \leq \nu n^k$, and $P$ a vertex-equitable $\Qart$-partition $(k-1)$-complex on $X$ with clusters $V_1,\dots,V_{rm}$ of size $n_1$. Suppose also that $J$ and $Z$ are $\Part F$-partite for some allocation $F$, and that $\delta_0^F(J) = n$. Then, with respect to $\Part_R$, the reduced $k$-system $R := R^{JZ}_{P\Qart}(\nu,\cb)$ satisfies
$$\delta^F(R) \geq \left((1-k\nu^{1/2})m, \left(\frac{\delta^F_1(J)}{n}-\alpha\right)m, \dots, \left(\frac{\delta^F_{k-1}(J)}{n}-\alpha\right)m\right).$$
\end{lemma}

\proof
Let $U_1, \dots, U_r$ be the parts of $\Part_R$ corresponding to $X_1, \dots, X_r$ respectively. Fix any $f \in F$, $j \in [k-1]$ and $S = \{u_1, \dots, u_j\}\in R_j$ with $u_i \in U_{f(i)}$ for $i \in [j]$. Then $S$ is $\Qart_R$-partite and $\Part_R F$-partite, and $|J[V_S]| \geq c_jn_1^j$. For any edge $e \in J[V_S]$, we may write $e = \{v_1,\dots,v_j\}$ with $v_i \in X_{f(i)}$ for $i \in [j]$. There are at least $\delta^F_j(J)$ vertices $v_{j+1} \in X_{f(j+1)}$ such that $\{v_1,\dots,v_{j+1}\} \in J$, and of these at most $jrn/h$ belong to the same part of $\Qart$ as one of $v_1,\dots,v_j$. Thus we obtain at least $|J[V_S]|(\delta^F_j(J)-jrn/h)$ edges in sets $J[V_T]$, $T \in \mc{T}$,
where $\mc{T}$ denotes the collection of $\Qart_R$-partite sets $T = S \cup \{u\}$ for some $u \in U_{f(j+1)} \sm S$. 
At most $mc_{j+1}n_1^{j+1} = c_{j+1}n n_1^j$ of these edges belong to sets $J[V_T]$ of size less than $c_{j+1}n_1^{j+1}$, so there are at least  $\delta^F_j(J)m/n - \alpha m/2$ sets $T \in \mc{T}$ with $|J[V_T]| \geq c_{j+1}n_1^{j+1}$.

Such a set $T$ is an edge of $R$, unless for some $T' \sub T$ there are more than $\nu^{2^{-{j'}-1}} n_1^{j'+1}n^{k-j'-1}$ edges of $Z$ which intersect each part of $V_{T'}$, where $|T'|=j'+1$ for some $0 \le j' \le j$; in this case we say that $T$ is \emph{$T'$-bad}. Note that such a set $T'$ is not contained in $S$, otherwise $S$ would be $T'$-bad, contradicting $S \in R$. So we can write $T' = S' \cup \{u\}$, where $S' \sub S$ and $T = S \cup \{u\}$. For fixed $S' \sub S$, there can be at most $\nu^{2^{-{j'}-1}}m$ such vertices $u$, otherwise $S$ would be $S'$-bad, contradicting $S \in R$. Summing over all $S' \sub S$ we find that there are at most $2^j \nu^{2^{-k}} m$ vertices $u \in [U_{f(j+1)}]$ such that $T = S \cup \{u\}$ is $T'$-bad for some $T' \subseteq T$. It follows that $$\delta^F_j(R) \geq \delta^F_j(J)m/n - \alpha m/2 - 2^j \nu^{2^{-k}} m \ge \delta^F_j(J)m/n - \alpha m.$$

It remains to show that $\delta^F_0(R) \geq (1 - k\nu^{1/2})m$, i.e.\ for each $i \in [r]$ which lies in the image of some $f \in F$, the set $\{u\}$ is an edge of $R$ for all but at most $k\nu^{1/2}m$ vertices $u \in U_i$. For any $u \in U_i$, the set $\{u\}$ is $\Qart_R$-partite, and since $\delta_0^f(J) = n$ we have $|J[V_u]| = n_1$. So~$\{u\}$ is an edge of~$R$ unless more than $\nu^{1/2}n_1n^{k-1}$ edges of~$Z$ intersect~$V_u$; since $|Z| \leq \nu n^k$ this is true for at most $k\nu^{1/2}m$ vertices $u$, as required.
\endproof

The next proposition shows that the $k$-edges of the reduced $k$-system are useful, in that the corresponding clusters contain a sub-$k$-graph of $J_k$ to which the blow-up lemma can be applied; here we note that the complex $G'$ obtained in the proposition meets the conditions of Theorem~\ref{blowup} (with $\theta=\nu^{1/3}$ and $d$ replaced by $d/2$).

\begin{prop}\label{findregularcomplex}
Suppose that $1/n \ll \eps \ll d_a \ll 1/a \ll \nu \ll d \ll 1/k$. Let $X$ be a set of $n$ vertices, $\Qart$ partition $X$ into parts of equal size, $G$ and $Z$ be $\Qart$-partite $k$-graphs on $X$, and $P$ be a vertex-equitable $a$-bounded $\Qart$-partition $(k-1)$-complex on $X$, with clusters of size $n_1$. Suppose that $G$ is perfectly $\eps$-regular with respect to~$P$, and $U_1, \dots, U_k$ are clusters of $P$ such that $U := U_1 \cup \dots \cup U_k$ satisfies $|Z[U]| \leq \nu n_1^k$ and $|G[U]| \geq dn_1^k$. Let $\Part$ denote the partition of $U$ into parts $U_1, \dots, U_k$. Then there exists a $\Part$-partite $k$-complex $G'$ on $U$ such that $G'_k \sub G$, $G'$ is $\eps$-regular, $d_{[k]}(G') \geq d/2$, $d(G') \geq d_a$, and $Z' = Z \cap G'_k$ has $|Z'| \leq \nu^{1/3} |G'_k|$.
\end{prop}

\proof
To find $G'$ we select a suitable cell of $P^*$, which we recall is the $a^k$-bounded $\Qart$-partition $k$-complex formed from $P$ by weak equivalence. Note that one of the partitions forming $P^*$ is a partition of $K_{[k]}(U)$ into cells $C_1, \dots, C_s$, where $s \le a^k$. So at most $dn_1^k/3$ edges of $K_{[k]}(U)$ lie within cells $C_{i}$ such that $|C_{i}| \leq dn_1^k/(3a^k)$.  Also, since $|Z[U]| \leq \nu n_1^k$, at most $\nu^{1/2}n_1^k$ edges of $K_{[k]}(U)$ lie within cells $C_{i}$ such that $|Z \cap C_{i}| \geq \nu^{1/2}|C_{i}|$. Then, since $|G[U]| \geq dn_1^k$, at least $dn_1^k/2$ edges of $G$ must lie within cells $C_{i}$ with $|C_{i}| > dn_1^k/(3a^k)$ and $|Z \cap C_{i}| < \nu^{1/2}|C_{i}|$. By averaging, there must exist such a cell $C_{i}$ that also satisfies $|G \cap C_{i}| > d|C_{i}|/2$. Fix such a choice of $C_{i}$, which we denote by $C$.

We now define $G'$ to be the complex with top level $G'_k = G \cap C$ and lower levels $G'_{<k} = C_{<k} = \bigcup_{S' \subset S} C_{S'}$, where $S \in C$. Then $G'_k \sub G$, and $G'$ is $\eps$-regular, as $G$ is perfectly $\eps$-regular with respect to $P$. Furthermore, we have $$d_{[k]}(G') = \frac{|G'_{[k]}|}{|(G')^*_{[k]}|} = \frac{|G \cap C|}{|C|} > d/2$$ and $$d(G') = \frac{|G'_{[k]}|}{|K_{[k]}(U)|} = \frac{|G \cap C|}{|C|} \cdot \frac{|C|}{n_1^k} > \frac{d^2}{6a^k}>d_a.$$
Finally, since $|G'_k| = |G \cap C| > d|C|/2$ and $|Z'| \leq |Z \cap C| < \nu^{1/2}|C|$, we have $|Z'| < \frac{2\nu^{1/2}}{d} |G'_k| < \nu^{1/3}|G'_k|$.
\endproof

Finally, we need the following lemma, which states (informally) that index vectors where the original $k$-system $J$ is dense and the forbidden $k$-graph $Z$ is sparse are inherited as index vectors where the reduced $k$-system $R$ is dense. If $J$ is a $k$-system on $X$, and $\Part$ partitions $X$, then we denote by $J_\ib$ the set of edges in $J$ with index $\ib$ with respect to $\Part$.

\begin{lemma} \label{edgesofGtoedgesofR}
Suppose that $1/n \ll \nu \ll c_k \ll \dots \ll c_1 \ll \mu, 1/k$. Let $X$ be a set of $n$ vertices, $\Qart$ partition $X$ into parts of equal size, $J$ be a $\Qart$-partite $k$-system, $Z$ a $\Qart$-partite $k$-graph with $|Z| \leq \nu n^k$, and $P$ a vertex-equitable $\Qart$-partition $(k-1)$-complex on $X$, with clusters $V_1,\dots,V_m$ of size $n_1$. Let $R := R^{JZ}_{P\Qart}(\nu,\cb)$ be the reduced $k$-system. Then for any partition $\Part$ of $V(J)$ which is refined by the partition of $V(J)$ into clusters of $P$ we have the following property:
\begin{description}
\item[(Inheritance of index vectors)]
If $|J_\ib| \ge \mu n^p$, where $p=|\ib|$, then $|R_\ib| \ge \mu m^p/2$, where index vectors are taken with respect to $\Part$ (for $J$) and $\Part_R$ (for $R$).
\end{description}
\end{lemma}

\proof
Let $\mc{B}$ be the set of $\Qart_R$-partite sets $S$ of size $p$ with $S \notin R_p$. We estimate the number of edges contained in all $p$-graphs $J_p[V_S]$ with $S \in \mc{B}$. There are two reasons for which we may have $S \notin R_p$. The first is that $|J_p[V_S]| < c_pn_1^p$; this gives at most $c_p n^p$ edges in total. The second is that $S$ is $S'$-bad for some $S' \sub S$, i.e.\ more than $\nu^{2^{-j'}} n_1^{j'}n^{k-j'}$ edges of~$Z$ intersect each member of $\{V_i : i \in S'\}$, where $j'=|S'|$. For any given $j' \le p$, there can be at most $\binom{p}{j'}\nu^{1/2} m^p$ such sets $S$, otherwise we would have at least $\binom{p}{j'}\nu^{1/2} m^{j'}$ sets $S'$ of size $j'$ for which more than $\nu^{1/2} n_1^{j'}n^{k-j'}$ edges of~$Z$ intersect each member of $\{V_i : i \in S'\}$, contradicting $|Z| \leq \nu n^k$. Summing over $j' \le p$, this gives at most $2^p \nu^{1/2} m^pn_1^p \leq c_pn^p$ edges in total. Thus there are at most $2c_pn^p$ edges contained in all $p$-graphs $J_p[V_S]$ with $S \in \mc{B}$. Now, if at least $\mu n^p$ edges $e \in J_p$ have $\ib_\Part(e) = \ib$, then at least $\mu n^p/2$ of these edges lie in some $J_p[V_S]$ with $S \in R_p$. For each such $S$ we have $\ib_{\Part_R}(S) = \ib$ and $|J_p[V_S]| \le n_1^p$. Thus there are at least $\mu m^p/2$ edges $S \in R_p$ with $\ib_{\Part_R}(S) = \ib$.
\endproof

\section{Proof of Lemma~\ref{RANDOMEDGESELECTION}.}

In this section it is more convenient to use the (much simpler) Weak Regularity Lemma, in the context of a simultaneous regularity partition for several hypergraphs. Suppose that $\Part$ partitions a set $V$ into $r$ parts $V_1, \dots, V_r$ and $G$ is a $\Part$-partite $k$-graph on $V$. For $\eps>0$ and $A \in \binom{[r]}{k}$, we say that the $k$-partite sub-$k$-graph $G_A$ is \emph{$(\eps,d)$-vertex-regular} if for any sets $V'_i \sub V_i$ with $|V'_i| \ge \eps |V_i|$ for $i \in A$, writing $V' = \bigcup_{i \in A} V'_i$, we have $d(G_A[V']) = d \pm \eps$. We say that $G_A$ is \emph{$\eps$-vertex-regular} if it is $(\eps,d)$-vertex-regular for some $d$. The following lemma has essentially the same proof as that of the Szemer\'edi Regularity Lemma \cite{Sz}, namely iteratively refining a partition until an `energy function' does not increase by much; in this case the energy function is the sum of the mean square densities of the $k$-graphs with respect to the partition.

\begin{theo} (Weak Regularity Lemma) \label{weakrl}
Suppose that $1/n \ll 1/m \ll \eps \ll 1/t, 1/r \le 1/k$. Suppose that $G^1,\dots,G^t$ are $k$-graphs on a set $V$ of $n$ vertices. Then there is a partition $\mc{P}$ of $V$ into $m' \le m$ parts, such that there is some $n_0$ so that each part of $\mc{P}$ has size $n_0$ or $n_0+1$, and for each $i \in [t]$, all but at most $\eps n^k$ edges of $G^i$ belong to $\eps$-vertex-regular $k$-partite sub-$k$-graphs.
\end{theo}

We also require the fact that regularity properties are inherited by random subsets with high probability. We use the formulation given by Czygrinow and Nagle \cite[Theorem 1.2]{CN2} (a similar non-partite statement was proved earlier by Mubayi and R\"odl \cite{MR}).

\begin{theo} \label{inheritreg}
Suppose that $1/n \ll 1/s, \eps, c \ll \eps', 1/k$, that $\Part$ partitions a set $V$ into $k$ parts $V_1, \dots, V_k$ of size at least $n$, and $G$ is an $(\eps,d)$-vertex-regular $k$-partite $k$-graph on $V$. Suppose that $s_i \ge s$ for $i \in [k]$, and sets $S_i \sub V_i$ of size $s_i$ are independently chosen uniformly at random for $i \in [k]$. Let $S = \bigcup_{i \in [k]} S_i$. Then $G[S]$ is $(\eps',d)$-vertex-regular with probability at least $1 - e^{-c\min_i s_i}$.
\end{theo}

Here and later in the paper we will need the following inequality known as the Chernoff bound, as applied to binomial and hypergeometric random variables. We briefly give the standard definitions. The binomial random variable with parameters $(n,p)$ is the sum of $n$ independent copies of the $\{0,1\}$-valued variable $A$ with $\mb{P}(A=1)=p$. The hypergeometric random variable $X$ with parameters $(N,m,n)$ is defined as $X = |T \cap S|$, where $S \sub [N]$ is a fixed set of size $m$, and $T \sub [N]$ is a uniformly random set of size $n$. If $m=pN$ then both variables have mean $pn$.

\begin{lemma} \cite[Corollary 2.3 and Theorem 2.10]{JLR} \label{chernoff}
Suppose $X$ has binomial or hypergeometric distribution and $0 < a < 3/2$. Then
$\mb{P}(|X - \mb{E}X| \ge a\mb{E}X) \le 2e^{-\frac{a^2}{3}\mb{E}X}$.
\end{lemma}

Loosely speaking, the proof of Lemma~\ref{RANDOMEDGESELECTION} proceeds as follows. We divide the edges of $M$ into `blocks', where each block contains one edge of each index vector in $I(F)$. Since $M$ $\alpha'$-represents $F$, we can do this so that only a few edges of $M$ are not included in some block (we discard these edges). We then choose $M'$ to consist of the edges in the blocks covering $u$ and $v$, as well as the edges in a further $n'/b - 2$ blocks chosen uniformly at random, where $b$ is the number of edges in each block. The set $V'$ of vertices covered by $M'$ then includes $n'$ vertices from each part of $\Part$, and induces a matched $k$-system $(J', M')$. Also, a fairly straightforward application of the Chernoff bound implies that with high probability the induced $k$-system $(J', M')$ inherits the minimum degree sequence of $J$ (with slightly greater error terms) giving condition (i).

The main difficulty is therefore to show that condition (ii) also holds with high probability, for which we use weak hypergraph regularity. This argument is similar in spirit to methods used for property-testing (see e.g.~\cite{RS2}), though it does not follow directly from them. To see the main ideas of this part of the proof, suppose that we apply Theorem~\ref{weakrl} to partition each $V_j$ into $m$ equally-sized clusters, and assume for the sake of simplicity that any $k$ clusters $U_1, \dots, U_k$ from different parts of $\Part$ induce a vertex-regular sub-$k$-graph of $J_k$ (in reality this will only be true of most such sets of $k$ clusters). With high probability our randomly chosen $V'$ will include approximately the same number of vertices from each cluster, and by Theorem~\ref{inheritreg} the induced subclusters $U_1', \dots, U'_k$ of any $k$ clusters $U_1, \dots, U_k$ from different parts of $\Part$ will also induce a vertex-regular subgraph of $G$ with similar density.
Now suppose some set $S' \subseteq V'$ contains $\lfloor pn'/k\rfloor$ vertices from each $V_j$. Then $S'$ must have not-too-small intersection with approximately $pm/k$ subclusters from each part of $\Part$. But if $S'$ has not-too-small intersection with each of $U_1', \dots, U_k'$, then the vertex-regularity of $U'_1, \dots, U'_k$ implies the density of the induced sub-$k$-graph $G[S' \cap (U_1' \cup \dots \cup U'_k)]$ is approximately equal to the density of $G[U_1' \cup \dots \cup U_k']$, which in turn is approximately equal to the density of $G[U_1 \cup \dots \cup U_k]$. It follows that the density of $G[S']$ is approximately equal to the density of $G[S]$, where $S$ is the union of the clusters whose subclusters have not-too-small intersection with $S'$. However, $S$ must contain close to $pn/k$ vertices in each part of $\Part$, and so condition (ii) on $J$ implies that $G[S']$ is not too small, from which we deduce that condition (ii) on $J'$ holds for (our arbitrary choice of)~$S'$.

In the above sketch we glossed over several significant difficulties. Perhaps foremost among these is the fact that Theorem~\ref{inheritreg} requires the choice of vertices to be uniformly random, which is certainly not the case here. To resolve this difficulty, rather than apply Theorem~\ref{weakrl} directly to $G$, we instead apply it to a family of auxiliary $k$-graphs $J^{TC}$, each of which has the blocks of edges of $M'$ as vertices. Ignoring the blocks which cover the vertices of $u$ and $v$, our random choice is then to take a subset of the vertex set of each $J^{TC}$ uniformly at random, so we may indeed apply Theorem~\ref{inheritreg}.
 
We now give the full details of the the proof of Lemma~\ref{RANDOMEDGESELECTION}. We introduce new constants $\eps, \eps', c, d, m$ such that $1/n' \ll 1/m \ll \eps,c \ll \eps' \ll d \ll \alpha'$. Let $e_u$ and $e_v$ be the edges of $M$ containing $u$ and $v$ respectively. Since $M$ $\alpha'$-represents $F$, we can partition $M$ as $M_0 \cup M_1$, where $M_1$ is $F$-balanced, $\{e_u,e_v\} \sub M_1$ and $|M_0| \le \alpha' |M|$. Letting $V_0$ be the set of vertices covered by $M_0$, we have $|V_0|=k|M_0| \le \alpha'rn$. Since $|M|=rn/k$ we can write $|M_1|=rn_1/k$ with $n_1 \ge (1-\alpha')n$. By uniformity of $F$, $M_1$ covers $n_1$ vertices in each part of $\Part$. Let $I := I(F)$ (recall that if $I(F)$ is a set then this is the set $\{\ib(f): f \in F\}$). So $|I| = |F|/k!$. We partition $M_1$ into blocks, where each block contains one edge $e'$ with $\ib(e') = \ib$ for each $\ib \in I$. So each block consists of $|I|$ edges, and uses $b := k|I|/r$ vertices from each part. We arbitrarily label the blocks as $B_i$, $i \in [n_B]$, where $n_B := |M_1|/|I| = n_1/b$, and for each $j \in [r]$ we let $B_{ij}$ be the vertices in $V_j$ used by edges of $B_i$, which we arbitrarily label by $[b]$. For $x \in [n_B]$, $t \in [r]$, $c \in [b]$ we let $v_{xtc}$ denote vertex $c$ in $B_{xt}$. We construct $M'$ as a union of blocks.  We start by taking two distinct blocks $B_{x_u}$ and $B_{x_v}$, where $e_u \in B_{x_u}$, and either $e_v \in B_{x_v}$ or $e_v \in B_{x_u}$ and $B_{x_v}$ is arbitrary. Then we let $X$ be a subset of $[n_B] \sm \{x_u,x_v\}$ of size $n'/b - 2$ chosen uniformly at random, and let $X' = X \cup \{x_u,x_v\}$, so $|X'|=n'/b$. Take 
$$M' = \bigcup_{x \in X'} B_x, \hspace{.5cm} V' = \bigcup_{e \in M'} e, \hspace{.5cm} J' = J[V'] \mbox{ and } V'_i = V_i \cap V' \mbox{ for each $i \in [r]$}.$$ Then we have $u,v \in V'_a$, $|V'_1| = \dots = |V'_r| = n'$, and $(J', M')$ is a matched $\Part F$-partite $k$-system on $V'$.

Fix $f \in F$ and consider the minimum $f$-degree sequence $\delta^f(J')$. We clearly have $\delta^f_0(J') = n'$, as for any $v \in V'$ we have $\{v\} \in J$, and so $\{v\} \in J'$. Now fix any $e = \{v_1,\dots,v_j\} \in J$ for some $j \in [k-1]$ such that $v_i \in V_{f(i)}$ for $i \in [j]$. For each $i \in [j]$ let $B_{d_i}$ be the block containing $v_i$, so $\mc{B}_e := B_{x_u} \cup B_{x_v} \cup \bigcup_{B_{d_i}}$ is the union of $\ell \leq j+2$ distinct blocks. Let $A$ be the set of vertices covered by edges in $\mc{B}_e$. Write $q = f(j+1)$, and for $c \in [b]$, let $N(e)_c$ be the set of vertices $v \in V_q \sm (V_0 \cup A)$ such that $e \cup \{v\} \in J$ and $v=v_{xqc}$ for some $x \in [n_B]$. Then 
$\sum_{c \in [b]} |N(e)_c| \ge \left( \frac{k-j}{k} - \alpha \right)n - \alpha'rn - \ell br$ by assumption on $\delta^f(J)$. Writing $|N(e)_c| = \theta_c n_1$, and $C = \{c \in [b]: \theta_c > 2\eps\}$ we have 
$$\sum_{c \in C} \theta_c \ge 1-j/k - \alpha - 2\alpha'r - 2b\eps.$$ Note that, conditioning on the event $e \in J'$, the random variable $|N(e)_c \cap V'_q|$ is hypergeometric with mean $\frac{n'/b - \ell}{n_1/b - \ell} |N(e)_c| > (\theta_c - \eps)n'$ for each $c \in C$, so the Chernoff bound (Lemma~\ref{chernoff}) gives $|N(e)_c \cap V'_q| > (\theta_c - 2\eps)n'$ with probability at least $1 - 2e^{-\eps^3 n'/3}$. On these events we have 
$$m^f(e) := \sum_{c \in [b]} |N(e)_c \cap V'_q| > \sum_{c \in C} (\theta_c - 2\eps)n' > (1-j/k - 2\alpha)n'.$$ So, letting $Z(e, f)$ denote the event that $m^f(e) \leq (1-j/k - 2\alpha)n'$, and no longer conditioning on $e \in J'$, we have
$$\Prob(e \in J' \textrm{ and } Z(e,f)) = \Prob(e \in J') \Prob(Z(e,f) \mid e \in J' ) 
\leq \left(\frac{n'/b - \ell}{n/b - \ell }\right)^{\ell - 2}\cdot 2e^{-\eps^3 n'/3}.$$
Taking a union bound over at most $|F|(k-1)(k+1)$ choices of $f \in F$, $j \in [k-1]$ and $\ell \in [j+1]$, and at most $(n_B)^{\ell-2}(b\ell)^j \leq n^{\ell-2}(2bk)^k$ edges $e \in J_j$ such that $\mc{B}_e$ is the union of $\ell$ distinct blocks, we see that the minimum $F$-degree property of $J'$ holds with high probability.

For property (ii) we use weak hypergraph regularity. Fix $p \in [k-1]$ and consider the following auxiliary $(p+1)$-graphs $J^{TC}$ on $[n_B]$. Given $1 \le x_1 < \dots < x_{p+1} \le n_B$, $T = (t_1,\dots,t_{p+1}) \in [r]^{p+1}$ such that there is some $f \in F$ for which $f(i)=t_i$ for $i \in [p+1]$, and $C = (c_1,\dots,c_{p+1}) \in [b]^{p+1}$, we say that $\{x_1, \dots, x_{p+1}\} \in J^{TC}$ if and only if $\{ v_{x_it_ic_i} : i \in [p+1] \} \in J_{p+1}$. Thus each edge of $J_{p+1}$ that has at most one vertex in any block corresponds to a unique edge in at most $(p+1)!$ of the $(p+1)$-graphs $J^{TC}$ (these are given by holding $C$ constant and permuting $T$). Observe that this accounts for most edges of $J_{p+1}$, as at most $br^{p+1}n^p$ edges have more than one vertex in some block.  By Theorem~\ref{weakrl}, there is a partition $\mc{P}'$ of $[n_B]$ into $m' \le m$ parts, such that there is some $n_0$ so that each part of $\mc{P}'$ has size $n_0$ or $n_0+1$, and for each $T$ and $C$ as above, all but at most $\eps n^{p+1}$ edges of $J^{TC}$ belong to $\eps$-vertex-regular $\Part'$-partite sub-$(p+1)$-graphs. Note that $m'n_0 \le n_B \le m'(n_0+1)$. Similarly to the previous argument, we will see that with high probability $X'$ represents all parts of $\mc{P}'$ approximately equally. Indeed, fix $U \in \mc{P'}$, and note that, since $n_B = n_1/b$, the random variable $|X' \cap U \sm \{x_u,x_v\}|$ is hypergeometric with mean $$\frac{n'/b - 2}{n_1/b - 2} (n_0 \pm 2) = \frac{n'}{n_1}(n_B/m' \pm 4) = \frac{n'}{bm'} \pm 1.$$ By the Chernoff bound, we have $|X' \cap U| = (1 \pm \eps) n'/bm'$ with probability at least $1 - 2e^{-\eps^2 n'/6bm'}$. We restrict attention to choices of $X'$ such that this estimate holds for all $U \in \mc{P'}$. Note that conditional on any specified values of $|X' \cap U|$ for $U \in \mc{P'}$, the choices of $X' \cap U \sm \{x_u,x_v\}$, $U \in \mc{P'}$ are independent uniformly random $|X' \cap U \sm \{x_u,x_v\}|$-sets in $U \sm \{x_u,x_v\}$ for $U \in \mc{P'}$. Now consider the `reduced' $(p+1)$-graphs $R^{TC}$, where $V(R^{TC}) = \mc{P}'$, and $E(R^{TC})$ consists of all $(p+1)$-tuples $(U_1,\dots,U_{p+1})$ of parts of $\mc{P}'$, such that writing $U = \bigcup_{i \in [p+1]} U_i$, $J^{TC}[U]$ is $(\eps,d')$-vertex-regular for some $d' \ge d$. Observe that at most $(d+\eps)n^{p+1}$ edges of $J^{TC}$ do not belong to $J^{TC}[U]$ for some such edge. Given an edge $(U_1,\dots,U_{p+1}) \in E(R^{TC})$, and a random choice of $X'$ such that $|X' \cap U| = (1 \pm \eps) n'/bm'$ for all $U \in \mc{P'}$, Theorem~\ref{inheritreg} implies that $J^{TC}[U \cap X']$ is $(\eps',d')$-vertex-regular with probability at least $1 - e^{-cn'/2bm'}$. We can assume that this holds for all $T$, $C$ and $(U_1,\dots,U_{p+1})$, as there are at most $(rbm)^{p+1}$ choices, so we can take a union bound. 

Now consider any sets $S'_t \sub V'_t$ such that $|S'_t| = \lfloor pn'/k \rfloor$ for each $t \in [r]$. We need to show that there are at least $\beta' (n')^{p+1}$ edges in $J_{p+1}[S']$, where $S' := \bigcup_{t \in [r]} S'_t$. Let $S'_{tc} = \{x \in [n_B]: v_{xtc} \in S'\}$ for $t \in [r]$ and $c \in [b]$. Let $Y_{tc} = \{U \in \mc{P'}: |S'_{tc} \cap U| \ge \beta^2|X' \cap U|\}$. For each $t \in [r]$ we have
\[pn'/k = |S'_t| = \sum_{c \in [b]} \sum_{U \in \mc{P'}} |S'_{tc} \cap U| \le \sum_{c \in [b]} |Y_{tc}|(1+\eps) n'/bm' + b\beta^2|X'|,\]
so 
$$\sum_{c \in [b]} |Y_{tc}| > \frac{pn'/k - b\beta^2 n'}{(1+\eps) n'/bm'} > (p/k - 2b\beta^2)bm'.$$
Write $S^0_t = \bigcup_{c \in [b]} \bigcup_{U \in Y_{tc}} U$. Then $|S^0_t| > (p/k - 2b\beta^2)bm'n_0 > (p/k - 3b\beta^2)n$, since $m'(n_0+1) \ge n_B = n/b$. Let $S_t \sub V_t$ be an arbitrary set of size $\lfloor pn/k \rfloor$ such that $S^0_t \sub S_t$ if $|S^0_t| \le pn/k$ or $S_t \sub S^0_t$ if $|S^0_t| \ge pn/k$.  Let $S = S_1 \cup \dots \cup S_r$ and $S_{tc} = \{x \in [n_B]: v_{xtc} \in S\}$ for $t \in [r]$ and $c \in [b]$. By assumption (ii), there are at least $\beta n^{p+1}$ edges in $J_{p+1}[S]$. Of these edges, we discard a small number of edges that are `bad' for one of the following reasons: at most $3rk\beta^2 n^{p+1}$ edges which are incident with $\bigcup_{t \in [r]} (S_t \sm S^0_t)$, at most $br^{p+1}n^p$ edges which have more than one vertex in some block, and at most $(br)^{p+1} (d+\eps) n^{p+1}$ edges which correspond to an edge in some $J^{TC}$ belonging to a $(p+1)$-partite sub-$(p+1)$-graph that is not $(\eps,d')$-vertex-regular with $d' \ge d$. This still leaves at least $(\beta-\beta^{3/2}) n^{p+1}$ edges, each of which corresponds to an edge $\{x_1, \dots, x_{p+1}\}$ in some $J^{TC}$ such that $x_i \in S_{t_ic_i} \cap U_i$ for some $U_i \in Y_{t_ic_i}$ for $i \in [p+1]$ and $(U_i: i \in [p+1])$ is an edge of $R^{TC}$. Let $Q^{TC}$ be the set of edges $(U_i: i \in [p+1])$ in $R^{TC}$ such that $U_i \in Y_{t_ic_i}$ for $i \in [p+1]$. Then we have
\[\sum_{T,C} |Q^{TC}| \ge (\beta-\beta^{3/2}) n^{p+1}/(n_0+1)^{p+1} \ge \frac{1}{2}\beta(bm')^{p+1}.\]
Recall that each edge of $J$ with at most one vertex in any block corresponds to at most $(p+1)!$ edges in the $(p+1)$-graphs $J^{TC}$. So $(p+1)!|J_{p+1}[S']|$ is at least the number of edges $\{x_1, \dots, x_{p+1}\}$ in some $J^{TC}$ such that $x_i \in S'_{t_ic_i} \cap U_i$ for $i \in [p+1]$ where $(U_i: i \in [p+1]) \in Q^{TC}$. For each $(U_i: i \in [p+1]) \in Q^{TC}$ we have at least 
$$(d-\eps')\prod_{i \in [p+1]} |S'_{t_ic_i} \cap U_i| > \frac{1}{2}d (\beta^2 n'/bm')^{p+1}$$ such edges of $J^{TC}$, using the definition of $Y_{t_ic_i}$ and $\eps'$-vertex-regularity. Since $\sum_{T,C} |Q^{TC}| \ge \frac{1}{2}\beta(bm')^{p+1}$ we obtain at least $\beta' (n')^{p+1}$ edges in $J_{p+1}[S'] = J'_{p+1}[S']$. \qed

\chapter{Matchings in $k$-systems} \label{sec:proofs}

In this chapter we prove our theorems on matchings in $k$-systems with a minimum degree sequence condition. In fact, we prove theorems in the setting of minimum $F$-degree sequences, which simultaneously generalise both our non-partite and partite theorems. In the first section we prove the general form of the fractional perfect matching result. We combine this with hypergraph regularity in the second section to prove a common generalisation of Theorems~\ref{almostpacking} and~\ref{almostpackingpartite} on almost perfect matchings. In the third section we apply transferrals to prove a common generalisation of Theorems~\ref{newmain} and~\ref{newmainpartite} on perfect matchings; we will also see that essentially the same proof gives a common generalisation of Theorems~\ref{prestability} and~\ref{prestabilitypartite}. 
 
\section{Fractional perfect matchings.} 

In this section we prove a lemma that generalises Lemma~\ref{fractionalmatching} to the minimum $F$-degree setting. Let $J$ be a $\Part F$-partite $k$-system on $V$, where $\Part$ is a balanced partition of $V$ into $r$ parts of size $n$, and $F$ is a $(k,r)$-uniform allocation. Recall that a fractional perfect matching in $J_k$ is an assignment of a weight $w_e \geq 0$ to each edge $e \in J_k$ such that for any $v \in V(J_k)$ we have $\sum_{e \ni v} w_e = 1$. The lemma will show that if $J$ satisfies our minimum $F$-degree condition and has no space barrier then $J_k$ admits a fractional perfect matching (in fact, with a slightly stronger minimum degree the latter condition is not required). We actually prove something stronger, namely that $J_k$ contains a fractional perfect matching which is \emph{$F$-balanced}, in that $\sum_{e \in J_k :~ \ib(e) = \ib} w_e$ is constant over all $\ib \in I(F)$; this can be seen as a fractional equivalent of an $F$-balanced matching as previously defined. We also say that a multiset $E$ in $J_k$ is \emph{$F$-balanced} if the number of edges in $E$ of index $\ib$ (counted with multiplicity) is the same for any $\ib \in I(F)$. First we need the following proposition. Let $\T$ be the collection of all sets $T \subseteq J_k$ which contain one edge of index $\ib$ for each $\ib \in I(F)$.

\begin{prop} \label{fracmatchequiv} The following statements are equivalent.
\begin{enumerate}[(i)]
\item $\1 \in PC(\{\chi(T) : T \in \T\})$.
\item $k|I(F)|\1/rn \in CH(\{\chi(T) : T \in \T\})$.  
\item $J_k$ admits an $F$-balanced fractional perfect matching in which at most $|I(F)|(rn+1)$ edges have non-zero weight.
\item $J_k$ admits an $F$-balanced fractional perfect matching.
\end{enumerate}
\end{prop}

\proof
Let $X = \{\chi(T) : T \in \T\}$. Suppose first that (i) holds, so that $\1 = \sum_{\xb \in X} c_\xb \xb$ with $c_\xb \ge 0$ for $\xb \in X$. Then we have $$k|I(F)| \sum_{\xb \in X} c_\xb = \sum_{\xb \in X} c_\xb \xb \cdot \1 = \1 \cdot \1 = rn.$$  Since any $\xb \in X$ has non-negative integer coordinates we also have $c_\xb \leq 1$ for each $\xb$.  Then by multiplying each $c_\xb$ by $k|I(F)|/rn$ we obtain (ii). Now suppose that (ii) holds. Then by Theorem~\ref{caratheodory} we may write $k|I(F)|\1/rn = \sum_{\xb \in X} c_\xb \xb$ with $c_\xb \ge 0$ for $\xb \in X$ so that at most $rn+1$ of the $c_\xb$'s are non-zero. For each $\xb \in X$ assign weight $w_T = rnc_\xb/k|I(F)|$ to some $T \in \T$ with $\chi(T)=\xb$. Then assigning to each edge $e \in J_k$ the weight $w_e := \sum_{T \in \T :~e \in T} w_T$ gives an $F$-balanced fractional perfect matching in which at most $|I(F)|(rn+1)$ edges have non-zero weight, so we have (iii). Trivially (iii) implies (iv), so it remains to show that (iv) implies (i). 

Consider an $F$-balanced fractional perfect matching in $J_k$, where $w_e$ denotes the weight of an edge $e$, so $\sum_{e \in J_k :~\ib(e) = \ib} w_e$ is constant over $\ib \in I(F)$. We assign weights to sets $T \in \T$ and modify the weights of edges $e \in J_k$ according to the following algorithm. Suppose at some step we have weights $w'_e$ for $e \in J_k$, such that $\sum_{e \in J_k :~\ib(e) = \ib} w'_e$ is constant over $\ib \in I(F)$. Suppose that not all weights are zero, and choose $e_0 \in J_k$ with the smallest non-zero weight. For every $\ib \in I(F)$ with $\ib \neq \ib(e_0)$, since 
$$\sum_{e \in J_k :~\ib(e) = \ib} w'_e = \sum_{e \in J_k :~\ib(e) = \ib(e_0)} w'_e \ne 0,$$
we may choose $e_\ib$ in $J_k$ with $\ib(e_\ib) = \ib$ and $w'_{e_\ib} \ne 0$. Note that $w'_{e_\ib} \geq w'_{e_0}$, by minimality of $w'_{e_0}$. Let $T \in \T$ consist of $e_0$ and the edges $e_\ib$ for $\ib \in I(F)$ with $\ib \neq \ib(e_0)$. We assign weight $w'_{e_0}$ to $T$, and define new weights by $w''_{e_0}=0$, $w''_{e_\ib} = w'_{e_\ib} - w'_{e_0}$, and $w''_e=w'_e$ otherwise. Then, with these new weights, $\sum_{e \in J_k :~\ib(e) = \ib} w''_{e}$ remains constant over $\ib \in I(F)$, and $w''_e \geq 0$ for every $e \in J_k$.  Since the number of edges of zero weight has increased by at least one, after at most $|J_k|$ iterations every edge will have zero weight, at which point we stop. Then for any $\xb \in X$, we let $c_\xb$ be the sum of the weights assigned to any $T$ with $\chi(T) = \xb$. By construction we have $\1 = \sum_{\xb \in X} c_\xb \xb$, as required.
\endproof

The following lemma generalises Lemma~\ref{fractionalmatching} to the minimum $F$-degree setting. Indeed,  Lemma~\ref{fractionalmatching} is the $\alpha=0$ `furthermore' statement of Lemma~\ref{fractionalmatchingFpartite}, applied with $r=1$ and $F$ generated by the unique function $f:[k]\to [1]$. The main statement in Lemma~\ref{fractionalmatchingFpartite} shows that the same conclusion holds under a slightly weaker $F$-degree sequence if there is no space barrier.

\begin{lemma} \label{fractionalmatchingFpartite}
Suppose that $1/n \ll \alpha' \ll \beta, 1/D_F, 1/k, 1/r$ and that $\alpha \leq \alpha'$. Let $V$ be a set partitioned into parts $V_1, \dots, V_r$ each of size $n$, and $F$ be a $(k,r)$-uniform connected allocation with $|F| \leq D_F$. Also let $J$ be a $\Part F$-partite $k$-system on $V$ such that 
\begin{enumerate}[(i)]
\item $\delta^F(J) \geq \left(n, \frac{(k-1)n}{k} - \alpha n, \frac{(k-2)n}{k} - \alpha n, \dots, \frac{n}{k} - \alpha n \right)$, and
\item for any $p \in [k-1]$ and sets $S_i \subseteq V_i$ with $|S_i| = \lfloor pn/k \rfloor$ for $i \in [r]$ we have $|J_{p+1}[S]| \geq \beta n^{p+1}$, where $S:= \bigcup_{i \in [r]} S_i$.
\end{enumerate}  
Then $J_k$ admits an $F$-balanced fractional perfect matching. Furthermore, if $\alpha = 0$ then this conclusion holds even without assuming condition (ii), any lower bound on $n$, or that $F$ is connected.
\end{lemma}

\proof
First we give a construction that reduces to the case when $k \mid n$. For $i \in [r]$, let $V'_i$ be a set of size $kn$ consisting of copies $v_i(j)$, $j \in [k]$ of each $v_i \in V_i$. Let $J'$ be the $k$-system whose edges are all possible copies of the edges of $J$. Then $\delta^F(J') = k\delta^F(J)$, so~(i) holds for $J'$. To see that~(ii) holds also, fix $p \in [k-1]$ and sets $S'_i \subseteq V'_i$ with $|S'_i| = pn$ for $i \in [r]$. Let $S_i \subseteq V_i$ consist of all vertices $v_i \in V_i$ with a copy $v_i(j)$ in $S'_i$. Then $|S_i| \geq |S'_i|/k \geq \lfloor pn/k \rfloor$. So $|J_{p+1}[S]| \geq \beta n^{p+1}$ by (ii), where $S:= \bigcup_{i \in [r]} S_i$. Each edge of $J_{p+1}[S]$ has at least one copy in $J'_{p+1}[S']$, where $S' := \bigcup_{i \in [r]} S'_i$. We deduce that $|J'_{p+1}[S']| \geq \beta n^{p+1} \geq (\beta/k^{p+1})(kn)^{p+1}$, so~(ii) holds for $J'$ with $\beta/k^{p+1}$ in place of $\beta$. Assuming the result when $k \mid n$, we find that $J'_k$ admits an $F$-balanced fractional perfect matching. From this we obtain an $F$-balanced fractional perfect matching in $J_k$, where the weight of an edge is obtained by combining the weights of its copies in $J'_k$ and dividing by $k$. Thus we can assume $k \mid n$.

Now suppose for a contradiction that $J_k$ has no $F$-balanced fractional perfect matching. Let $\T$ be as above; then Proposition~\ref{fracmatchequiv} implies that $\1 \notin PC(\{\chi(T) : T \in \T\})$. So by Farkas' Lemma (Lemma~\ref{farkas}), there is some $\ab \in \R^{rn}$ such that $\ab \cdot \1 < 0$ and $\ab \cdot \chi(T) \geq 0$ for every $T \in \T$. Note that any $F$-balanced multiset $E$ in $J_k$ can be expressed as $\sum_{i=1}^z T_i$ for some $T_1,\dots,T_z$ in $\T$, and so satisfies $\ab \cdot \chi(E) \ge 0$. For $i \in [r]$ let $v_{i,1}, \dots, v_{i,n}$ be the vertices of~$V_i$, and let $a_{i,1}, \dots, a_{i,n}$ be the corresponding coordinates of~$\ab$, where the vertex labels are chosen so that $a_{i,1} \leq a_{i,2} \leq \dots \leq a_{i,n}$ for each $i \in [r]$.

For any multisets $S$ and $S'$ in $V$ of equal size $m$, we say that~$S$ \emph{dominates}~$S'$, and write $S' \le S$, if we may write $S = \{v_{i_1, j_1}, \dots, v_{i_{m}, j_{m}}\}$ and $S' = \{v_{i'_1, j'_1}, \dots, v_{i'_{m}, j'_{m}}\}$ so that for each $\ell \in [m]$ we have $i_\ell = i'_\ell$ and $j'_\ell \leq j_\ell$. Note that $\le$ is a transitive binary relation. We also observe that if~$S' \le S$ then $\ab \cdot \chi(S') \leq \ab \cdot \chi(S)$. As usual, for a multiset $E$ in $J_k$ we write $\chi(E) = \sum_{e \in E} \chi(e)$, thus identifying $E$ with the multiset in $V$ in which the multiplicity of $v \in V$ is the number of edges in $E$ containing it, counting with repetition. We recall that $S+T$ denotes the multiset union of two multisets $S$ and $T$, and for $i \in \N$, $i S$ denotes the multiset union of $i$ copies of $S$. We extend our `arithmetic' of multisets to include subtraction, writing $S-T$ for the multiset $A$ such that $A+T=S$, if it exists. It will also be convenient to manipulate formal expressions $S-T$ that do not correspond to multisets, via the rule $(S-T) + (S'-T') = (S+S') - (T+T')$, which can be understood as a shorthand for $\chi(S)-\chi(T) + \chi(S')-\chi(T') = \chi(S)+\chi(S') - \chi(T)+\chi(T')$.

We start by proving the $\alpha=0$ statement. For a mental picture, it is helpful to think of the vertices arranged in a grid, with $r$ columns corresponding to the parts $V_i$, $i \in [r]$, and $k$ rows, where the subsquare in column $i$ and row $j$ contains the vertices $v_{i,(j-1)n/k+s}$, $s \in [n/k]$. We partition $V$ into sets $(W_s: s \in [n/k])$, where 
$$W_s = \{v_{i,(j-1)n/k+s}: i \in [r], j \in [k]\}$$ 
consists of all vertices at \emph{height} $s$, for each column $i$ and row $j$. Note that for each $s \in [n/k]$ we have $W_s \ge W_1$. Thus we have 
$$0 > \ab \cdot \1 = \sum_{s \in [n/k]} \ab \cdot \chi(W_s) \ge (n/k) \ab \cdot \chi(W_1).$$ To obtain the required contradiction, we will show that a constant multiple of $W_1$ dominates an $F$-balanced multiset of edges of $J_k$. To see that this suffices, define 
$$X^f := \{v_{f(j),(j-1)n/k+1}: j \in [k]\}$$ 
for $f \in F$. Since $F$ is $(k,r)$-uniform we have $\sum_{f \in F} X^f = |F|r^{-1}W_1$. Now, by the minimum $F$-degree of $J$, we may greedily form an edge $e^f = \{v_{f(1),d_1}, \dots, v_{f(k), d_k}\} \in J$ with $d_1 = 1$ and $d_j \leq (j-1)n/k + 1$ for each $2 \leq j \leq k$. Then $X^f$ dominates $e^f$, so $|F|r^{-1}W_1 = \sum_{f \in F} X^f$ dominates the $F$-balanced multiset $\{e^f : f \in F\}$. It follows that $\ab \cdot \chi(W_1) \ge 0$, so we have the required contradiction to the assumption that $J_k$ has no $F$-balanced fractional perfect matching.

Now consider the case $0 < \alpha < \alpha'$. We will obtain a contradiction by a similar strategy to that used when $\alpha=0$, namely partitioning $V$ into `dominating' sets, where we say a multiset $S$ is \emph{dominating} if some constant multiple of $S$ dominates an $F$-balanced multiset of edges in $J_k$. Note that if $S$ is dominating we have $\ab \cdot \chi(S) \ge 0$, so this will give the contradiction $\ab \cdot \1 \ge 0$. We may assume that $C := \alpha n$ is an integer. We also let $N$ be an integer with $\alpha' \ll 1/N \ll 1/D_F, 1/k, 1/r$ and $|F| \mid N$. Now we define sets $W_s$, $s \in [n/k-CN]$ of size $rk$ by
\[ W_s := \{v_{i, s} : i \in [r]\} \cup \{v_{i, (j-1)n/k + C + s} : i \in [r], 2 \leq j \leq k\}.\]
Note that this agrees with our previous definition in the case when $C=\alpha=0$; now we have increased by $C$ the height of the vertices in rows $2$ to $k$. Again we have $W_s \ge W_1$, and we will show that a constant multiple of $W_1$ dominates an $F$-balanced multiset of edges. To see this, define 
$$X^f := \{v_{f(1), 1}\} \cup \{v_{f(j),(j-1)n/k+C+1}: 2 \leq j \leq k\}$$ for $f \in F$ 
(again, this agrees with our previous notation when  $C=\alpha=0$). Since $F$ is $(k,r)$-uniform we have $\sum_{f \in F} X^f = |F|r^{-1}W_1$. By the minimum $F$-degree of $J$, we may greedily form an edge $e^f = \{v_{f(1),d_1}, \dots, v_{f(k), d_k}\} \in J$ with $d_1 = 1$ and $d_j \leq (j-1)n/k + C + 1$ for each $2 \leq j \leq k$. Then $X^f$ dominates $e^f$, so $\sum_{f \in F} X^f = |F|r^{-1}W_1$ dominates $\{e^f: f \in F\}$.

Thus we have arranged that most vertices in $V$ belong to dominating sets $W_s$, but in each column  we still need to deal with the $CN$ highest vertices in row $1$, and the $C$ lowest and $C(N-1)$ highest vertices in rows $2$ to $k$. We will partition these into sets $Z_s$, $s \in [C]$ of size $rkN$, so that each $Z_s$ contains $N$ vertices in each of the $rk$ subsquares of the grid of vertices, and in rows $2$ to $k$ these comprise $1$ `low' vertex and $N-1$ `high' vertices. The formal definition is as follows. Define $P := [k] \times \{0, \dots, N-1\} \cup \{(1, N)\} \sm \{(k, 0)\}$. Then for each $s \in [C]$, let 
$$ Z_s := \left\{v_{i, jn/k - tC + s} : (j, t) \in P, i \in [r]\right\}.$$
Note that the sets $W_s$, $s \in [n/k-CN]$ and $Z_s$, $s \in [C]$ partition $V$. Next we define a multiset that is dominated by each $Z_s$. Write $y_{i,j} = v_{i, jn/k - CN}$ for $i \in [r]$, $j \in [k]$,
\[ Y = \{ y_{i, j} : i \in [r], j \in [k]\}, \text{ and } D = NY + \{ y_{i, 1} : i \in [r]\} - \{ y_{i, k} : i \in [r]\}.\]
Thus the $CN$ highest vertices in the subsquare of column $i$ and row $j$ are each above the corresponding vertex $y_{i,j}$ of $Y$. We claim that each $Z_s$ dominates $D$. To see this, note that for each $i \in [r]$ and $2 \le j \le k-1$, the $N$ copies of $y_{i,j}$ in $D$ are dominated by the $N-1$ high vertices $v_{i, jn/k - tC + s}$, $t \in [N-1]$ in row $j$ and the low vertex $v_{i, jn/k + s}$ in row $j+1$. In row $k$ we have removed one copy of each $y_{i,k}$ from $D$, so the remaining $N-1$ copies are dominated by the high vertices in row $k$. In row $1$ we have $N+1$ copies of $y_{i,j}$ in $D$, which are dominated by the $N$ vertices $v_{i, n/k - tC + s}$, $t \in [N]$ in row $1$ and the low vertex $v_{i, n/k + s}$ in row $2$. Thus $Z_s \ge D$. 

The remainder of the proof is showing that $D$ is a dominating multiset; we divide this into 3 claims. The first claim exploits the absence of a space barrier to find edges that are lower than those guaranteed by the minimum degree condition. Whereas the minimum degree condition gives edges in which the $i$th vertex is (close to being) in row $i$ or below for $i \in [k]$, for each $p \in [k-1]$ we can find an edge where the first $p+1$ vertices are in row $p$ or below, and the $i$th vertex is (close to being) in row $i$ or below for $p+2 \le i \le k$. Intuitively, we can think of such an edge as having a `$p$-demoted' vertex, in that the $(p+1)$st vertex is lower than guaranteed by the minimum degree condition (although the first $p$ vertices may be higher). We also require that a demoted vertex is not too near the top of the row it has demoted too, in that it is below the corresponding vertex $y_{i,j}$. In the first claim we have no control over which of the parts $V_i$ contains a demoted vertex, so in the second claim we exploit the connectivity of $F$ to construct a multiset with a demoted vertex in any desired part. Finally, in the third claim we take an appropriate non-negative linear combination of multisets with $p$-demoted vertices for all $p$ to obtain a multiset dominated by $D$.

For the first claim we need edges dominated by one of the sets
\[B_p^f := \{y_{f(j), p} : j \in [p+1]\} \cup \{y_{f(j), j} : p+2 \leq j \leq k\}.\]

\begin{claim} \label{domination1}
For any $p \in [k-1]$ there is some $f_p \in F$ and $e_p \in J_k$ with $e_p \le B_p^{f_p}$.
\end{claim} 

To prove the claim, for each $i \in [r]$ let $S_i = \{v_{i, d} : d \leq pn/k - CN\}$, and arbitrarily choose $S_i'$ of size $pn/k$ containing $S_i$. Let $S = \bigcup_{i \in [r]} S_i$ and $S' = \bigcup_{i \in [r]} S'$. Then by condition (ii) we have $|J_{p+1}[S']| \geq \beta n^{p+1}$. At most $rCN (rn)^p < \beta n^{p+1}$ edges of $J_{p+1}$ intersect $S' \sm S$, so $|J_{p+1}[S]| > 0$. We may therefore choose an edge $e = \{v_{f(1), d_1}, \dots, v_{f(p+1), d_{p+1}}\} \in J_{p+1}$ where $d_j \leq pn/k - CN$ for $j \in [p+1]$, for some function $f : [p+1] \to [r]$. Since $J$ is $\Part F$-partite, $f$ must be the restriction of some $f_p \in F$. Then by the minimum $f_p$-degree assumption on $J$ we can greedily extend $e$ to an edge $e_p = \{v_{f(1), d_1}, \dots, v_{f(k), d_k}\}$ with $d_j \leq (j-1)n/k + C + 1$ for $p+2 \leq j \leq k$. Thus $B_p^{f_p}$ dominates $e_p$, which proves Claim \ref{domination1}.

\medskip

For the second claim we will find $F$-balanced multisets of edges dominated by the multisets 
\[D_p^\ell = 2k|F|r^{-1} Y + (p+1) \{y_{\ell,p}\} - \{y_{\ell,j} : j \in [p+1]\}.\]

\begin{claim} \label{domination2}
For any $p \in [k-1]$ and $\ell \in [r]$ there is an $F$-balanced multiset $E^\ell_p$ in $J_k$ with $E^\ell_p \le D^\ell_p$.
\end{claim}

To prove the claim, we start by applying Claim \ref{domination1}, obtaining $e_p \in J_k$ and $f_p \in F$ with $e_p \le B_p^{f_p}$. Since $F$ is connected, there is a connected graph $G_F$ on $[r]$ such that for every $ii' \in E(G_F)$ and $j, j' \in [k]$ with $j \ne j'$ there is $f\in F$ with $f(j)=i$ and $f(j')=i'$. Choose for each $j \in [p+1]$ a path $f_p(j) = i_1^j, \dots, i_{s_j+1}^j = \ell$ in $G_F$ from $f_p(j)$ to $\ell$. For each $z \in [s_j]$, let $f^j_z \in F$ be such that $f^j_z(j) = i_z^j$ and $f^j_z(p) = i_{z+1}^j$, and let $\hat{f}^j_z$ be obtained from $f^j_z$ by swapping the values of $f^j_z(j)$ and $f^j_z(p)$. Since $F$ is invariant under permutation we have $\hat{f}^j_z \in F$. Now recall that for each $f \in F$ we have an edge $e^f \le X^f$, where $X^f = \{v_{f(j),(j-1)n/k+C+1}: j \in [k]\}$. We define $Y^f := \{y_{f(j),j}: j \in [k]\}$, and note that $Y^f \ge X^f \ge e^f$ and $\sum_{f \in F} Y^f = |F|r^{-1}Y$. Next we show that we can write $D_p^\ell$ as
\begin{equation} \label{eq:dpl}
D_p^\ell = 2k\sum_{f \in F} Y^f + (B_p^{f_p} - Y^{f_p}) + \sum_{j \in [p+1]} \sum_{z \in [s_j]} (Y^{{f}^j_z} - Y^{\hat{f}^j_z}).
\end{equation}
To see this, note that $Y^{{f}^j_z} - Y^{\hat{f}^j_z} = \{y_{f^j_z(j), j}, y_{f^j_z(p), p}\} - \{y_{f^j_z(j), p}, y_{f^j_z(p), j}\} =  (\{y_{i^j_z, j} \} - \{y_{i^j_z, p} \}) - (\{ y_{i^j_{z+1}, j}\} - \{ y_{i^j_{z+1}, p} \})$ for each $j \in [p+1], z \in [s_j]$, so 
\[\sum_{z \in [s_j]} (Y^{{f}^j_z} - Y^{\hat{f}^j_z}) = (\{y_{f_p(j), j} \} - \{y_{f_p(j), p} \}) - (\{ y_{\ell, j}\} - \{ y_{\ell, p} \})\] 
for each $j \in [p+1]$. Then
\begin{align*}
\sum_{j \in [p+1]} \sum_{z \in [s_j]} (Y^{{f}^j_z} - Y^{\hat{f}^j_z}) 
& = \sum_{j \in [p+1]} (\{y_{f_p(j), j} \} - \{y_{f_p(j), p} \}) - \sum_{j \in [p+1]} (\{ y_{\ell, j}\} - \{ y_{\ell, p} \}) \\
& = Y^{f_p} - B_p^{f_p} + (p+1) \{y_{\ell,p}\} - \{y_{\ell,j} : j \in [p+1]\},
\end{align*}
which proves (\ref{eq:dpl}). To define $E^\ell_p$, we take $2k$ copies of $e^f$ for each $f \in F$, replace one copy of $e^{f_p}$ by a copy of $e_p$, and replace one copy of $e^{\hat{f}^j_z}$ by one copy of $e^{{f}^j_z}$ for each $j \in [p+1], z \in [s_j]$; note that there are enough copies for these replacements, as $f^j_z$, $z \in [s_j]$ are distinct. Thus $E^\ell_p$ contains $2k|F|$ edges and is $F$-balanced, since each edge was replaced by another of the same index. To see that $E^\ell_p \le D^\ell_p$, we assign edges in $E^\ell_p$ to terms on the right hand side of (\ref{eq:dpl}), counted with multiplicity according to their coefficient, such that each edge of $E^\ell_p$ is dominated by its assigned term. Before the replacements, we have $2k$ copies of $e^f$ for each $f \in F$, which we assign to the $2k$ copies of $Y^f$. To account for the replacement of a copy of $e^{f_p}$ by $e_p$, we remove one of the assignments to $Y_{f_p}$ and assign the replaced copy of $e_p$ to $B_p^{f_p}$. To replace a copy of $e^{\hat{f}^j_z}$ by $e^{{f}^j_z}$, we remove one of the assignments to $Y^{\hat{f}^j_z}$ and assign the replaced copy of $e^{{f}^j_z}$ to $Y^{{f}^j_z}$. Thus we have an assignment showing that $E^\ell_p \le D^\ell_p$, which proves Claim \ref{domination2}.

\medskip

Finally, we will establish the required property of $D$ by taking an appropriate non-negative linear combination of the multisets $E^\ell_p$ and $D^\ell_p$. This is given by the following claim.

\begin{claim} \label{domination3}
There are non-negative integer coefficients $m_j$, $j \in [k]$ such that  $$M := \sum_{j=1}^{k} m_j \leq k^{k} \mbox{ and } D = (N-2k|F|M)Y +  \sum_{\ell \in [r]} \sum_{p \in [k-1]} m_p D^\ell_p.$$ 
\end{claim}

To prove the claim, we define the coefficients recursively by $m_k = 0$, $m_{k-1} = 1$, $m_{k-2} = k-1$ and 
$$m_j = (j+2)m_{j+1} - (j+3)m_{j+2} \textrm{ for } 1 \leq j \leq k-3.$$
Since $m_j \leq k m_{j+1}$ for any $j \in [k-1]$ we have $M \le k^k$. Next we show that 
\begin{equation} \label{eq:M}
(s+1)m_s - \sum_{j = s-1}^{k} m_j =
\begin{cases}
-1 & s = k, \\
0 & k-1 \ge s \ge 2, \\
1 & s = 1.
\end{cases}
\end{equation}
For $s=k$ we have $(k+1)m_k - m_{k-1} - m_k = -1$, and for $s=k-1$ we have $km_{k-1} - \sum_{j = k-2}^{k} m_j = k - (k-1) - 1 = 0$. 
Also, for $k-2 \ge s \ge 2$ we have
\[(s+1)m_s - \sum_{j = s-1}^{k} m_j = m_{s-1} + (s+2)m_{s+1} - \sum_{j = s-1}^{k}m_j = (s+2)m_{s+1} - \sum_{j = s}^{k} m_j, \]
so these cases follow from the case $s=k-1$. Finally, writing $m_0 = 0$, we have 
\[\sum_{s=1}^k \left((s+1)m_s - \sum_{j=s-1}^k m_j\right) = m_k - m_0 = 0;\] 
this implies the case $s=1$, and so finishes the proof of (\ref{eq:M}). Now we show by induction that $m_j \geq (j+2)m_{j+1}/2$ for $j=k-2,k-3,\dots,2$. For the base case $j=k-2$ we have $m_{k-2} = k \ge k/2 = km_{k-1}/2$. Now suppose $2 \le j \le k-3$ and $m_{j+1} \geq (j+3)m_{j+2}/2$. Then $m_j - (j+2)m_{j+1}/2 \geq (j+2)m_{j+1}/2 - (j+3)m_{j+2} \geq (j-2)m_{j+1}/2 \geq 0$. Therefore $m_j \ge 0$ for $2 \le j \le k$. Then by (\ref{eq:M}), we also have $m_1 = (1+\sum_{j = 0}^k m_j) \ge 0$. Now by definition of $D_p^\ell$, for any $\ell \in [r]$ we have
\begin{align*} 
\sum_{p \in [k-1]} m_p D^\ell_p & = \sum_{p \in [k-1]} m_p(2k|F|r^{-1} Y + (p+1) \{y_{\ell,p}\} - \{y_{\ell,j} : j \in [p+1]\}) \\
& = 2Mk|F|r^{-1}Y + \sum_{s \in [k]} \brac{(s+1)m_s - \sum_{j = s-1}^{k} m_j} \{y_{\ell,s}\}.
\end{align*}
Summing over $\ell$ and applying (\ref{eq:M}), we obtain $\sum_{\ell \in [r]} \sum_{p \in [k-1]} m_p D^\ell_p = 2Mk|F|Y + \{y_{\ell,1}: \ell \in [r]\} - \{y_{\ell,k}: \ell \in [r]\}$. By definition of $D$, this proves Claim \ref{domination3}.

\medskip

To finish the proof of the lemma, we show that $D$ is dominating, by defining an $F$-balanced multiset $E$ in $J_k$ with $E \le D$. We let $E$ be the combination of $Nr/|F| - 2Mrk$ copies of $\{e^f : f \in F\}$ with $m_p$ copies of $E^\ell_p$ for each $\ell \in [r]$ and $p \in [k]$ (note that $Nr/|F| - 2Mrk$ is positive since $N \gg D_F, k, r$). Then $E$ is $F$-balanced, since each $E^\ell_p$ is $F$-balanced. Now recall from the proof of Claim \ref{domination2} that we can write $\sum_{f \in F} Y^f = |F|r^{-1}Y$ with $Y^f \ge e^f$ for $f \in F$. Substituting this in the expression for $D$ in Claim \ref{domination3} we can see that $D$ dominates $E$ termwise: each $Y^f$ in $D$ dominates an $e^f$ in $E$, and each $D^\ell_p$ in $D$ dominates an $E^\ell_p$ in $E$. Now recall that $V$ is partitioned into sets $W_s$, $s \in [n/k-CN]$ and $Z_s$, $s \in [C]$, where each $W_s$ dominates the dominating set $W_1$, and each $Z_s$ dominates $D$, which we have now shown is dominating. We deduce that $\ab \cdot \1 \ge 0$, which gives the required contradiction to the assumption that $J_k$ has no $F$-balanced fractional perfect matching.
\endproof

\section{Almost perfect matchings.}

In this section we prove a lemma that will be used in the following section to prove a common generalisation of Theorems~\ref{prestability} and~\ref{prestabilitypartite}, in which we find a matching covering all but a constant number of vertices, and a common generalisation of Theorems~\ref{newmain} and~\ref{newmainpartite}, where we find a perfect matching under the additional assumption that there is no divisibility barrier. The following lemma is a weaker version of the common generalisation of Theorems~\ref{prestability} and~\ref{prestabilitypartite}. It states that if a $k$-system $J$ satisfies our minimum $F$-degree condition and does not have a space barrier, then $J_k$ contains a matching which covers all but a small proportion of $V(J)$. One should note that the proportion $\psi$ of uncovered vertices can be made much smaller than the deficiency $\alpha$ in the $F$-degree sequence. The proof will also show that a slightly stronger degree condition yields the same conclusion even in the presence of a space barrier, from which we shall deduce a common generalisation of Theorems~\ref{almostpacking} and~\ref{almostpackingpartite} on almost perfect matchings.

\begin{lemma} \label{approxalmostpackingFpartite}
Suppose that $1/n \ll \psi \ll \alpha \ll \beta, 1/D_F, 1/r, 1/k$. Let $\Part$ be a partition of a set $V$ into parts $V_1, \dots, V_r$ of size $n$ and $F$ be a connected $(k,r)$-uniform allocation with $|F| \leq D_F$. Suppose that $J$ is a $\Part F$-partite $k$-system on $V$ such that 
\begin{enumerate}[(i)]
\item $\delta^F(J) \geq \left(n, \left(\frac{k-1}{k} - \alpha\right)n, \left(\frac{k-2}{k} - \alpha\right) n, \dots, \left(\frac{1}{k} - \alpha\right) n\right)$, and
\item for any $p \in [k-1]$ and sets $S_i \subseteq V_i$ with $|S_i| = \lfloor pn/k \rfloor$ for $i \in [r]$ we have $|J_{p+1}[S]| \geq \beta n^{p+1}$, where $S:= \bigcup_{i \in [r]} S_i$.
\end{enumerate} 
Then $J_k$ contains an $F$-balanced matching $M$ which covers all but at most $\psi n$ vertices of $J$.
\end{lemma}

\proof
Introduce new constants with
\begin{align*}
1/n \ll \eps \ll d_a \ll 1/a \ll \nu, 1/h \ll c_k \ll \dots \ll c_1 \ll \gamma \ll \psi \ll \alpha \ll \beta, 1/D_F, 1/r, 1/k.
\end{align*}
We may additionally assume that $r \mid h$. Since $r \mid |V(J)|$ and $r \mid a!h$, we may delete up to $a!h$ vertices of $J$ so that equally many vertices are deleted from each vertex class and the number of vertices remaining in each part is divisible by $a!h$. By adjusting the constants, we can assume for simplicity that in fact $a!h$ divides $n$, so no vertices were deleted. Fix any partition $\Qart$ of $V$ into $h$ parts of equal size which refines $\Part$, and let $J'$ be the $k$-system on $V$ formed by all edges of $J$ which are $\Qart$-partite. Then $J'_k$ is a $\Qart$-partite $k$-graph on $V$, so by Theorem~\ref{eq-partition} there is an $a$-bounded vertex-equitable $\Qart$-partition $(k-1)$-complex $P$ on~$V$ and a $\Qart$-partite $k$-graph $G$ on $V$ that is $\nu/r^k$-close to $J'_k$ and perfectly $\eps$-regular with respect to $P$. Let $Z := G \triangle J'_k$. Then since $G$ is $\nu/r^k$-close to $J'_k$ we have $|Z| \leq \nu n^k$. Also let $W_1, \dots, W_{rm_1}$ be the clusters of $P$ (note that $rm_1 \leq ah$), and let $n_1$ be their common size, so $n_1m_1 = n$.

Consider the reduced $k$-system $R := R^{J'Z}_{P\Qart}(\nu,\cb)$. Recall that $R$ has vertex set $[rm_1]$, where vertex $i$ corresponds to cluster $W_i$ of $P$, and that $\Part_R$ and $\Qart_R$ are the partitions of $[rm_1]$ corresponding to $\Part$ and $\Qart$ respectively. So $R$ is $\Part_R F$-partite and $\Qart$-partite. For each $i \in [r]$ let $U_i$ be the part of $\Part_R$ corresponding to part $V_i$ of $\Part$; since $P$ was vertex-equitable each part $U_i$ has size $m_1$. We now show that since $J$ had no space barrier, $R$ also does not have a space barrier, even after the deletion of a small number of vertices from each $U_i$. 

\begin{claim} \label{matchingclaim}
Suppose that sets $U'_i \subseteq U_i$ for $i \in [r]$ satisfy $|U'_1| = \dots = |U'_r| = m' \geq (1-\alpha) m_1$. Let $U' := \bigcup_{i \in [r]} U'_i$, $R' := R[U']$ and $\Part_{R'}$ be the restriction of $\Part_R$ to $U'$. Then for any $p \in [k-1]$ and sets $S'_i \subseteq U'_i$ with $|S'_i| = \lfloor pm'/k \rfloor$ for $i \in [r]$ we have $|R'_{p+1}[S']| \geq \beta (m')^{p+1}/10$, where $S' := \bigcup_{i \in [r]} S'_i$.
\end{claim}

To prove the claim, let $S = \bigcup_{i \in [r]} S_i$, where $S_i = \bigcup_{j \in S'_i} W_j \sub V_i \sub V$ for $i \in [r]$. Then $|S'_i| \ge n_1\lfloor pm'/k\rfloor \geq (1-2\alpha)pn/k$ for $i \in [r]$. Let $S'' = \bigcup_{i \in [r]} S''_i$, where for $i \in [r]$ we take $S''_i$ to be any set of size $\lfloor pn/k \rfloor$ such that $S_i \sub S''_i$ if $|S_i| \le \lfloor pn/k \rfloor$ and $S''_i \sub S_i$ if $|S_i| \ge \lfloor pn/k \rfloor$. By assumption (ii) on $J$ we have $|J_{p+1}[S'']| \geq \beta n^{p+1}$, and so $|J'_{p+1}[S'']| \geq \beta n^{p+1}/2$, since at most $n^{p+1}/h$ edges of $J_{p+1}$ are not edges of $J'_{p+1}$. Since at most $2\alpha rn$ vertices of $S''$ are not vertices of $S$, it follows that $|J'_{p+1}[S]| \geq \beta n^{p+1}/4$. Letting $\Sart$ denote the partition of $V(J')$ into parts $S$ and $V(J') \sm S$, we can rephrase this as at least $\beta n^{p+1}/4$ edges $e \in J'_{p+1}$ have $\ib_{\Sart}(e) = (p+1, 0)$. By Lemma~\ref{edgesofGtoedgesofR} it follows that at least $\beta m_1^{p+1}/8$ edges $e \in R_{p+1}$ have $\ib_{\Sart_{R}}(e) = (p+1, 0)$. Since $\Sart_{R}$ is the partition of $[m_1]$ into $S'$ and $[m_1] \sm S'$, we conclude that $|R_{p+1}[S']| \geq \beta m_1^{p+1}/8$, and therefore $|R'_{p+1}[S']| \geq \beta (m')^{p+1}/10$, proving the claim.\medskip

Now, since any edge $e \in J'$ has $d^F_{J'}(e) \ge d^F_J(e) - \alpha n/2$, by Lemma~\ref{reducedgraphminimumdegree} we have
\begin{equation*}
\delta^F(R) \geq ((1-k\nu^{1/2})m_1, (\delta_1(J)/n- \alpha/2)m_1, \dots, (\delta_{k-1}(J)/n - \alpha/2)m_1)
\end{equation*}
with respect to $\Part_R$. So there are at most~$k\nu^{1/2} m_1$ vertices $i$ in each part of $\Part_R$ for which~$\{i\}$ is not an edge of~$R$. Thus we can delete $k\nu^{1/2} m_1$ vertices from each part to obtain $R'$ with $m' = (1-k\nu^{1/2})m_1 \geq (1-\alpha)m_1$ vertices in each part and
\begin{equation} \label{eq:Fdeg}
\delta^F(R') \ge (m', (\delta_1(J)/n - \alpha)m', \dots, (\delta_{k-1}(J)/n - \alpha)m').
\end{equation}
Claim~\ref{matchingclaim} and (\ref{eq:Fdeg}) together show that $R'_k$ satisfies the conditions of Lemma~\ref{fractionalmatchingFpartite}, and so admits an $F$-balanced fractional perfect matching.
By Proposition~\ref{fracmatchequiv} it follows that $R'_k$ admits a $F$-balanced fractional perfect matching in which there are at most $|I(F)|(rm'+1)$ edges $e \in R'_k$ of non-zero weight. For $e \in R'_k$ let $w_e$ be the weight of $e$ in such a fractional matching.
So $\sum_{e \ni v} w_e = 1$ for any $v \in V(R')$ and $\sum_{e \in R'_k : ~\ib(e) = \ib} w_e$ has common value $m'r/k|I(F)|$ for every $\ib \in I(F)$. Next, partition each cluster~$W_i$ into parts $\{W_i^e : e \in R'_k\}$ such that
$$|W_i^e| = \begin{cases} w_e n_1 & \text{if $e$ is incident to vertex $i$ of $R'$,}
\\ 0& \text{otherwise.} \end{cases}$$
The lemma will now follow easily from the following claim.

\begin{claim}\label{claim_me}
For any $e \in R'_k$ there exists a matching $M_e$ in $J_k[\bigcup_{i \in e} W_i^e]$ of size at least $(w_e-\gamma)n_1$.
\end{claim}

To prove the claim, first note that if $w_e \leq \gamma$ then there is nothing to prove, so we may assume that $w_e > \gamma$. Let $M$ be a maximal matching in $J_k[\bigcup_{i \in e} W_i^e]$, and suppose for a contradiction that $|M| < (w_e - \gamma)n_1$. For each $i \in e$ let $W'_i$ consist of the vertices in $W^e_i$ not covered by $M$; then $|W'_i| \geq \gamma n_1$ for each $i \in e$. Now observe that since $e \in R'_k$, and therefore $e \in R_k$, we know that $e$ is $\Qart_R$-partite, that $|Z[\bigcup_{i \in e} W_i^e]| \leq \nu^{2^{-k}} n_1^k$, and that $|J'_k[\bigcup_{i \in e} W_i^e]| \geq c_kn_1^k$. Since $J'_k \sm Z \sub G_k$ we therefore have $|G[\bigcup_{i \in e} W_i^e]| \geq c_kn^1_k/2$. Since $G$ is perfectly $\eps$-regular with respect to $P$, by Proposition~\ref{findregularcomplex} there is a $k$-partite $k$-complex $G'$ whose vertex classes are~$W_i^e$ for $i \in e$ such that $G'$ is $\eps$-regular, $d_{[k]}(G') \geq c_k/4$, $d(G') \geq d_a$, $G'_k \sub G$, and $|Z \cap G'_k| \leq \nu^{2^{-k}/3} |G'_k|$.
Writing $W' := \bigcup_{i \in e} W'_i$, it follows by Lemma~\ref{regularrestriction} that $G'[W']$ has $d(G'[W']) \geq d(G')/2$. So $G_k'[W']$ contains at least $\gamma^k |G'_k|/2$ edges. Since $\nu \ll \gamma$, some edge of $G'_k[W']$ is not an edge of $Z$, and is therefore an edge of $J_k[W']$, contradicting the maximality of $M$. This proves Claim~\ref{claim_me}.
\medskip

To finish the proof of Lemma~\ref{approxalmostpackingFpartite}, we apply Claim~\ref{claim_me} to find matchings $M_e$ for
each edge $e \in R'_k$ with $w_e > 0$. Then the union $M := \bigcup_{e \in R'_k :~w_e > 0} M^e$ of all these matchings is a matching in $J_k$. Furthermore, by choice of the weights $w_e$, for any $\ib \in I(F)$ the number of edges $e' \in M$ with $\ib(e') = \ib$ is at least
\begin{align*}
& \sum_{e \in R'_k : w_e > 0,~\ib(e) = \ib } (w_e - \gamma)n_1 \geq \frac{rm'n_1}{k|I(F)|} -  \gamma |I(F)|(rm'+1)n_1 \\ 
& \geq \frac{(1-k\nu^{1/2})rm_1n_1 - 2\gamma kD_F^2 rm_1 n_1}{k|I(F)|}\geq \frac{|V(J)|- \psi n}{k|I(F)|}.
\end{align*}
So for each $\ib \in I(F)$ we may choose $(|V(J)|- \psi n)/k|I(F)|$ edges $e' \in M$ with $\ib(e') = \ib$; these edges together form an $F$-balanced matching in $J_k$ which covers all but at most $\psi n$ vertices of~$J$.
\endproof

Examining the above proof, we note that condition (ii) and the connectedness of $F$ were only used in the proof of Claim~\ref{matchingclaim}, which in turn was only used to show that $R'_k$ admits an $F$-balanced fractional perfect matching. If we instead assume the stronger $F$-degree condition $\delta^F(J) \geq \left(n, \left(\frac{k-1}{k} + \alpha\right)n, \left(\frac{k-2}{k} + \alpha\right) n, \dots, \left(\frac{1}{k} + \alpha\right) n\right)$ then we have $\delta^F(R') \geq \left(m', \frac{(k-1)m'}{k}, \dots, \frac{m'}{k}\right)$ by~(\ref{eq:Fdeg}). Then the `furthermore' statement of Lemma~\ref {fractionalmatchingFpartite} implies that $R'_k$ admits an $F$-balanced fractional perfect matching. Thus we have also proved the following lemma.

\begin{lemma} \label{approxalmostpackingFpartite2}
Suppose that $1/n \ll \psi \ll \alpha \ll 1/D_F, 1/r, 1/k$. Let $\Part$ be a partition of a set $V$ into parts $V_1, \dots, V_r$ of size $n$ and $F$ be a $(k,r)$-uniform allocation with $|F| \leq D_F$. Suppose that $J$ is a $\Part F$-partite $k$-system on $V$ such that 
\[\delta^F(J) \geq \left(n, \left(\frac{k-1}{k} + \alpha\right)n, \left(\frac{k-2}{k} + \alpha\right) n, \dots, \left(\frac{1}{k} + \alpha\right) n\right).\] 
Then $J_k$ contains an $F$-balanced matching $M$ which covers all but at most $\psi n$ vertices of $J$.
\end{lemma}

Now we deduce the common generalisation of Theorems~\ref{almostpacking} and~\ref{almostpackingpartite}. To do this, we add some `fake' edges to increase the minimum degree  of $J$, so that we may use Lemma~\ref{approxalmostpackingFpartite2}. This gives an almost perfect matching in the new system. The fake edges are chosen so that only a small number of them can appear in any matching, so we then remove them to give an almost perfect matching in the original system. 

\begin{theo} \label{almostpackingF}
Suppose that $1/n \ll \alpha \ll 1/D_F, 1/r, 1/k$. Let $\Part$ be a partition of a set $V$ into sets $V_1, \dots, V_r$ each of size $n$ and $F$ be a $(k,r)$-uniform allocation with $|F| \leq D_F$. Suppose that $J$ is a $\Part F$-partite $k$-system on $V$ with
$$\delta^F(J) \geq \left(n, \left(\frac{k-1}{k} - \alpha\right)n, \left(\frac{k-2}{k} - \alpha\right) n, \dots, \left(\frac{1}{k} - \alpha\right) n\right).$$
Then $J$ contains a matching covers all but at most $9k^2r\alpha n$ vertices of $J$.
\end{theo}

\proof
Choose a set $X$ consisting of $8k\alpha n$ vertices in each part $V_i$ for $i \in [r]$ uniformly at random. We form a $k$-system $J'$ on $V(J)$ whose edge set consists of every edge of $J$, and {\em fake edges}, which are every $S \in \binom{V(J)}{\leq k}$ which is $\Part F$-partite and intersects $X$. We claim that with high probability we have
\begin{equation} \label{randomchoiceinalmostpackingF}
\delta^F(J') \geq \left(n, \left(\frac{k-1}{k} + \alpha\right)n, \left(\frac{k-2}{k} + \alpha\right) n, \dots, \left(\frac{1}{k} + \alpha\right) n\right).
\end{equation}
To see this, first observe that $\delta^F_0(J') = \delta^F_0(J) = n$. Also, if $e$ is a fake edge, then $d^F_{J'}(e) \geq n - k$. Now consider any $f \in F$, $j \in [k-1]$ and $e \in J_j$ such that we may write $e=\{v_1,\dots,v_j\}$ with $v_i \in V_{f(i)}$ for $i \in [j]$. Then there are at least $\left(\frac{k-j}{k} - \alpha\right)n$ vertices $v_{j+1} \in V_{f(j+1)}$ such that
$\{v_1,\dots,v_{j+1}\} \in J$. We can assume that there is a set $S$ of $n/2k$ vertices $v_{j+1} \in V_{f(j+1)}$ such that
$\{v_1,\dots,v_{j+1}\} \notin J$ (otherwise we are done). Then $Y = |S \cap X|$ is hypergeometric with mean $4\alpha n$, so by Lemma~\ref{chernoff} we have $Y \ge 2\alpha n$ with probability at least $1-2e^{-\alpha n/3}$. On this event, there are at least $(k-j)n/k + \alpha n$ vertices $v_{j+1} \in V_{f(j+1)}$ such that
$\{v_1,\dots,v_{j+1}\} \in J'$.
Taking a union bound over at most $|F|kn^{k-1}$ choices of $f, j$ and $e$ we see that (\ref{randomchoiceinalmostpackingF}) holds with high probability. Fix a choice of $X$ such that~(\ref{randomchoiceinalmostpackingF}) holds. Then by Lemma~\ref{approxalmostpackingFpartite2} $J'$ contains a matching $M'$ which covers all but at most $\alpha n$ vertices of $J'$. Since every fake edge intersects $X$, there can be at most $|X| = 8kr\alpha n$ fake edges in $M$; deleting these edges we obtain a matching $M$ which covers all but at most $\alpha n + 8k^2r\alpha n \leq 9k^2r\alpha n$ vertices of $J$.
\endproof

To deduce Theorem~\ref{almostpacking} we apply Theorem~\ref{almostpackingF} with $r=1$ and $F$ consisting of ($k!$ copies of) the unique function $f:[k] \to [1]$, so $\delta^F(J)=\delta^f(J)=\delta(J)$. Similarly, to deduce Theorem~\ref{almostpackingpartite} we apply Theorem~\ref{almostpackingF} with $F$ consisting of all injections $f:[k] \to [r]$, so $\delta^F(J)=\delta^*(J)$ is the partite minimum degree sequence. Note that in the latter case an $F$-balanced matching is exactly our notion of a balanced matching from earlier.

\section{Perfect matchings.}

In this section we prove a common generalisation of Theorems~\ref{newmain} and~\ref{newmainpartite} on perfect matchings. Essentially the same proof will also give a common generalisation of Theorems~\ref{prestability} and~\ref{prestabilitypartite}.  Before giving the proof, for the purpose of exposition we sketch an alternative argument for the case of graphs $(k=2)$. For simplicity we just consider the cases $r=2$ and $r=1$. Suppose first that $r=2$, i.e.\ we have a bipartite graph $G$ with parts $V_1$ and $V_2$ of size $n$, with $\delta(G) \ge (1/2-\alpha)n$. Suppose that $G$ does not have a perfect matching. Then by Hall's theorem, there is $S'_1 \sub V_1$ with $|N(S'_1)|<|S'_1|$. By the minimum degree condition, each of $S'_1$ and $N(S'_1)$ must have size $(1/2 \pm \alpha)n$. Let $S_1$ be a set of size $\bfl{n/2}$ that either contains or is contained in $S'_1$ and let $S_2$ be a set of size $\bfl{n/2}$ that either contains or is contained in $V_2 \sm N(S'_1)$. Then $S_1 \cup S_2$ contains at most $2\alpha n^2$ edges, so we have a space barrier. 

Now suppose that $r=1$, i.e.\ we have a graph $G$ on $n$ vertices with $\delta(G) \ge (1/2-\alpha)n$, where $n$ is even.  Suppose that $G$ does not have a perfect matching. Then by Tutte's Theorem, there is a set $U \sub V(G)$ so that $G \sm U$ has more than $|U|$ odd components (i.e.\ connected components with an odd number of vertices). Suppose first that $|U| < (1/2-2\alpha)n$. Then $\delta(G \sm U) \ge \alpha n$, so $G \sm U$ has at most $\alpha^{-1}$ components. It follows that $|U| < \alpha^{-1}$. Then $\delta(G \sm U) \ge (1/2-2\alpha)n$, so $G \sm U$ has at most $2$ components. It follows that $|U|$ is $0$ or $1$. Since $n$ is even, $|U|=0$, and $G$ has two odd components, i.e.\ a divisibility barrier. Now suppose that $|U| > (1/2-2\alpha)n$. Then all but at most $2\alpha n$ of the odd components of $G \sm U$ are isolated vertices, i.e.\ $G$ contains an independent set $I$ of size $(1/2-4\alpha)n$. Let $S$ be a set of size $n/2$ containing $I$.  Then $S$ contains at most $2\alpha n^2$ edges, so we have a space barrier.

\begin{theo} \label{perfectmatchingF}
Let $1/n \ll \gamma \ll \alpha \ll \beta, \mu \ll 1/D_F, 1/r, 1/k$. Suppose $F$ is a $(k,r)$-uniform connected allocation with $|F| \leq D_F$, and that $b = k|I(F)|/r$ divides $n$. Let $\Part$ be a partition of a set $V$ into parts $V_1, \dots, V_r$ of size $n$, and $J$ be a $\Part F$-partite $k$-complex on $V$ such that
\begin{enumerate}[(i)]
\item $\delta^F(J) \geq \left(n, \left(\frac{k-1}{k} - \alpha \right)n, \left(\frac{k-2}{k} - \alpha \right) n, \dots, \left(\frac{1}{k} - \alpha \right) n\right)$,
\item for any $p \in [k-1]$ and sets $S_i \sub V_i$ such that $|S_i| = \bfl{pn/k}$ for each $i \in [r]$ there are at least $\beta n^k$ edges of $J_k$ with more than $p$ vertices in $S :=  \bigcup_{i \in [r]} S_i$, and
\item $L^\mu_{\Part'}(J_k)$ is complete with respect to $\Part$ for any partition~$\Part'$ of~$V(J)$ which refines~$\Part$ and whose parts each have size at least $n/k - \mu n$.
\end{enumerate}
Then $J_k$ contains a perfect matching which $\gamma$-represents $F$.
\end{theo}

\proof
We follow the strategy outlined in Section~\ref{sec:outline}. The first step is to use hypergraph regularity to decompose $J$ into an exceptional set and some clusters with a matched reduced $k$-system. Introduce new constants with
\begin{align*}
 1/n &\ll \eps \ll d^* \ll d_a \ll 1/a \ll \nu, 1/h \ll \theta \ll d, c \ll c'_k \ll \dots \ll c'_1  \\ & \ll c_k \ll \dots \ll c_1 \ll \psi \ll \gamma \ll 1/C \ll \alpha \ll \alpha' \ll \mu, \beta \ll 1/D_F,1/r,1/k.
\end{align*}
We also assume that $r \mid h$. By uniformity of $F$, any matching in $J_k$ which contains one edge of index $\ib$ for each $\ib \in I(F)$ covers $b$ vertices in each part of $\Part'$. Since $b \mid a!h$ and $b \mid n$, we can arbitrarily delete at most $a!h/b$ such matchings to make the number of vertices remaining in each part divisible by $a!h$. By adjusting the constants, we can assume for simplicity that in fact $a!h$ divides $n$, so no vertices were deleted. Fix any balanced partition~$\Qart$ of~$V$ into~$h$ parts which refines $\Part$, and let~$J'$ be the subcomplex obtained by deleting from $J$ all those edges which are not $\Qart$-partite. Then~$J'_k$ is a $\Qart$-partite $k$-graph with order divisible by~$a!h$, so by Theorem~\ref{eq-partition} there exists an $a$-bounded $\eps$-regular vertex-equitable $\Qart$-partition $(k-1)$-complex~$P$ on~$V(J)$ and a $\Qart$-partite $k$-graph~$G$ on~$V$ that is $\nu/{r^k}$-close to~$J'_k$ and perfectly $\eps$-regular with respect to~$P$. Let $Z = G \bigtriangleup J'_k$, so $|Z| \leq \nu n^k$ and any edge of $G \sm Z$ is also an edge of~$J_k$. Note that since $r \mid h$ the number of clusters of $P$ is divisible by $r$; let $W_1, \dots, W_{rm_1}$ be the clusters of $P$. Since $P$ is $a$-bounded we must have $rm_1 \leq ah$. In addition, since $P$ is vertex-equitable, each cluster $W_i$ has the same size; let $n_1 := |W_1| = \dots = |W_{rm_1}| = n/m_1$ be this common size.

Let $R^1 := R^{J'Z}_{P\Qart}(\nu,\cb)$ and $R^2 := R^{J'Z}_{P\Qart}(\nu,\cb')$ be reduced $k$-systems. So $R^1$ and $R^2$ have the common vertex set $[rm_1]$ partitioned into $r$ parts $U_1^1, \dots, U_r^1$ by $\Part_{R^1} = \Part_{R^2}$, where for each $i \in [r]$ part $U^1_i$ of $\Part_{R^1}$ corresponds to part $V_i$ of $\Part$. Note that since $P$ is vertex-equitable each $U_i^1$ has size $m_1$. Also note that since $R^2$ has weaker density parameters than $R^1$, any edge of $R^1$ is also an edge of $R^2$. Now, since $J$ is a $k$-complex, condition (ii) implies that  for any $p \in [k-1]$ and sets $S_i \subseteq V_i$ with $|S_i| = \lfloor pn/k \rfloor$ for $i \in [r]$ we have $|J_{p+1}[S]| \geq \beta n^{p+1}$, where $S:= \bigcup_{i \in [r]} S_i$. So the conditions of Theorem~\ref{approxalmostpackingFpartite} hold. Since $R^1$ is defined here exactly as $R$ was in the proof of Theorem~\ref{approxalmostpackingFpartite}, we have the following claim, whose proof is identical to that of Claim~\ref{matchingclaim}. 

\begin{claim} \label{matchingclaim2}
Suppose that sets $U_i \subseteq U^1_i$ for $i \in [r]$ satisfy $|U_1| = \dots = |U_r| = m \geq (1- 2k\alpha) m_1$. Let $U := \bigcup_{i \in [r]} U_i$, $R := R^1[U]$ and let $\Part_{R}$ be the restriction of $\Part_{R^1}$ to~$U$. Then for any $p \in [k-1]$ and sets $S_i \subseteq U_i$ with $|S_i| = \lfloor pm/k \rfloor$ for $i \in [r]$ we have $|R_{p+1}[S]| \geq \beta m^{p+1}/10$, where $S := \bigcup_{i \in [r]} S_i$.
\end{claim}

Since any edge $e \in J'$ has $d^F_{J'}(e) \ge d^F_J(e) - \alpha n/2$, by Lemma~\ref{reducedgraphminimumdegree} we have
\begin{equation}\label{eq:mindegRF}
\delta^F(R^1), \delta^F(R^2) \geq \left((1-k\nu^{1/2})m_1, \left(\frac{k-1}{k}-2 \alpha\right)m_1, \dots, \left(\frac{1}{k} - 2\alpha\right)m_1\right).
\end{equation}
So there are at most~$k\nu^{1/2} m_1$ vertices $i$ in each part of $\Part_{R^1}$ for which~$\{i\}$ is not an edge of~$R^1$. We can therefore delete $k\nu^{1/2} m_1$ vertices from each part to obtain $R^0 \sub R^1$ with $m_0 = (1-k\nu^{1/2})m_1$ vertices in each part and
\begin{equation}\label{eq:mindegR0F}
\delta^F(R^0) \geq \left(m_0, \left(\frac{k-1}{k}-3 \alpha\right)m_0, \dots, \left(\frac{1}{k} - 3\alpha\right)m_0\right).
\end{equation} 
By (\ref{eq:mindegR0F}) and Claim~\ref{matchingclaim2}, $R^0$ satisfies the conditions of Lemma~\ref{approxalmostpackingFpartite}, so contains an $F$-balanced matching $M$ which covers all but at most $\psi m_0$ vertices of $R^0$.
For each $i \in [r]$ let $U_i$ be the vertices in the $i$th part of $\Part_{R^1}$ covered by $M$; by uniformity of $F$ each $U_i$ has  size $m \geq (1-\psi/r)m_0 \geq (1-2\psi/r)m_1$. Let $\Part_R$ be the restriction of $\Part_{R^1}$ to $U := \bigcup_{j \in [r]} U_j$, so that the parts of $\Part_R$ are $U_i$ for $i \in [r]$. Note that $M$ is also a matching in $R^2$. Let $R$ be the restriction of $R^2$ to the vertices covered by $M$; then $(R, M)$ is a matched $\Part_R F$-partite $k$-system on $U$. 
For each edge $e \in M$, we arbitrarily relabel the clusters $W_i$ for $i \in e$ as $V_1^e, \dots, V_k^e$, and let $V^e = \bigcup_{j=1}^k V_j^e$. Then the following claim provides the partition of $V(J)$ required for the first step of the proof.

\begin{claim} \label{claim_partitionJ}
There is a partition of $V(J)$ into an exceptional set $\GG$ and sets $X$, $Y$, where $X$, $Y$ are partitioned into $X^e$, $Y^e$ for $e \in M$, and $X^e$, $Y^e$ are partitioned into sets $X^e_j$, $Y^e_j$, $j \in [k]$, such that
\begin{enumerate}[(i)]
\item $\GG$ is partitioned into sets $\GG_1,\dots,\GG_t$ with $t\le 3\psi n$, where each $\GG_i$ has
$b = k|I(F)|/r$ vertices in each $V_i$, $i \in [r]$,
\item for any sets $X_j^e \sub \Lambda_j^e \sub X_j^e \cup Y_j^e$, $j \in [k]$ with $|\Lambda_1^e| = \dots = |\Lambda_k^e|$, writing $\Lambda^e = \Lambda_1^e \cup \dots \cup \Lambda_k^e$, there is a perfect matching in $J'[\Lambda^e]$,
\item for any $j \in [r]$, $v \in V_j$, $\ib \in I(F)$ with $i_j > 0$ and $S \sub V(J) \sm X$ with $|S \cap V_i| \geq n/2 - n/6k$ for each $i$, there is an edge $e \in J[S \cup \{v\}]$ with $v \in e$ and $\ib_\Part(e) = \ib$, and
\item $|X^e_j| = n_1/2$ and $|Y^e_j| = (1/2 - \psi/r)n_1$ for each $e \in M$ and $j \in [k]$.
\end{enumerate}
\end{claim}

To prove the claim, we start by applying Proposition~\ref{findregularcomplex} and Theorem~\ref{blowup} to each $e \in M$, deleting at most $\theta n_1$ vertices from each $V^e_j$ to obtain $\hat{V}^e_1, \dots, \hat{V}^e_k$ and a $k$-partite $k$-complex $G^e$ on $\hat{V}^e := \hat{V}_1^e \cup \dots \cup \hat{V}_k^e$ such that 
\begin{enumerate}[(i)]
\item $G^e_k \sub G \sm Z \sub J_k$,
\item $G^e$ is $c$-robustly $2^k$-universal, and
\item $|G^e_k(v)| \geq d^* n^{k-1}$ for any $v \in \hat{V}^e$.
\end{enumerate}
For each $j \in [k]$ choose $X^e_j \sub \hat{V}_j^e$ of size $n_1/2$ uniformly at random and independently of all other choices. Then
Claim~\ref{claim_partitionJ} (ii) holds with probability $1-o(1)$: this follows from Lemma~\ref{randomsplitkeepsmatching}, for which condition (i) is immediate, and condition (ii) holds by the following lemma from \cite{KKMO} (it is part of Lemma 4.4), which is proved by a martingale argument. 

\begin{lemma} \label{randomsplit}
Suppose that $1/n \ll d^*, b_2, 1/k, 1/b_1 < 1$. Let $H$ be a $k$-partite $k$-graph with vertex classes $V_1, \dots, V_k$, where $n \leq |V_i| \leq b_1n$ for each $i\in [k]$. Also suppose that
$H$ has density $d(H) \ge d^*$ and that $b_2 |V_i| \leq t_i \leq |V_i|$ for each $i$.
If we choose a subset $X_i \sub V_i$ with $|X_i| = t_i$ uniformly at random
and independently for each $i$, and let $X = X_1 \cup \dots \cup X_k$, then the probability that
$H[X]$ has density $d(H[X]) > d^*/2$ is at least $1-1/n^2$.
\end{lemma}

Next, for each $e \in M$, arbitrarily delete up to $\psi n_1/r$ vertices from each $\hat{V}_j^e \sm X_j^e$ to obtain $Y_j^e$ such that $|X_j^e \cup Y_j^e| = (1-\psi/r) n_1$ for each $j \in [k]$. Note that Claim~\ref{claim_partitionJ} (iv) is then satisfied. Let $\GG$ consist of all vertices of $J$ which do not lie in some $X^e$ or $Y^e$. Then $(\GG,X,Y)$ is a partition of $V(J)$. Note that $\GG$ consists of all vertices in clusters $W_j$ which were not covered by edges of~$M$, and all vertices deleted in forming the sets $X^e$ and $Y^e$. There are at most $2\psi m_1 \cdot n_1 = 2 \psi n$ vertices of the first type, and at most $rm_1 \cdot \psi n_1/r = \psi n$ vertices of the second type, so $|\GG| \le 3\psi n$. Also, $b$ divides $|\GG \cap V_i|$ for $i \in [r]$, as $b$ divides $n$, and $M$ is $F$-balanced. (We also use the fact that each $X^e_j \cup Y^e_j$, $e \in M$, $j \in [k]$ has the same size.) We fix an arbitrary partition of $\Gamma$ into $\Gamma_1, \dots, \Gamma_t$ that satisfies Claim~\ref{claim_partitionJ} (i). It remains to satisfy Claim~\ref{claim_partitionJ} (iii). Since $|X \cap V_i| \leq n/2$ for each $i$, applying the Chernoff bound, with probability $1-o(1)$ we have
$\delta^F_\ell(J[V(J) \sm X]) \geq (k-\ell)n/3k$ for each $\ell \in [k-1]$
(we omit the calculation, which is similar to others in the paper, e.g.\ that in Theorem~\ref{almostpackingF}.) Now consider
any $j$, $v$, $\ib$ and $S$ as in Claim ~\ref{claim_partitionJ} (iii). Starting with $v$, we greedily construct an edge $e \in J[S \cup \{v\}]$ with $\ib_\Part(e) = \ib$ and $v \in e$. This is possible since $|X \cap V_i|,|S \cap V_i| > n/2 - n/6k$ for each $i$, so
$\delta^F_\ell(J[V(J) \sm X]) \ge n/3k > |V_i \sm (X \cup S)|$ for $\ell \in [k-1]$ and $i \in [r]$. This completes the proof of Claim~\ref{claim_partitionJ}. \medskip

The second step of the proof is the following claim, which states that $(R,M)$ contains small transferrals between any two vertices in any part of $\Part_R$, even after deleting the vertices of a small number of edges of $M$.

\begin{claim}\label{completetransdigraphF}
Let $R'$ and $M'$ be formed from $R$ and $M$ respectively by the deletion of the vertices covered by an $F$-balanced submatching $M_0 \sub M$ which satisfies $|M_0| \leq \alpha m$. Let $\Part_{R'}$ be the partition of $V(R')$ into $U'_1,\dots,U'_r$ obtained by restricting $\Part_R$. Then $D_C(R',M')[U'_i]$ is complete for each $i \in [r]$.
\end{claim}

To prove the claim, we consider a matched $k$-system $(R',M')$ formed in this manner and verify the conditions of Lemma~\ref{main-reduction-partiteF}. Note that each $U'_i$ has size $m'$, where $m \geq m' \geq (1-k\alpha/r)m \geq (1-2k\alpha)m_1$, and that $M'$ is $F$-balanced. Note also that if $e \in M'$, then $e \in M$, so $e \in R^1$. By Lemma~\ref{edgescloseddownwards} we therefore have $e \sm \{v\} \in R^2$ for any $v \in e$; this was the reason for using two reduced systems. Since $R'$ is the restriction of $R^2$ to the vertices of $M'$, we therefore have $e \sm \{v\} \in R'$ for any $v \in e$. Thus by Lemma~\ref{main-reduction-partiteF}, it suffices to verify the following conditions:
\begin{enumerate}[(a)]
\item $\delta^F(R') \geq \left(m', \left(\frac{k-1}{k} - \alpha' \right)m', \left(\frac{k-2}{k} - \alpha' \right) m', \dots, \left(\frac{1}{k} - \alpha' \right) m'\right)$,
\item for any $p \in [k-1]$ and sets $S_i \sub U'_i$ such that $|S_i| = pm'/k$ for each $i \in [r]$ there are at least $\beta (m')^{p+1}/10$ edges in $R'_{p+1}[S]$, where $S := \bigcup_{i \in [r]} S_i$, and
\item $L_{\Part^*_{R'}}(R'_k)$ is complete with respect to $\Part_{R'}$ for any partition~$\Part^*_{R'}$ of~$V(R')$ which refines~$\Part_{R'}$ and whose parts each have size at least $m'/k - 2\alpha' m'$.
\end{enumerate}

Condition (a) is immediate from (\ref{eq:mindegRF}), as $R'$ was formed from $R^2$ by deleting at most $2k\alpha m_1$ vertices from each part including every vertex $i$ for which $\{i\} \notin R^2$. For the same reason, since every edge of $R^1$ is also an edge of $R^2$, by Claim~\ref{matchingclaim2} we have (b). So it remains to verify condition (c). Consider any partition $\Part^*_{R'}$ of $V(R')$ refining $\Part_{R'}$ into $d$ parts $U^*_1, \dots, U^*_d$ each of size at least $m'/k - 2\alpha' m'$. Form a partition $\Part^{**}_{R^2}$ of $[rm_1]$ into $d$ parts $U^{**}_1, \dots, U^{**}_d$, where for each $i \in [r]$ the at most $2k\alpha m_1$ vertices of $[rm_1] \sm V(R')$ in the same part of $\Part_{R^2}$ as $U^*_i$ are inserted arbitrarily among the parts of $\Part^{**}_{R'}$ in that part of $\Part_{R^2}$. Then let $\Part^{\#}$ be the partition of $V(J')=V(J)$ into $d$ parts $U^{\#}_j := \bigcup_{i \in U_j^{**}} W_i$ for each $j$. Since any part of $\Part^*_{R'}$ or $\Part^{**}_{R^2}$ has size at least $m'/k - 2\alpha' m' \geq m_1/k - \mu m_1$, any part of $\Part^{\#}$ has size at least $n/k - \mu n$. Then $L^\mu_{\Part^{\#}}(J_k)$ is complete with respect to $\Part$ by assumption (iii) on $J$. Now consider any $\ib$ for which at least $\mu n^k$ edges $e \in J_k$ have $\ib_{\Part^{\#}}(e) = \ib$ (i.e. $\ib$ is in the generating set of $L^\mu_{\Part^{\#}}(J_k)$). At least $\mu n^k/2$ edges $e \in J'_k$ then have $\ib_{\Part^{\#}}(e) = \ib$, and so by Lemma~\ref{edgesofGtoedgesofR} it follows that at least $\mu m_1^k/4$ edges $e \in R^2_k$ have $\ib_{\Part^{**}_{R^2}}(e) = \ib$. Since $R'$ was formed from $R^2$ by deleting at most $2k\alpha m_1$ vertices from each part, there is at least one edge of $R'_k$ with $\ib_{\Part^*_{R'}}(e) = \ib$, so $\ib$ is in the generating set of $L_{\Part^*_{R'}}(R'_k)$. This verifies condition (c), so we have proved Claim~\ref{completetransdigraphF}. \medskip

Continuing with the proof of Theorem~\ref{perfectmatchingF}, the third step is to find a matching covering the exceptional set $\GG$. We proceed through $\GG_1,\dots,\GG_t$ in turn, at each step covering a set $\GG_i$ and using transferrals to rebalance the cluster sizes. Suppose we have chosen matchings $E_1, \dots, E_{s-1}$ for some $s \in [t]$, where $E_i$ covers the vertex set $V(E_i)$, with the following properties.
\begin{itemize}
\item[(a)] The sets $V(E_i)$ for $i \in [s-1]$ are pairwise-disjoint and have size at most $brk^2C$,
\item[(b)] $\Gamma_i \sub V(E_i) \sub \Gamma_i \cup Y$ for each $i \in [s-1]$, and
\item[(c)] $|\bigcup_{i \in [s-1]} V(E_i) \cap Y_1^e| = \dots = |\bigcup_{i \in [s-1]} V(E_i) \cap Y_k^e| \leq n_1/6k$ for any $e \in M$.
\end{itemize}
Then the following claim will enable us to continue the process.

\begin{claim}\label{greedy_matching}
There is a matching $E_s$ in $J$ satisfying properties (a), (b), and (c) with $s$ in place of $s-1$.
\end{claim}

To prove the claim, we need to first remove any clusters that have been too heavily used by $E_1,\dots,E_{s-1}$. Note that 
$$\sum_{i \in [s-1]} |V(E_i)| \le brk^2Ct \le 3brk^2C\psi n.$$
Thus there are at most $\alpha m_1/|I(F)|$ edges $e \in M$ such that $|\bigcup_{i \in [s-1]} V(E_i) \cap Y_j^e| \ge n_1/7k$ for some $j \in [k]$. Choose arbitrarily an $F$-balanced matching $M_0$ in $R$ of size at most $\alpha m_1$ which includes every such $e \in M$, and let $(R', M')$ be the matched $k$-system obtained by deleting the vertices covered by $M_0$ from both $R$ and $M$. Also let $\Part_{R'}$ be the partition of $V(R')$ into $U'_1,\dots,U'_r$ obtained by restricting $\Part_R$. Then by Claim~\ref{completetransdigraphF}, $D_C(R',M')[U'_i]$ is complete for each $i \in [r]$. Let $Y' \sub Y$ consist of all vertices of $Y$ except for those which lie in clusters deleted in forming $R'$ and those which lie in some $V(E_i)$. Then for any $\ell \in [k]$ we have 
$$|Y' \cap V_\ell| \geq |Y \cap V_\ell| - k\alpha m_1n_1 - tbrk^2C \geq n/2 - n/8k.$$
Now choose for every $x \in \Gamma_s$ an edge $e_x \in J$ with $x \in e_x$ so that $\{e_x : x \in \Gamma_s\}$ is an $F$-balanced matching in $J$. Since $|\Gamma_s \cap V_i| = b = k|I(F)|/r$ for each $i$, by uniformity of $F$ this will be accomplished by including $k$ edges $e_x$ of index $\ib$ for each $\ib \in I(F)$. By repeatedly applying Claim~\ref{claim_partitionJ}~(iii) we may choose these edges to lie within $Y' \cup \GG_s$ and to be pairwise-disjoint. The matching $\{e_x : x \in \Gamma_s\}$ then covers every vertex of $\GG_s$, but 
we may have unbalanced some clusters. To address this, consider the extra vertices $Q := \bigcup_{x \in \GG_s} e_x \sm \{x\}$ that were removed. Since each $\ib \in I(F)$ was represented $k$ times, $Q_j := Q \cap V_j$ has size $(k-1)b$ for each $j \in [r]$. Divide each $Q_j$ arbitrarily into $b$ parts of size $k-1$, and for each $\ib \in I(F)$ form $Q_{\ib}$ of size $k(k-1)$ by taking $i_j$ parts in $Q_j$, so that $(Q_{\ib}: \ib \in I(F))$ partitions $Q$ and for each $\ib \in I(F)$ the set $Q_\ib$ has index $\ib(Q_\ib) = (k-1)\ib$.

Next, for each $\ib \in I(F)$ we arbitrarily pick a `target' edge $e(\ib) \in M'$ of index $\ib$ to which we will transfer the imbalances caused by $Q_{\ib}$. Then for any $\ib \in I(F)$ and $x \in Q_{\ib}$ we may choose $i(x) \in [m_1]$ and $i'(x) \in e(\ib)$ such that $x$ is in the cluster $W_{i(x)}$, and $i(x)$ and $i'(x)$ are in the same part of $\Part_{R'}$. Furthermore, we may choose the $i'(x)$'s so that the multiset $\{i'(x): x \in  Q_{\ib}\}$ contains $k-1$ copies of each vertex of $e(\ib)$. Since $D_C(R',M')[U'_i]$ is complete for each $i \in [r]$, there is a simple $(i'(x),i(x))$-transferral of size at most $C$ in $(R', M')$. Choose such a transferral for every $x \in Q$, and let $(T,T')$ be the combination of these transferrals. To implement the transferral we need to select a matching $E^*_s = \{e^*: e \in T\}$ in $J[Y]$, whose edges correspond to the edges of $T$ (counted with multiplicity), in that $e^*$ contains one vertex in each cluster $W_i$ with $i \in e$ for each $e \in T$. We can construct $E^*_s$ greedily, using Lemma~\ref{regularrestriction}, since by Claim~\ref{claim_partitionJ}~(iv) $Y$ has $n_1/2 - \psi n_1/r$ vertices in each cluster, and we have used at most $n_1/5k$ of these in $E_1,\dots,E_{s-1}$ and edges so far chosen for $E^*_s$. 

Now we let $E_s = \{e_x\}_{x \in \Gamma_s} \cup E^*_s$. Then $|E_s| = br + |T| \le br + |Q| C \leq brkC$, so $E_s$ covers at most $brk^2C$ vertices. By construction $E_s$ is disjoint from $E_1,\dots,E_{s-1}$ and satisfies $\Gamma_s \sub V(E_s) \sub \Gamma_s \cup Y$. Furthermore, $E_s$ avoids all clusters in which $\bigcup_{i \in [s-1]} V(E_i)$ covers at least $n_1/7k$ vertices. Thus it remains to show that $|V(E_s) \cap Y_1^e| = \dots = |V(E_s) \cap Y_k^e|$ for any $e \in M$. To see this, note that the transferral $(T,T')$ was chosen so that 
$$\chi(T)-\chi(T') = \sum_{x \in Q} \chi(\{i'(x)\}) - \chi(\{i(x)\}) \in \R^{m_1}.$$ 
For each edge $e \in E_s$ we write $\chi^R(e) \in \R^{m_1}$ for the vector with $\chi^R(e)_i = |e \cap W_i|$. (Note that vertices in the exceptional set do not contribute here.) We also have $\sum_{e \in E^*_s} \chi^R(e) = \chi(T)$ and $\sum_{x \in \Gamma_s} \chi^R(e_x) = \sum_{x \in Q} \chi(\{i(x)\})$. It follows that 
$$\sum_{e \in E_s} \chi^R(e) = \chi(T') + \sum_{x \in Q} \chi(\{i'(x)\}) = \chi(T') + \sum_{\ib \in I(F)} (k-1)\chi(e(\ib)).$$ 
Since $T'$ and $T' + (k-1)\sum_{\ib \in I(F)}e(\ib)$ are multisets in $M$, this establishes the required property of $E_s$ to prove Claim~\ref{greedy_matching}. \medskip

Thus we can greedily complete the third step of covering the exceptional set $\GG$, in such a way that for each $e \in M$ we use an equal number of vertices in each $Y^e_j$, $j \in [k]$. To finish the proof of the theorem, for each $e \in M$ and $j \in [k]$ let $\Lambda^e_j$ consist of those vertices of $X^e_j \cup Y^e_j$ which are not covered by any of the sets $V(E_i)$ for $i \in [t]$. Then $X^e_j \sub \Lambda^e_j \sub X^e_j \cup Y^e_j$ and $|\Lambda^e_1| = \dots = |\Lambda^e_k|$. By Claim~\ref{claim_partitionJ} (ii), writing $\Lambda^e = \Lambda^e_1 \cup \dots \cup \Lambda^e_k$, there is a perfect matching in $J'_k[\Lambda^e]$. Combining these matchings for $e \in M$ and the matchings $E_1,\dots,E_t$ we obtain a perfect matching $M^*$ in $J_k$. Furthermore, all but at most $tbrkC \leq \gamma rn/ kD_F^2$ edges of this matching lie in $J'[\Lambda^e]$ for some $e \in M$. Since $M$ is $F$-balanced, it follows that $M^*$ $\gamma$-represents $F$.
\endproof

We shall see shortly that Theorems~\ref{newmain} and~\ref{newmainpartite} follow by a straightforward deduction from Theorem~\ref{perfectmatchingF}. However, the divisibility barriers considered there had the additional property of being transferral-free; we obtain this property using the next proposition.
 
\begin{prop} \label{obtaintransferralfree}
Let $1/n \ll \mu_1 \ll \mu \ll 1/r, 1/k$. Let $\Part$ be a partition of a set $V$ into parts $V_1, \dots, V_r$ of size $n$, and $J$ be a $k$-complex on $V$ such that 
there exists a partition $\Part_1$ of~$V(J)$ which refines~$\Part$ into parts of size at least $n/k - \mu_1 n$ for which 
$L^{\mu_1}_{\Part_1}(J_k)$ is incomplete with respect to $\Part$. Then there there exists a partition $\Part'$ of~$V(J)$ which refines~$\Part$ into parts of size at least $n/k - \mu n$ such that $L^{\mu}_{\Part'}(J_k)$ is transferral-free and incomplete with respect to~$\Part$.
\end{prop}

\proof
Introduce constants $\mu_1 \ll \dots \ll \mu_{kr} = \mu$, and repeat the following step for $t \geq 1$. If $L^{\mu_t}_{\Part_t}(J_k)$ is transferral-free, then terminate; otherwise  $L^{\mu_t}_{\Part_t}(J_k)$ contains some difference of index vectors $\ub_{i} - \ub_j$ with $i \neq j$. In this case, form a new partition $\Part_{t+1}$ of $V$ from $\Part_t$ by merging parts $V_i$ and $V_j$. Since $\Part_1$ has at most $kr$ parts, this process must terminate with a partition $\Part' = \Part_T$ for some $T \leq kr$. Observe that for any $t \geq 1$, if $L^{\mu_{t+1}}_{\Part_{t+1}}(J_k)$ is complete with respect to $\Part$ then the same must be true of $L^{\mu_{t}}_{\Part_{t}}(J_k)$, since $L^{\mu_1}_{\Part_1}(J_k)$ is incomplete with respect to $\Part$ it follows that $L^{\mu_T}_{\Part'}(J_k)$ is incomplete with respect to $\Part$. Furthermore, the fact that the process terminated with $\Part'$ implies that $L^{\mu_T}_{\Part'}(J_k)$ is transferral-free; since $L^{\mu'}_{\Part'}(J_k) \subseteq L^{\mu_T}_{\Part'}(J_k)$ this completes the proof.
\endproof

Now, to deduce Theorem~\ref{newmain}, we assume that properties 2 (Space barrier) and 3 (Divisibility barrier) do not hold, and show that property 1 (Matching) must hold. For this, introduce a new constant $\mu'$ with $\alpha \mu' \ll \mu$. We will apply Theorem~\ref{perfectmatchingF} with $r=1$ and $F$ consisting of ($k!$ copies of) the unique function $f:[k] \to [1]$, so that $\delta^F(J)=\delta^f(J)=\delta(J)$ and $b=k|I(F)|/r=k$ divides $n$. Now the conditions of Theorem~\ref{perfectmatchingF} are satisfied. Indeed, (i) holds by the minimum degree sequence, (ii) holds because there is no space barrier, and (iii) holds (with $\mu'$ in place of $\mu$) by our assumption that there is no divisibility barrier combined with Proposition~\ref{obtaintransferralfree}. Then $J_k$ contains a perfect matching, so this proves Theorem~\ref{newmain}.

Similarly, we can deduce Theorem~\ref{newmainpartite} for the case where $b = k\binom{r}{k}/r = \binom{r-1}{k-1}$ divides $n$ by taking $F$ to consist of all injective functions $f : [k] \to [r]$; then $\delta^F(J)=\delta^*(J)$ is the partite minimum degree sequence, and the fact that $M$ $\gamma$-represents $F$ implies that $M$ is $\gamma$-balanced. For the general case, write $d = \gcd(r,k)$, and note that $k/d$ divides $b$ and $n$, since we assume that $k \mid rn$. Thus we can choose $0 \le a \le b$ such that $b$ divides $n-ak/d$. By choosing a matching with one edge of each index in $I = \{\sum_{i \in [k]} \ub_{i + jd} : j \in [r/d]\}$ (where addition in the subscript is modulo $r$) we can remove $k/d$ vertices from each part. We can delete the edges of $a$ vertex-disjoint such matchings from $J$, which only slightly weakens the conditions of the theorem, so the general case follows from the case where $b \mid n$. We also note that we did not need the full strength of the degree assumption in Theorem~\ref{newmainpartite} except to obtain that $M$ is $\gamma$-balanced. For example, we could instead have taken the $F$ generated by $I$ as above, if $r>k$ is not divisible by $k$ (so that $F$ is connected).

Theorems~\ref{prestability} and~\ref{prestabilitypartite} on matchings covering all but a constant number of vertices follow in the same way from the following common generalisation.

\begin{theo} \label{n-O(1)matchingF}
Let $1/n \ll 1/\ell \ll \gamma \ll \alpha \ll \beta, \mu \ll 1/D_F, 1/r, 1/k$. Suppose $F$ is a $(k,r)$-uniform connected allocation with $|F| \leq D_F$, and $\Part$ partitions a set $V$ into sets $V_1, \dots, V_r$ each of size~$n$. Suppose that $J$ is a $\Part  F$-partite $k$-complex on $V$ such that
\begin{enumerate}[(i)]
\item $\delta^F(J) \geq \left(n, \left(\frac{k-1}{k} - \alpha \right)n, \left(\frac{k-2}{k} - \alpha \right) n, \dots, \left(\frac{1}{k} - \alpha \right) n\right)$,
\item for any $p \in [k-1]$ and sets $S_i \sub V_i$ such that $|S_i| = \lfloor pn/k \rfloor$ for each $i \in [r]$ there are at least $\beta n^k$ edges of $J'_k$ with more than $p$ vertices in $S :=  \bigcup_{i \in [r]} S_i$, and
\end{enumerate}
Then $J_k$ contains a matching which $\gamma$-represents $F$ and covers all but at most $\ell$ vertices.
\end{theo}
 
\proof
The proof is very similar to that of Theorem~\ref{perfectmatchingF}, so we just indicate the necessary modifications. Take the same hierarchy of constants as in the proof of Theorem~\ref{perfectmatchingF}, and also take $\psi \ll \gamma \ll 1/B, 1/C \ll \alpha$. We start by forming the same partition of $V(J)$ as in Claim~\ref{claim_partitionJ} (the divisibility conditions can be ensured by discarding a constant number of vertices). Instead of Claim~\ref{completetransdigraphF}, the corresponding claim here is that $(R',M')$ is $(B,C)$-irreducible with respect to $\Part_{R'}$. The proof is the same, except that we apply Lemma~\ref{gen_trans_2_partite} instead of Lemma~\ref{main-reduction-partiteF}, and there is no lattice condition to check. Next, in the analogue of Claim~\ref{greedy_matching}, we can no longer ensure property (c), but instead we maintain the property
\begin{equation}\tag{c'} ||\bigcup_{i \in [s]} V(E_i) \cap Y_j^e| - |\bigcup_{i \in [s]} V(E_i) \cap Y_{j'}^e|| \le 4B\mbox{ for all $e \in M$ and $j,j' \in [k]$.}
\end{equation} 
In constructing the sets $E_s$, we at first cover $\Gamma_s$ by $e_x$, $x \in \Gamma_s$ but do not attempt to balance the cluster sizes. Instead, whenever we have formed a set $E_s$ that would cause property (c') to fail, we remedy this using suitable $b$-fold transferrals with $b \le B$. To do this, write $$m^e_s = \sum_{j \in [k]} |\bigcup_{i \in [s]} V(E_i) \cap Y_j^e|\mbox{ and }d^{ej}_s = |\bigcup_{i \in [s]} V(E_i) \cap Y_j^e| - m^e_s/k$$ for $e \in M$, $j \in [k]$. Suppose we have formed $E_s$ which has $d^{ej}_s > 2B$ for some $e \in M$ and $j \in [k]$. Let $W$ be the cluster containing $Y_j^e$. Since each $\{e_x\}_{x \in \Gamma_s}$ uses the same number of vertices in each part $V_i$, $i \in [r]$, we can find another cluster $W'$ in the same part as $W$, which contains some $Y_{j'}^{e'}$ for which $d^{e'j'}_s < 0$. By $(B,C)$-irreducibility there is a $b$-fold $(W',W)$-transferral $(T,T')$ of size at most $C$, for some $b \le B$. To implement the transferral, we include a matching in $E^*_s$ whose edges correspond to the edges of $T$, using the same procedure as described in the proof of Theorem~\ref{perfectmatchingF}. Similarly, if $d^{ej}_s < -2B$ for some $e \in M$ and $j \in [k]$, then we can find another cluster $W'$ in the same part as $W$ which contains some $Y_{j'}^{e'}$ for which $d^{e'j'}_s > 0$. We then implement a $b$-fold $(W,W')$-transferral $(T,T')$ of size at most $C$, for some $b \le B$, in the same way as before. By repeating this process while there is any $d^{ej}_s > 2B$ for some $e \in M$ and $j \in [k]$, we can construct $E^*_s$ which when added to $E_s$ satisfies property (c'). Thus we can greedily cover the exceptional set $\GG$, in such a way that for each $e \in M$ and $j, j' \in [k]$ the number of vertices used from $Y^e_j$ and $Y^e_{j'}$ differ by at most $4B$. Now we discard at most $4B$ vertices arbitrarily from each $Y^e_j$ to balance the cluster sizes, so that Claim~\ref{claim_partitionJ} (ii) gives a perfect matching on the remaining sets. Combining all the matchings, we have covered all vertices, except for some number that is bounded by a constant independent of $n$.
\endproof

\chapter{Packing Tetrahedra} \label{sec:tetra} 

\begin{figure}[t]
\centering
\psfrag{1}{$V_1$}
\psfrag{2}{$V_2$}
\psfrag{3}{$V_3$}
\psfrag{4}{$V_4$}
\psfrag{V}{}
\includegraphics[width=8cm]{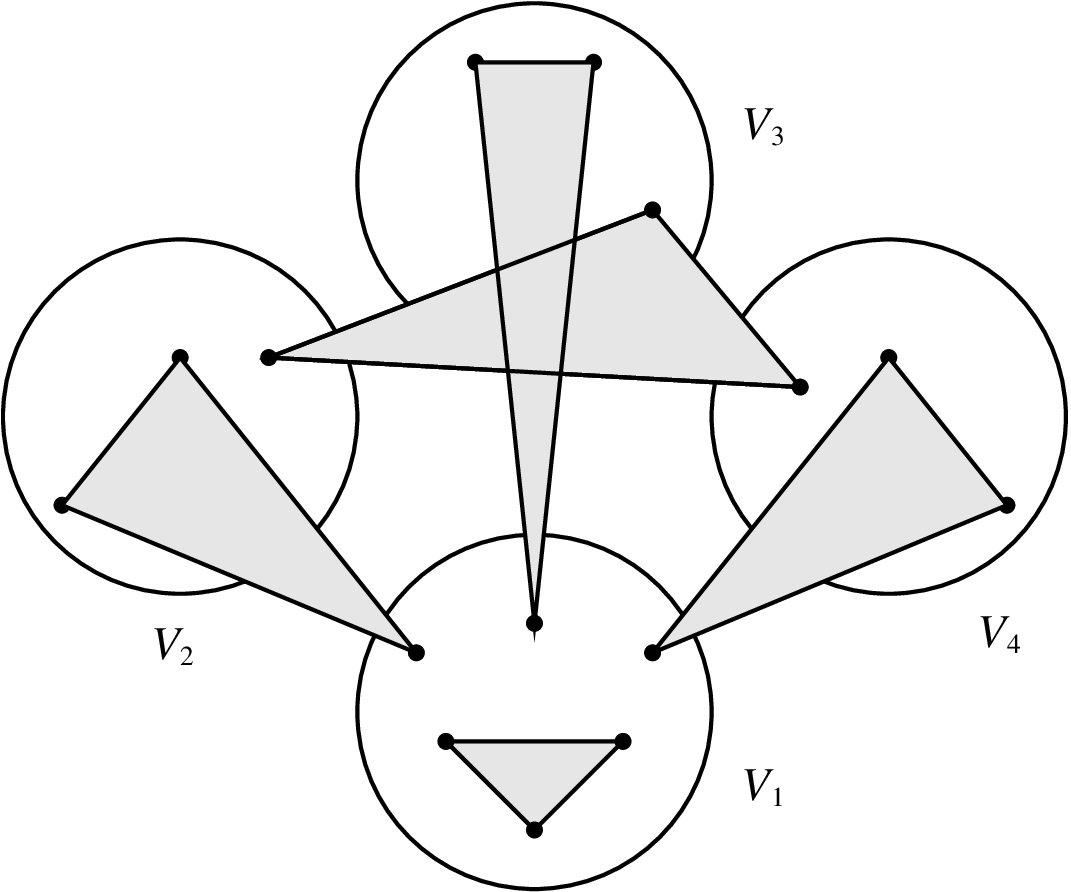} 
\caption{The 3-graph $G$ formed in Proposition~\ref{tetrahedrapackextrex} consists of all edges except for those of the types shown above. Note that any copy of $K^3_4$ in $G$ must have an even number of vertices in $V_1$.}
\label{fig:tetraex}
\end{figure}

In this chapter we prove Theorem~\ref{TETRAPACK}, which determines precisely the codegree threshold for a perfect tetrahedron packing in a $3$-graph $G$ on $n$ vertices, where $4 \mid n$ and $n$ is sufficiently large. The theorem states that if $8 \mid n$ and $\delta(G) \ge 3n/4-2$, or if $8 \nmid n$ and $\delta(G) \ge 3n/4-1$, then $G$ contains a perfect $K^3_4$-packing. We start with a construction due to Pikhurko \cite{P} showing that this minimum degree bound is best possible; the complement of this construction is illustrated in Figure~\ref{fig:tetraex}.

\medskip

\begin{prop}~\label{tetrahedrapackextrex}
For any $n$ divisible by 4 there exists a 3-graph $G$ on $n$ vertices with
$$\delta(G) \geq
\begin{cases}
3n/4 - 3 & \textrm{ if $8 \mid n$}\\ 3n/4 - 2 & \textrm{ otherwise,}
\end{cases}$$
which does not contain a perfect $K^3_4$-packing.
\end{prop}

\proof
Partition a set $V$ of $n$ vertices into parts $V_1, V_2, V_3, V_4$ such that $|V_1|$ is odd and the set sizes $|V_j|$ are as equal as possible. So $|V_j| \leq n/4 + 1$ for every $j \in [4]$, and if $8 \nmid n$ then $|V_j| = n/4$ for each $n$. Let $G$ be the $3$-graph on vertex set $V$ whose edges are all $e \in \binom{V}{3}$ \emph{except} for those $e$ such that
\begin{enumerate}[(i)]
\item $|e \cap V_1| = 3$,
\item $|e \cap V_1| = 1$ and $|e \cap V_i| = 2$ for some $i \in \{2,3,4\}$, or
\item $|e \cap V_2| = |e \cap V_3| = |e \cap V_4| = 1$.
\end{enumerate}
To calculate $\delta(G)$, consider any pair of vertices $xy$. We claim that there is precisely one class $V_i$, $i \in [4]$ such that $xyz$ is not an edge for $z \in V_i$. For this we just check all cases: if $x,y$ are in the same class then $i=1$, if $x \in V_1$ and $y \in V_j$ with $j \ne 1$ then $i=j$, and if $x,y$ are in different parts among $\{V_2,V_3,V_4\}$ then $V_i$ is the third of these parts. Thus $\delta(G) = n - 2 - \max\{|V_1|, |V_2|, |V_3|, |V_4|\}$, which has value $3n/4-3$ if $8 \mid n$ and $3n/4 - 2$ otherwise. Furthermore, we claim that any copy $X$ of $K^3_4$ in $G$ must have an even number of vertices in $V_1$. To see this, note first that we cannot have $|X \cap V_1|=3$, since there are no edges contained in $V_1$. Now suppose that $|X \cap V_1|=1$. Of the remaining $3$ vertices, either $2$ lie in the same part among $\{V_2,V_3,V_4\}$, or all $3$ parts are represented; either way some triple in $X$ is not an edge of $G$. Therefore $|X \cap V_1|$ is even. Since $|V_1|$ is odd, we conclude that $G$ does not contain a perfect $K^3_4$-packing.
\endproof

Now we start on the proof of Theorem~\ref{TETRAPACK}, which proceeds as follows. We will apply our results on perfect matchings in simplicial complexes to find a perfect matching in the clique $4$-complex $J = J_4(G)$. Note that this is precisely a perfect tetrahedron packing in $G$. We shall see in the next section that $J$ satisfies the usual minimum degree assumption, and so we consider the possible space and divisibility barriers to such a perfect matching. One such barrier is the divisibility barrier in the $4$-graph constructed in Proposition~\ref{tetrahedrapackextrex}. Indeed, in this construction every edge of $J_4$ (that is, every tetrahedron) has even intersection with $V_1$. In fact, we shall see that this is the only space or divisibility barrier to a perfect matching in $J_4$ (that is, a tetrahedron packing in $G$). 

In Section~\ref{tetraspace}, we show that there is no space barrier to a perfect matching in $G$. Furthermore, this is true even under the weaker assumption $\delta(G) \ge (3/4-c)n$ for some absolute constant $c>0$. In particular, it then follows from Theorem~\ref{prestability} that $G$ contains a tetrahedron packing that covers all but at most $\ell$ vertices, for some absolute constant $\ell$, even under this weaker assumption. Thus we demonstrate that the threshold for covering all but a constant number of vertices is asymptotically different from that for covering all vertices.

Next, in Section~\ref{tetraindex} we prove some results on how the densities of edges in a $k$-complex $J$ of different index vectors are related to one another. Indeed, for a fixed partition of $V(J)$ we define the density $d_\ib(G)$ to be the proportion of possible edges of $G$ of index $\ib$ which are in fact actual edges of $G$. Then, for example, if $V(G)$ is partitioned into parts $V_1, V_2$ and $V_3$, we should be able to say that if $d_{(1, 1, 2)}$ is large, then $d_{(1, 1, 1)}$, $d_{(0, 1, 2)}$ and $d_{(1, 0, 2)}$ should each be large (since $G$ is a $k$-complex). Similarly, if $d_{(1, 1, 1)}$ is large, and $G$ has high codegree, then at least one of $d_{(2, 1, 1)}$, $d_{(1, 2, 1)}$ and $d_{(1, 1, 2)}$ should be large. In Section~\ref{tetraindex} we prove general forms of these and other arguments, which we make extensive use of in Sections~\ref{tetratwoparts} and~\ref{tetramoreparts}.

We then proceed to the most difficult part of the proof, namely the analysis of divisibility barriers, for which we must consider Tur\'an-type problems for $K^3_4$ (fortunately we can solve these ones, unlike Tur\'an's famous conjecture!) If we do not have a perfect matching in $J_4$, then having shown that space barriers are impossible, we conclude from Theorem~\ref{newmain} that there must be a divisibility barrier. In Section~\ref{tetratwoparts} we consider divisibility barriers with two parts, and then in Section~\ref{tetramoreparts} we consider divisibility barriers with more than two parts. Taken together, in these sections we show that all but one of the potential divisibility barriers are impossible, and furthermore that the possible divisibility barrier implies a stability result, that $G$ is structurally close to the construction described in Proposition~\ref{tetrahedrapackextrex}. This sets up an application of the stability method in Section~\ref{tetraproof}: either $G$ is not structurally close to the construction, in which case Theorem~\ref{newmain} gives a perfect matching in $J_4$, or $G$ is structurally close to the construction, in which case we give a separate argument exploiting this structure to see that again we have a perfect matching in $J_4$, that is, a perfect tetrahedron packing in~$G$.

\section{Packing to within a constant.} \label{tetraspace}

As a prelude to our result on perfect tetrahedron packings, we prove the following result that demonstrates that there is a different threshold for packing to within a constant; the ingredients of its proof will also be used for the main result.

\begin{theo} \label{almosttetra}
Suppose $1/n \ll 1/\ell \ll c \ll 1$ and $G$ is a $3$-graph on $n$ vertices with $\delta(G) \ge (3/4-c)n$.
Then $G$ contains a tetrahedron packing that covers all but at most $\ell$ vertices.
\end{theo}

A reformulation of the desired conclusion is that there is a matching covering all but at most $\ell$ vertices in $J_4$, where $J=J_4(G)$ is the clique $4$-complex. First we express the minimum degree sequence of $J$ in terms of the minimum degree of $G$. Since $J_i$ is complete for $i \le 2$, we have $\delta_0(J)= n$, $\delta_1(J)=n-1$ and $\delta_2(J)=\delta(G)$. For $\delta_3(J)$, consider any $e = xyz \in J_3 = G$. For each of the pairs $xy$, $xz$, $yz$, there are at most $(n-3)-(\delta(G)-1)$ vertices $a \in V(G) \sm \{x,y,z\}$ that do not form an edge of $G$ with that pair. Thus $xyza$ is a tetrahedron in $G$ for at least $(n-3)-3(n-\delta(G)-2) = 3\delta(G)-2n+3$ vertices $a$. Therefore
\begin{equation}\label{eq:4complexdeg}
\delta(J_4(G)) \ge (n,n-1,\delta(G),3\delta(G)-2n+3).
\end{equation}

With our current assumption $\delta(G) \ge (3/4-c)n$, this gives $\delta(J) \ge (n,\,n-1,\,(3/4-c)n,\,(1/4-3c)n+3)$, so our usual minimum degree sequence assumption is satisfied; in fact, the bounds on $\delta_1(J)$ and $\delta_2(J)$ happen to be much stronger than necessary. Next, we need to consider space barriers, of which there are $3$ possible types, corresponding to sets $S$ of size $n/4$, $n/2$ or $3n/4$. The first two cases are covered by the next proposition, which is a straightforward application of the minimum degree sequence.

\begin{prop}\label{nospacebar}
Suppose $c \ll 1$ and $G$ is a $3$-graph on $n$ vertices with $\delta(G) \ge (3/4-c)n$. Then
\begin{enumerate}[(i)]
\item for any $S \sub V(J)$ of size $n/4$, the number of edges of $J_4$ with at least $2$ vertices in $S$ is at least $n^4/3000$, i.e.\ $J_4$ is not $1/3000$-contained in $J(S,1)_4$, and
\item for any $S \sub V(J)$ of size $n/2$, the number of edges of $J_4$ with at least $3$ vertices in $S$ is at least $n^4/2000$, i.e.\ $J_4$ is not $1/2000$-contained in $J(S,2)_4$.
\end{enumerate}
\end{prop}

\proof
For (i), we construct edges $v_1v_2v_3v_4$ with $v_1,v_2 \in S$. There are $n/4$ choices for $v_1$, $n/4-1$ choices for $v_2$, at least $(3/4-c)n$ choices for $v_3$, and at least $(1/4-3c)n$ choices for $v_4$.  Each such edge is counted at most $24$ times, so the bound follows. Similarly, for (ii) we construct edges $v_1v_2v_3v_4$ with $v_1,v_2,v_3 \in S$. There are $n/2$ choices for $v_1$, $n/2-1$ choices for $v_2$, at least $|S| - (n-\delta(G)) \ge (1/4-c)n$ choices for $v_3$,  and at least $(1/4-3c)n$ choices for $v_4$, so again the bound follows.
\endproof

It remains to consider the case $|S|=3n/4$, for which we need to show that many edges of $J_4$ (that is, many tetrahedra in $G$) are contained in $S$. This follows from known bounds on the Tur\'an density of $K^3_4$ and `supersaturation'. Here we only quote what suffices for our purposes, so we refer the reader to the survey \cite{K2} for more information. Given an $r$-graph $F$, the Tur\'an number $\text{ex}(n,F)$ is the maximum number of edges in an $F$-free $r$-graph on $n$ vertices. The Tur\'an density is $\pi(F) := \lim_{n \to \infty} \binom{n}{r}^{-1} \text{ex}(n,F)$. We use the following result of Chung and Lu \cite{CL}.

\begin{theo} \label{partialturan}
$\pi(K^3_4) \le (3+\sqrt{17})/12 (\approx 0.5936)$.
\end{theo}

We combine this with the following supersaturation result of Erd\H{o}s and Simonovits \cite{ES}.

\begin{theo} \label{supersat}
For any $r$-graph $F$ and $a>0$ there are $b,n_0 > 0$ so that if $G$ is an $r$-graph on $n>n_0$ vertices with $e(G) > (\pi(F)+a)\binom{n}{r}$ then $G$ contains at least $bn^{|V(F)|}$ copies of $F$.
\end{theo}

Given a set $S \sub V(J)$ of size $sn$, any pair $xy$ in $S$ is contained in at least $\delta(G)-(1-s)n$ edges of $G[S]$. Under the assumption $\delta(G) \ge (3/4-c)n$, we get $|G[S]| \ge \frac{1}{3}\binom{|S|}{2}(s-1/4-c)n \ge (1 - 1/4s - c/s)\binom{|S|}{3}$. Note that if $s \ge 5/8$ then $1 - 1/4s \ge 3/5$, so $1 - 1/4s - c/s > \pi(K^3_4) + 1/200$ for small $c$. Then Theorems~\ref{partialturan} and~\ref{supersat} imply the following result.

\begin{prop}\label{largesetstetra}
Suppose $1/n \ll b,c \ll 1$ and $G$ is a $3$-graph on $n$ vertices with $\delta(G) \ge (3/4-c)n$. Then for any $S \sub V(J)$ with $|S| \ge 5n/8$, the number of edges of $J_4$ contained in $S$ is at least $bn^4$. In particular, if $|S|=3n/4$ then $J_4$ is not $b$-contained in $J(S,3)_4$.
\end{prop}

Now Theorem~\ref{almosttetra} follows by applying Theorem~\ref{prestability} to the clique $4$-complex $J_4(G)$. Indeed, the degree sequence assumption holds by (\ref{eq:4complexdeg}), and there are no space barriers by Propositions~\ref{nospacebar} and~\ref{largesetstetra}. Thus there is a matching covering all but at most $\ell$ vertices in $J_4(G)_4$, as required.

\section{Properties of index vectors.} \label{tetraindex}

Now we need some notation and simple properties of index vectors. Let $X$ be a set which is partitioned into $r$ non-empty parts $X_1, \dots, X_r$. We let $K(X)$ denote the complete $k$-complex on $X$, where any $e \in \binom{X}{\leq k}$ is an edge of $K(X)$. Recall that the index vector $\ib(e) \in \Z^r$ of a set $e \sub X$ has coordinates $i(e)_j = |e \cap X_j|$. If $\ib(e)=\ib$, then the possible index vectors of subsets of $e$ are precisely those $\ib' \in \Z^r$ with $\0 \le \ib' \le \ib$, where as usual $\le$ is the pointwise partial order on vectors. Also recall that if $J$ is a $k$-system on $X$, then $J_\ib$ denotes the set of edges in $J$ with index $\ib$. We define the \emph{density of~$J$ at $\ib$} as
$$d_\ib(J) := \frac{|J_\ib|}{|K(X)_\ib|}.$$
Note that this notation should not be confused with that for the relative density of a complex used in Chapter~\ref{sec:regularity} (which we will not use in this chapter). Similarly to the notation used for complexes, we let $J^*_\ib$ denote the set of $e \in K(X)_\ib$ such that $e' \in J$ for every strict subset $e' \subset e$. We also recall that we write $|\ib|=\sum_{j\in [r]} i_j$, and that $\unit_i$ denotes the standard basis vector with $i$th coordinate $1$ and all other coordinates $0$. Given an index vector $\ib = (i_1,\dots,i_r)$, we write $\pl\ib$ for the multiset that contains $i_j$ copies of $\ib-\unit_j$ for each $j \in [r]$. The next proposition sets out some useful properties linking the densities of different indices.

\begin{prop} \label{lymetc}
Let $J$ be a $k$-complex on a set $X$ partitioned as $(X_1, \dots, X_r)$. Then
\begin{enumerate}[(i)]
\item for any $\ib, \ib'$ with $\ib' \leq \ib$ we have $d_{\ib'}(J) \geq d_{\ib}(J)$,
\item for any $\ib$ we have $\sum_{j \in [r]} d_{\ib+\unit_j}(J)|X_j| \geq \delta_{|\ib|}(J) d_\ib(J)$,
\item if $J_\ib = J^*_\ib$ for some $\ib = (i_1, \dots, i_r)$, then
$d_\ib(J) \geq 1 - \sum_{j \in [r]} i_j (1-d_{\ib - \unit_j}(J))$, and
\item for any $\ib$ we have $\prod_{\ib' \in \pl\ib} d_{\ib'}(J) \ge d_{\ib}(J)^{|\ib|-1} + O(1/|X|)$.
\end{enumerate}
\end{prop}

\proof
For (i), it is sufficient to consider the case when $|\ib'| = |\ib|-1$. We can write $\ib' = \ib - \unit_j$ for some $j \in [r]$ such that $i_j \geq 1$. The required density inequality is a variant of the Local LYM inequality (see e.g. \cite[Theorem 3.3]{B}). We briefly give the standard double-counting argument, which is as follows. Consider the pairs $(e',e)$ with $e' \in J_{\ib'}$, $e \in J_\ib$ and $e' \sub e$. The number of such pairs is at least $i_j|J_\ib|$, and at most $(|X_j|-i_j+1)|J_{\ib'}|$. It follows that 
$$\frac{|J_{\ib'}|}{|J_\ib|} \ge \frac{i_j}{(|X_j|-i_j+1)} =  \left.\binom{|X_j|}{i_j-1}\middle/\binom{|X_j|}{i_j}\right. = \frac{|K(X)_{\ib'}|}{|K(X)_\ib|},$$ which gives (i).

For (ii), we consider the pairs $(e,v)$ with $e \in J_\ib$ and $e \cup \{v\} \in J$. The number of such pairs is at least $\delta_{|\ib|}(J) |J_\ib|$, and at most $\sum_{j \in [r]} (i_j+1)|J_{\ib+\unit_j}|$. This gives the stated inequality, using $|J_\ib|=d_\ib(J)|K(X)_\ib|$,  $|J_{\ib+\unit_j}|=d_{\ib+\unit_j}(J)|K(X)_{\ib+\unit_j}|$, and $|X_j||K(X)_\ib| = (i_j+1)|K(X)_{\ib+\unit_j}|$.

For (iii), choose $e$ in $K(X)_\ib$ uniformly at random. For each $e' \subset e$ with $|e'|=|e|-1$, writing $\ib(e')=\ib-\unit_j$,
the probability that $e'$ is \emph{not} an edge of $J$ is $1 - d_{\ib - \unit_j}(J)$. Thus the probability that every strict subset $e' \subset e$ is an edge of $J$ is at least $1 - \sum_{j \in [r]} i_j (1-d_{\ib - \unit_j}(J))$. Since $J_\ib = J_\ib^*$, this event implies that $e \in J$, so we have the stated density.

For (iv), we first consider the multipartite case, when each $i_j$ is $0$ or $1$, and let $I = \{j: i_j=1\}$. We apply Shearer's Lemma (see \cite{CGFS}), which states that if $\mc{A}$ and $\mc{F}$ are families of subsets of $S$ such that each element of $S$ is contained in at least $\ell$ elements of $\mc{A}$, and $\mc{F}_A = \{A \cap F: F \in \mc{F}\}$ for $A \in \mc{A}$, then $\prod_{A \in \mc{A}} |\mc{F}_A| \ge |\mc{F}|^\ell$. Here we take $S=X_I$, $\mc{F}=J_\ib$, and $\mc{A} = \{X_{I \sm j}\}_{j \in I}$, so $\ell=|\ib|-1$. Then $\prod_{j \in I} |J_{\ib-\unit_j}| \ge |J_\ib|^{|\ib|-1}$, which gives the result (even without the $O(1/|X|)$ term). For the general case we apply the multipartite case to an auxiliary $k$-partite $k$-complex $H$, whose parts $X'_1,\dots,X'_{|\ib|}$ consist of $i_j$ copies of $X_j$ for $j \in [r]$, and whose edges are all $k$-partite copies of each edge in $\bigcup_{\ib' \leq \ib} J_{\ib'}$ (so, for example, each edge of $J_\ib$ corresponds to $\prod_{j \in [r]} i_j!$ of $H_k$). Then $d_{\1}(H) = d_\ib(J) + O(1/|X|)$ and 
$\{d_{\1 - \ub_j}(H)\}_{j \in [|\ib|]}$ has $i_j$ copies of $d_{\ib-\unit_j}(J) + O(1/|X|)$ for $j\in [r]$, so the result follows.
\endproof

\section{Divisibility barriers with two parts.} \label{tetratwoparts}

For the perfect packing result, we need to consider divisibility barriers. We start with the case of a bipartition, for which there are two natural candidates for a barrier: (i) almost all tetrahedra are `even', i.e.\ have an even number of vertices in both parts, and (ii) almost all tetrahedra are `odd', i.e.\ have an odd number of vertices in both parts. We start by showing that these are indeed the only two possibilities. Note that here, and for most of the proof, we only assume the weaker degree hypothesis used in the previous section; the exact degree condition is only used in Lemma~\ref{obtainoddtetrahedron}.

\begin{lemma}\label{2split}
Suppose $1/n \ll \mu,c \ll 1$ and $G$ is a $3$-graph on $n$ vertices with $\delta(G) \ge (3/4-c)n$. Write $J=J_4(G)$ for the clique $4$-complex of $G$. Suppose $\Part$ is a partition of $V$ into two parts $V_1,V_2$ of size at least $n/4 - \mu n$ such that
$L^\mu_\Part(J_4)$ is incomplete. Then either (i) there are fewer than $3\mu n^4$ odd tetrahedra, or (ii) there are fewer than $3\mu n^4$ even tetrahedra.
\end{lemma}

\proof
Suppose for a contradiction that there are at least $3\mu n^4$ odd tetrahedra and at least $3\mu n^4$ odd tetrahedra. Then $L:=L^\mu_\Part(J_4)$ contains at least one of $(1,3)$ and $(3,1)$, and also at least one of $(4,0)$, $(2,2)$ and $(0,4)$.
By (\ref{eq:4complexdeg}), we have $\delta(J) \ge (n,n-1,(3/4-c)n,(1/4-3c)n)$. Next we apply Proposition~\ref{lymetc}(ii) to $(1,1)$, which gives $$d_{(2,1)}(J)|V_1| + d_{(1,2)}(J)|V_2| \ge (3/4-c)n,$$ and to $(2,0)$, which gives $$d_{(3,0)}(J)|V_1| + d_{(2,1)}(J)|V_2| \ge (3/4-c)n.$$ Summing the inequalities we deduce $$d_{(2,1)}(J)n \ge (3/2-2c)n - |V_2| - |V_1| = (1/2-2c)n.$$ Now applying Proposition~\ref{lymetc}(ii) to $(2,1)$ we get $$d_{(3,1)}(J)|V_1| + d_{(2,2)}(J)|V_2| \ge (1/2-2c)(1/4-3c)n,$$ so $L$ contains at least one of $(3,1)$ and $(2,2)$. Similarly, $L$ contains at least one of $(1,3)$ and $(2,2)$. Now we can deduce that $(-1,1) \in L$. Indeed, if $(2,2) \in L$ this holds since $(-1,1) = (1,3) - (2,2) = (2,2) - (3,1)$ and $L$ contains at least one of $(1,3)$ and $(3,1)$. On the other hand, if
$(2,2) \notin L$ then $L$ contains both $(3,1)$ and $(1,3)$. Also, $L$ contains at least one of $(4,0)$ and $(0,4)$. Since $(-1,1) = (3,1) - (4,0) = (0,4) - (1,3)$, again we get $(-1,1) \in L$. However, this contradicts the assumption that $L$ is incomplete.
\endproof

Next we show that it is impossible to avoid even tetrahedra.

\begin{lemma} \label{noeventetrahedron}
Suppose $1/n \ll b,c \ll 1$ and $G$ is a $3$-graph on $n$ vertices with $\delta(G) \ge (3/4-c)n$. Then for any bipartition of $V(G)$, there are at least $bn^4$ copies of $K^3_4$ with an even number of vertices in both parts.
\end{lemma}

\proof
Consider any partition $(A,B)$ of $V(G)$. Form a digraph $D$ on vertex set $V(D) := \binom{A}{2} \cup \binom{B}{2}$, where there is an edge from $P \in V(D)$ to $Q \in V(D)$ precisely if $P \cup \{q\}$ is an edge of $G$ for each $q \in Q$. Then there is a one-to-one correspondence between $2$-cycles $\{PQ,QP\}$ in $D$ and copies of $K^3_4$ in $G$ with an even number of vertices in $A$. We write $|A| = an$, where without loss of generality we have $0 \le a \le 1/2$. We can also assume that $a \ge 3/8$, as otherwise $|B| \ge 5n/8$, so the desired conclusion follows from Proposition~\ref{largesetstetra}. Then by convexity,
\[|V(D)| = \binom{|A|}{2}+\binom{|B|}{2} \leq \binom{3n/8}{2}+\binom{5n/8}{2} \leq \frac{17n^2}{64}.\]
Also, for any $P \in V(D)$, we can write $d^+_D(P) = \binom{d_A(P)}{2} + \binom{d_A(B)}{2}$, where $d_X(P)$ is the number of vertices $x \in X$ such that $P \cup \{x\} \in G$, for $X \in \{A,B\}$. Since $d_A(P)+d_B(P)= \delta(G) \ge (3/4-c)n$, by convexity
$$ \delta^+(D) \geq 2 \binom{(3/4-c)n/2}{2} > \frac{9n^2}{64} - c n^2.$$
So the number of $2$-cycles in $D$ is at least $(9|V(D)|/17-cn^2)|V(D)| - \binom{|V(D)|}{2} \ge bn^4$, as required.
\endproof

Now we consider odd tetrahedra. Recall that the construction in Proposition~\ref{tetrahedrapackextrex} gave a $3$-graph $G$ on $n$ vertices partitioned into two parts $V_1$ and $V_2 \cup V_3 \cup V_4$, so that $|V_1|$ is an odd number close to $n/4$, and all tetrahedra are even with respect to the bipartition. The next lemma shows that the minimum degree condition of Theorem~\ref{TETRAPACK} is sufficient to prevent any such partition in $G$. That is, for any $k$-graph $G$ which meets the conditions of Theorem~\ref{TETRAPACK}, and any partition of $V(G)$ into two not-too-small parts, there is some tetrahedron in $H$ which intersects each part in an odd number of vertices. This is the only part of the proof of Theorem~\ref{TETRAPACK} for which the minimum degree bound is tight.

\begin{lemma} \label{obtainoddtetrahedron}
Suppose $G$ is a $3$-graph on $n$ vertices with $4 \mid n$ and $\delta(G) \ge 3n/4-1$. Let $(A,B)$ be a partition of $V(G)$ with $3n/16 < |A| < 5n/16$. Then $G$ contains a copy of $K^3_4$ with an odd number of vertices in both parts.
\end{lemma}

\proof
Consider such a partition $(A,B)$, and suppose for a contradiction that there is no odd tetrahedron. By (\ref{eq:4complexdeg}), for any edge $e \in G$, the number of vertices $v$ such that $e \cup \{v\}$ induces a copy of $K^3_4$ is at least $3\delta(G) - 2n + 3 \ge n/4$. If $e$ has an odd number of vertices in $A$, then each such vertex $v$ must also lie in $A$.
Consider any $u \in A$ and $v \in B$. Since $\delta(G) > n-|B|$ we may choose $v' \in B$ such that $e = uvv'$ is an edge of $G$. Then $e$ has an odd number of vertices in $A$, so there must be at least $n/4$ vertices $u' \in A$ such that $e \cup \{u'\}$ induces a copy of $K^3_4$. Each such $u'$ is distinct from $u$, so we must have $|A| \geq n/4 +1$.

Now consider any distinct $x,y \in B$. There are at least $\delta(G) \geq 3n/4-1$ vertices $z \in V(G)$ such that $xyz \in G$. At most $|B| - 2 \leq 3n/4 - 3$ of these lie in $B$, so we may choose some such $z \in A$. Then $e' = xyz$ has an odd number of vertices in $A$, so again there must be at least $n/4$ vertices $w \in A$ such that $e' \cup \{w\}$ induces a copy of $K^3_4$ in $G$. But $x$ and $y$ were arbitrary, so we deduce that for any $x,y \in B$, there are at most $|A|-n/4 < n/16$ vertices $z \in A$ for which $xyz$ is not an edge of $G$. Since $|A| > 3n/16$, for any edge $x'y'z'$ in $G[B]$, there is some $u' \in A$ such that $u'x'y'z'$ induces a copy of $K^3_4$ in $G$, which contradicts the assumption that there is no odd tetrahedron.
\endproof

The previous lemma will be useful when considering the extremal case, but a single odd tetrahedron does not suffice to rule out a divisibility barrier. To find many odd tetrahedra we need a result relating the triangle density to the edge density in any graph. Given a graph $G$ on $n$ vertices, the edge density of $G$ is $d(G) = e(G) \binom{n}{2}^{-1}$, and the triangle density is $\triangle(G) = N_3(G) \binom{n}{3}^{-1}$, where $N_3(G)$ is the number of triangles in $G$. In the following theorem, the first part follows from a result of Goodman \cite{G}, and the second from a result of Lov\'asz and Simonovits \cite[Theorem 2]{LS}. (Razborov~\cite{R} proved an asymptotically tight general result.)

\begin{theo} \label{triangledensity} Suppose $G$ is a graph on $n$ vertices. Then
\begin{enumerate}[(i)]
\item $\triangle(G) \ge d(G)(2d(G)-1)$, and
\item if $1/n \ll \eps \ll \eps' \ll 1$, $d(G) = 2/3 \pm \eps$ and $\triangle(G) = 2/9 \pm \eps$, then we may partition $V(G)$ into sets $V_1, V_2, V_3$ each of size at least $n/3 - \eps' n$ so that $d(G[V_i]) \leq \eps'$ for each $i \in [3]$.
\end{enumerate}
\end{theo}

Now we can prove the following result, which states that, even under the weaker degree assumption $\delta(G) \ge (3/4-c)n$, either there is no divisibility barrier for odd tetrahedra, or $G$ is close to the construction of Proposition~\ref{tetrahedrapackextrex}. This lemma plays a key role in the stability argument outlined at the start of this chapter. Indeed, to prove Theorem~\ref{TETRAPACK} we will consider first two cases: either there is no divisibility barrier, in which case Theorem~\ref{newmain} implies that $J_4(G)_4$ contains a perfect matching, or there is a divisibility barrier, in which case the structure provided by Lemma~\ref{oddtetrahedraextremalstruct}, combined with the odd tetrahedron guaranteed by Lemma~\ref{obtainoddtetrahedron}, will enable us to construct a perfect matching in $J_4(G)_4$.

Let $\Part$ partition a set $V$ into parts $V_1$, $V_2$, $V_3$, $V_4$. We say that an index vector $\ib$ with respect to $\Part$ is {\em bad} if it is not used by edges of the extremal example, i.e.\ $\ib$ is one of $(3,0,0,0)$, $(1,2,0,0)$, $(1,0,2,0)$, $(1,0,0,2)$ or $(0,1,1,1)$; otherwise $\ib$ is {\em good}.

\begin{lemma} \label{oddtetrahedraextremalstruct}
Let $1/n \ll b,c \ll \gamma \ll 1$ and $G$ be a $3$-graph on $n$ vertices with $\delta(G) \geq (3/4-c)n$. Suppose that there is some partition of $V(G)$ into parts $U_1$ and $U_2$ each of size at least $n/4 - bn$ such that at most $bn^4$ copies of $K^3_4$ in $G$ are `odd', i.e.\ have an odd number of vertices in each part. Then there is a partition of $V(G)$ into parts $V_1$, $V_2$, $V_3$, $V_4$ each of size at least $n/4 - \gamma n$, such that $V_1$ is equal to one of $U_1$ or $U_2$, and $d_\ib(G) \geq 1-\gamma$ for each good $\ib$.
\end{lemma}

\proof
Introduce new constants with $b,c \ll \alpha \ll \beta \ll \beta' \ll \gamma$.
Let $J=J_4(G)$ be the clique $4$-complex of $G$. By (\ref{eq:4complexdeg}), we have $\delta(J) \ge (n,n-1,(3/4-c)n,(1/4-3c)n)$.
Write
\begin{align*}
\tau &:= d_{(3,0)}(J) = d(G[U_1]), &  \tau' &:= d_{(0,3)}(J) = d(G[U_2]) \\
\rho &:= d_{(2,1)}(J) - \alpha, & \rho' &:= d_{(1,2)}(J) - \alpha.
\end{align*}
Let $Z$ be the set of vertices that lie in more than $\sqrt{b} n^3$ odd tetrahedra. Then $|Z| < 4\sqrt{b}n$.
Let $U'_1 = \{v \in U_1: d(G(v)[U_2])) \geq \rho'\}$. Then
$$\rho'+\alpha = d_{(1,2)}(J) = \frac{1}{|U_1|} \sum_{v \in U_1} d(G(v)[U_2])
\le \frac{|U'_1| + |U_1 \sm U'_1|\rho'}{|U_1|},$$ so $|U'_1| \ge \alpha|U_1|$. Thus we can choose $v \in U'_1 \sm Z$. We claim that
$G(v)[U_2]$ must have triangle density at most $1-\tau'+\alpha$. Otherwise, we would have at least $\alpha \binom{|U_2|}{3} > \sqrt{b} n^3$ triples in $U_2$ that are both triangles in $G(v)[U_2]$ and edges of $G[U_2]$. But these triples form odd tetrahedra with $v$, which contradicts $v \notin Z$. So $G(v)[U_2]$ has triangle density at most $1-\tau'+\alpha$, and, since $v \in U'_1$, edge density at least $\rho'$. Now Theorem~\ref{triangledensity}(i) gives
\begin{equation}\label{eq:first}
\rho'(2\rho'-1) + \tau'  - \alpha \leq 1.
\end{equation}
The same argument with the roles of $U_1$ and $U_2$ switched shows that
\begin{equation}\label{eq:second}
\rho(2\rho-1) + \tau  - \alpha \leq 1.
\end{equation}
Next, write $a = |U_1|/n$, where without loss of generality
\begin{equation} \label{eq:third}
1/4 - b \leq a \leq 1/2.
\end{equation}
Since $d_\ib(J) = 1$ for any $\ib$ with $|\ib| = 2$, by Proposition~\ref{lymetc}(ii) applied to each such $\ib$ we have
\begin{align}
\label{eq:taurho1} \tau a + (\rho + \alpha)(1-a) &\geq 3/4 - c,  \\
\label{eq:taurho2} \tau' (1-a) + (\rho' + \alpha)a &\geq 3/4 - c, \\
\label{eq:taurho3} (\rho + \alpha)a + (\rho'+ \alpha) (1-a) &\geq 3/4 - c.
\end{align}
We claim that inequalities (\ref{eq:first})-(\ref{eq:taurho3}) imply the estimates
\begin{align} \label{eq:estimates}
a &= 1/4 \pm \beta, & \tau &\leq \beta, & \rho &\geq 1 - \beta, & \rho' &= 2/3 \pm \beta, & \tau' &= 7/9 \pm \beta.
\end{align}
The proof of this claim requires some calculations, which are tiresome, but not difficult. We first consider inequality (\ref{eq:taurho1}). The left hand side is linear in $a$, so its maximum is attained at (at least) one of the extreme values $a=1/2$ and $a=1/4-b$. Thus we either have (A) $3/4-c \le \tau/2 + (\rho + \alpha)/2$ or (B) $3/4-c \le (1/4-b)\tau + (\rho + \alpha)(3/4+b)$. We will show that the first case (A) leads to a contradiction. For it implies
\begin{align*}
1 - 2c &= 2(3/4-c)-1/2 \le \tau + \rho + \alpha - 1/2 = \tau + \rho(2\rho-1) - (2\rho-1)^2/2 + \alpha \\ &\le 1 - (2\rho-1)^2/2 + 2\alpha
\end{align*}
by (\ref{eq:second}), so $(2\rho-1)^2 \le 5\alpha$. Then $\rho = 1/2 \pm \sqrt{5\alpha}$, so (\ref{eq:taurho1}) gives
\[3/4 - c \le a + (1/2 + \sqrt{6\alpha})(1-a) = 1/2 + \sqrt{6\alpha} + (1/2 - \sqrt{6\alpha})a,\]
so $a \ge 1/2 - 9\sqrt{\alpha}$. Since $a \le 1/2$, we can write $a = 1/2 \pm  9\sqrt{\alpha}$, and apply the same reasoning to inequality (\ref{eq:taurho2}). This gives $3/4 - c \le (\tau'+\rho'+\alpha)(1/2+ 9\sqrt{\alpha})$, so
\begin{align*}
1 - 2c &\le (\tau' + \rho' + \alpha)(1 + 18\sqrt{\alpha}) - 1/2 \le \tau' + \rho'(2\rho'-1) - (2\rho'-1)^2/2 + 38\sqrt{\alpha} \\ &\le 1 - (2\rho'-1)^2/2 + 39\sqrt{\alpha},
\end{align*}
by (\ref{eq:first}). Then $(2\rho'-1)^2 \le 80\sqrt{\alpha}$, so $\rho' = 1/2 \pm 10\alpha^{1/4}$. But now substituting $a,\rho,\rho'$ in (\ref{eq:taurho3}) gives a contradiction. Thus case (A) is impossible, so we must have case (B). This implies
\begin{align*}
1 - 4c &= 4(3/4-c)-2 \le \tau + 3\rho + 4\alpha - 2 = \tau + \rho(2\rho-1) - 2(\rho-1)^2 + 4\alpha \\ &\le 1 - 2(\rho-1)^2 + 5\alpha 
\end{align*}
by (\ref{eq:second}). Then $(\rho-1)^2 \le 3\alpha$, so $\rho \ge 1 - \sqrt{3\alpha}$. Then (\ref{eq:second}) gives
\[\tau \le 1 + \alpha - (1 - \sqrt{3\alpha})(1 - 2\sqrt{3\alpha}) < 6\sqrt{\alpha}.\]
Now (\ref{eq:taurho1}) gives $3/4 - c \le 6\sqrt{\alpha} + (1 + \alpha)(1-a)$, so
\[a < 1/4 + 7\sqrt{\alpha}.\]
Then (\ref{eq:taurho3}) gives $3/4-c \le (1+\alpha)(1/4 + 7\sqrt{\alpha}) + (\rho'+\alpha)(3/4+b)$, so
\[\rho' \ge 2/3 - 10\sqrt{\alpha}.\]
Now (\ref{eq:first}) implies
\[\tau' \le 1 + \alpha - (2/3 - 10\sqrt{\alpha})(1/3 - 20\sqrt{\alpha}) < 7/9 + 20\sqrt{\alpha}.\]
Next, since $\tau' \leq 1 + \alpha  - \rho'(2\rho'-1)$ by (\ref{eq:first}), substituting in (\ref{eq:taurho2}) gives
\begin{align*}
3/4 - c &\le (1 + \alpha  - \rho'(2\rho'-1))(3/4+b) + (\rho' + \alpha)(1/4 + 7\sqrt{\alpha}) \\ &\le 3/4 - \frac{3}{2}\rho'(\rho'-2/3) + 8\sqrt{\alpha},
\end{align*}
so $\rho' \le 2/3 + 9\sqrt{\alpha}$. Finally, substituting this in (\ref{eq:taurho2}) gives
\[3/4 -c \le \tau'(3/4+b) + (2/3 + 10\sqrt{\alpha})(1/4 + 7\sqrt{\alpha}),\]
so $\tau' \ge 7/9 - 10\sqrt{\alpha}$. Thus we have verified all the estimates in (\ref{eq:estimates}).

Applying (\ref{eq:estimates}) to the vertex $v \in U'_1 \sm Z$ chosen above, we see that the graph $G(v)[U_2]$ has edge density at least $2/3 - \beta$ and triangle density at most $2/9 + 2\beta$. So by Theorem~\ref{triangledensity}(ii) we may partition $U_2$ into sets $V_2$, $V_3$ and $V_4$ each of size at least $n/4 - \gamma n$ such that for each $i \in \{2,3,4\}$ we have $d(G(v)[V_i]) \leq \beta'$.  By Theorem~\ref{triangledensity}(i), the triangle density of $G(v)[U_2]$ is at least $\rho' (2\rho' -1) \geq 2/9 - 2 \beta$. Since $d(G[U_2]) = \tau' \geq 7/9 - \beta$, we see that all but at most $4\beta \binom{|U_2|}{3}$ triples in $U_2$ either form triangles in $G(v)[U_2]$ or edges in $G[U_2]$; otherwise, we would have at least $\beta \binom{|U_2|}{3} > \sqrt{b} n^3$ odd tetrahedra containing $v$, contradicting $v \notin Z$. Now write $V_1 := U_1$ and consider index vectors with respect to the partition $(V_1,V_2,V_3,V_4)$. Then for each $i \in \{2,3,4\}$, since $d(G(v)[V_i]) \leq \beta'$ we have $d_{3\unit_i}(J) = d(G[V_i]) \geq 1- \gamma$. Also, since all but at most $\beta' n^3$ triples of index $(0,1,1,1)$ are triangles of $G(v)$ we have
$d_{(0,1,1,1)}(J) < 100\beta'$. Applying  Proposition~\ref{lymetc}(ii) to $\ib=(2,0,0,0)$ we have
$$ d_{(3,0,0,0)}(J)|V_1| + d_{(2,1,0,0)}(J)|V_2| + d_{(2,0,1,0)}(J)|V_3| + d_{(2,0,0,1)}(J)|V_4| \geq (3/4 - c)n.$$
Since $d_{(3,0,0,0)}(J) = d(G[V_1]) = \tau \leq \beta$ and $|V_1|=(1/4 \pm \beta)n$, each of $d_{(2,1,0,0)}(J)$, $d_{(2,0,1,0)}(J)$, $d_{(2,0,0,1)}(J)$ must be at least $1 - \gamma$. Similarly, since $d_{(0,1,1,1)}(J) < 100\beta'$, applying Proposition~\ref{lymetc}(ii) to each of $(0,1,1,0)$, $(0,1,0,1)$ and $(0,0,1,1)$ we see that $d_\ib(J) \ge 1 - \gamma$ for all remaining good index vectors.
\endproof

\section{Divisibility barriers with more parts.} \label{tetramoreparts}

Having considered divisibility barriers with two parts, we now consider the possibility of divisibility barriers with three or more parts. Indeed, in this section we show that, if there are no divisibility barriers with two parts, then there are no divisibility barriers with more than two parts. We consider the cases of three parts and four parts  separately (there cannot be more because of the minimum size of the parts). We will repeatedly use the following observation, which is immediate from Proposition~\ref{lymetc}(ii).

\begin{prop}\label{extend}
Let $\mu \ll d$ and $\Part$ partition a set $X$ of $n$ vertices into $r$ parts. Suppose $J$ is a $k$-complex on $X$ with $d_\ib(J) \ge d$ and $\delta_{|\ib|}(J) \ge dn$ for some index vector $\ib \in \Z^{r}$. Then there is some $j \in [r]$ such that $\ib+\unit_j \in L^\mu_\Part(J_{|\ib|+1})$, i.e.\ there are at least $\mu n^{|\ib|+1}$ edges in $J$ with index vector $\ib+\unit_j$.
\end{prop}

We start by considering divisibility barriers with three parts. Recall that a lattice $L \subseteq \Z^d$ is transferral-free if it does not contain any difference $\ub_i - \ub_j$ of distinct unit vectors $\ub_i, \ub_j$.

\begin{lemma}\label{3split}
Let $1/n \ll c \ll \mu \ll 1$, and $G$ be a $3$-graph on $n$ vertices with $\delta(G) \ge (3/4-c)n$. Write $J=J_4(G)$ for the clique $4$-complex of $G$. Then there is no partition~$\Part$ of~$V(G)$ into three parts of size at least $n/4 - \mu n$ such that $L^\mu_\Part(J_4)$ is transferral-free and incomplete.
\end{lemma}

\proof
Introduce constants $\alpha, \beta, \beta'$ with $\mu \ll \alpha \ll \beta \ll \beta' \ll 1$. We need the following claim.

\begin{claim} \label{claim:3split}
Suppose $(U_1,U_2)$ is a partition of $V$ with $|U_1|=an$ and $1/4 - \alpha \le a \le 1/2 + \alpha$.
Then at least one of the following holds:
\begin{enumerate}[(i)]
\item there are at least $\mu n^4$ tetrahedra with $3$ vertices in $U_1$,
\item $a = 1/4 \pm \beta$ and $d_{(2,1)}(J) > 1 - \beta'$,
\item $a = 1/2 \pm \beta$ and $d_{(1,2)}(J) > 1 - \beta'$.
\end{enumerate}
\end{claim}

To prove the claim, suppose option~(i) does not hold, i.e.\ there are at most $\mu n^4$ tetrahedra with $3$ vertices in $U_1$. We repeat the first part of the proof of Lemma~\ref{oddtetrahedraextremalstruct} with $\mu$ in place of $b$. Define $\tau,\tau',\rho,\rho'$ as in that proof. Then there is a vertex $v \in U_2$ with $d(G(v)[U_1]) \ge \rho$ that belongs to at most $\sqrt{\mu} n^3$ odd tetrahedra (interchanging the roles of $U_1$ and $U_2$ from before). It follows that $G(v)[U_1]$ has triangle density at most $1-\tau+\alpha$, so (\ref{eq:second}) holds. Also (\ref{eq:taurho1}) and (\ref{eq:taurho3}) hold as before.

First we show that $1/4 + \beta < a < 1/2 - \beta$ leads to a contradiction. Since (\ref{eq:taurho1}) is linear in $a$, we either have~(A) $3/4 -c \le (1/4 + \beta)\tau + (3/4 - \beta)(\rho+\alpha)$ or~(B) $3/4 -c \le (1/2 - \beta)\tau + (1/2 + \beta)(\rho+\alpha)$. Consider case~(A). It implies 
\begin{align*}1 - 4c &= 4(3/4-c)-2 \le \tau + 3\rho -2 - 4\beta(\rho - \tau - \alpha/\beta) - \alpha\\& = \tau - \alpha + \rho(2\rho-1) - 2(\rho-1)^2 - 4\beta(\rho - \tau - \alpha/\beta) \\& \le 1 - 2(\rho-1)^2 - 4\beta(\rho - \tau - \alpha/\beta),
\end{align*} so $2(\rho-1)^2 \le 4c - 4\beta(\rho - \tau - \alpha/\beta)$. However, this implies $\rho > 1 - 2\sqrt{\beta}$, so $\tau < 10\sqrt{\beta}$ by (\ref{eq:second}), which contradicts $4c - 4\beta(\rho - \tau - \alpha/\beta) \geq 0$. 

Now consider case~(B). It implies 
\begin{align*}
1 - 2c &= 2(3/4-c)-1/2 \le \tau + \rho - 1/2 - 2\beta(\tau-\rho-\alpha/\beta) \\&= \tau + \rho(2\rho-1) - (2\rho-1)^2/2 - 2\beta(\tau-\rho-\alpha/\beta) \\&\le 1 + \alpha - (2\rho-1)^2/2 - 2\beta(\tau-\rho-\alpha/\beta),
\end{align*} 
so $(2\rho-1)^2 + 4\beta(\tau-\rho-\alpha/\beta) \le 4c+2\alpha$.
Then $\rho = 1/2 \pm 3\sqrt{\beta}$ and $\beta(\tau-\rho-\alpha/\beta) \le c + \alpha/2$. Now (\ref{eq:taurho1}) gives $3/4 - c \le \tau(1/2 - \beta) + (1/2 + 3\sqrt{\beta} + \alpha)(3/4-\beta)$, so $\tau \ge 3/4$. But now $\beta/5 \le \beta(\tau-\rho-\alpha/\beta) \le c + \alpha/2$, which is a contradiction.

Next suppose that $a = 1/4 \pm \beta$. By (\ref{eq:taurho1}) we have $3/4 -c \le (1/4 + \beta)\tau + (3/4 + \beta)(\rho+\alpha)$. Similarly to case~(A) we have 
$$1 - 4c = 4(3/4-c)-2 \le \tau + 3\rho - 2 + 9 \beta \le 1 - 2(\rho-1)^2 + \alpha + 9\beta,$$ 
so $d_{(2,1)}(J) = \rho > 1 - \beta'$, which is option~(ii). Finally, suppose that $a = 1/2 \pm \beta$. By (\ref{eq:taurho1}) we have $3/4 -c \le (1/2 + \beta)(\tau + \rho+\alpha)$. Similarly to case~(B) we have 
$$1 - 2c = 2(3/4-c)-1/2 \le \tau + \rho - 1/2 + 5\beta   \le 1 +\alpha - (2\rho-1)^2/2 + 5\beta,$$ so $\rho = 1/2 \pm 2\sqrt{\beta}$. Now (\ref{eq:taurho3}) gives $3/4 - c \le (1/2 + \beta)(\rho + \alpha + \rho'+ \alpha)$, so $d_{(1,2)}(J) = \rho' > 1 - \beta'$, which is option~(iii). This proves Claim~\ref{claim:3split}.

\medskip

Returning to the proof of the lemma, suppose for a contradiction we have a partition $\Part$ of~$V$ into parts $V_1,V_2,V_3$ of size at least $n/4 - \mu n$ such that $L^\mu_\Part(J_4)$ is transferral-free and incomplete. Without loss of generality $|V_1| \ge n/3$.
Also, since each of $V_1,V_2,V_3$ has size at least $n/4 - \mu n$, each has size at most $n/2 + 2\mu n$. Now recall that $\delta(J) \ge (n,n-1,(3/4-c)n,(1/4-3c)n)$ by (\ref{eq:4complexdeg}). Since $J_2$ is complete, for any $\ib$ with $|\ib|=2$, Proposition~\ref{lymetc}(ii) gives 
$$(3/4-c)n \le \sum_{j \in [3]} d_{\ib+\unit_j}(J)|V_j| \le d_{\ib+\unit_1}(J)|V_1| + n-|V_1|.$$ 
This gives 
$$(1/4+c)n \ge (1-d_{\ib+\unit_1}(J))|V_1| \ge (1-d_{\ib+\unit_1}(J))n/3,$$
so $d_{\ib+\unit_1}(J) \geq 1/5$.
Thus $d_\ib(J) \geq 1/5$ for every $\ib = (i_1,i_2,i_3)$ with $|\ib|=3$ and $i_1 \ge 1$. For each such $\ib$, Proposition~\ref{extend} gives some $j \in [3]$ with $\ib+\unit_j \in L := L^\mu_\Part(J_4)$. In particular, there is some $j \in [3]$ such that $(1,1,1)+\unit_j \in L$. Without loss of generality $j=1$ or $j=2$, as if $j=3$ we can rename $V_2$ and $V_3$ to get $j=2$.

Suppose first that $j=1$, i.e.\ $(2,1,1) \in L$. We apply Claim~\ref{claim:3split} to the partition with $U_1 = V_1$ and $U_2 = V_2 \cup V_3$. Option (i) cannot hold, as then either $(3,1,0) \in L$, so  $(3,1,0)-(2,1,1)=(1,0,-1)\in L$, or $(3,0,1) \in L$, so $(3,0,1) - (2,1,1) = (1,-1,0) \in L$, which contradicts the fact that $L$ is transferral-free. Also, option (ii) cannot hold, since $|V_1| \ge n/3$. Thus option (iii) holds, i.e.\ $|U_1|=an$ with $a = 1/2 \pm \beta$ and $d_{(1,2)}(J) > 1 - \beta'$ with respect to $(U_1,U_2)$. It follows that $d_\ib(J) > 1 - 10\beta'$ with respect to $\Part$ for any $\ib = (i_1,i_2,i_3)$ with $i_1=1$, $i_2+i_3=2$ and $i_2,i_3 \ge 0$. Now Proposition~\ref{lymetc}(ii) applied to $(0,1,1)$ gives $(3/4-c)n \le (1/2+\beta)n + (d_{(0,2,1)}(J) + d_{(0,1,2)}(J))(1/4 + 2\beta)n$, so without loss of generality $d_{(0,2,1)}(J) > 1/3$. Since $J=J_4(G)$ is a clique complex, we have $J_{(1,2,1)}=J^*_{(1,2,1)}$, so we can apply Proposition~\ref{lymetc}(iii) to get $d_{(1,2,1)}(J) \geq 1 - 30\beta' - 2/3 \ge 1/4$. Thus $(1,2,1) \in L$, so $(1,2,1) - (2,1,1) = (-1,1,0)\in L$, again contradicting the fact that $L$ is transferral-free.

Now suppose that $j=2$, i.e.\ $(1,2,1) \in L$. We apply Claim~\ref{claim:3split} to the partition with $U_1 = V_2$ and $U_2 = V_1 \cup V_3$. As in the case $j=1$, option (i) cannot hold since $L$ is transferral-free. Also, option (iii) cannot hold, since $|V_1| \ge n/3$ and $|V_3| \ge n/4 - \mu n$. Thus option (ii) holds, i.e.\ $|U_2| = an$ with $a = 1/4 \pm \beta$ and $d_{(2,1)}(J) > 1 - \beta'$ with respect to $(U_1,U_2)$. It follows that $d_{(1,2,0)}(J)$ and $d_{(0,2,1)}(J)$ are both at least $1 - 10\beta'$ (say) with respect to $\Part$. Next note that $(2,2,0) \notin L$ and $(2,1,1) \notin L$ since $L$ is transferral-free. Thus $d_{(2,2,0)}(J)$ and $d_{(2,1,1)}(J)$ are both at most $500\mu$ (say). Applying Proposition~\ref{lymetc}(ii) to $(2,1,0)$ gives $(1/4-3c)n \le d_{(3,1,0)}(J)|V_1| + 1000\mu n$, so $d_{(3,1,0)}(J) \ge 1/3$ (say). Then Proposition~\ref{lymetc}(i) gives $d_{(2,1,0)}(J) \ge 1/3$. Now choose $e = xx'yy'$ in $K(V)_{(2,2,0)}$ uniformly at random, with $x,x' \in V_1$ and $y,y' \in V_2$. Given $xx'$, let $\rho_{xx'}|V_2|$ be the number of $v_2 \in V_2$ such that $xx'v_2 \in G$. Then $\mb{E}_{xx'} \rho_{xx'} = d_{(2,1,0)}(J) \ge 1/3$. So the probability that $xx'y$ and $xx'y'$ are both edges is at least $\mb{E}_{xx'} \rho_{xx'}(\rho_{xx'}-5/n) \ge 1/10$, by Cauchy-Schwartz. Also, each of $xyy'$ and $x'yy'$ is an edge with probability at least $d_{(1,2,0)}(J) \ge 1 - 10\beta'$. Since $J=J_4(G)$ is a clique complex, we deduce $d_{(2,2,0)}(J) = \mb{P}(e \in J_{(2,2,0)}) > 1/11$. But this contradicts $(2,2,0) \notin L$.

In either case we obtain a contradiction to the existence of the divisibility barrier $\Part$.
\endproof

Now we consider divisibility barriers with four parts.

\begin{lemma}\label{4split}
Let $1/n \ll c \ll \mu \ll 1$, and $G$ be a $3$-graph on $n$ vertices with $\delta(G) \ge (3/4-c)n$. Write $J=J_4(G)$ for the clique $4$-complex of $G$. Then there is no partition $\Part$ of $V(G)$ into four parts of size at least $n/4 - \mu n$ such that $L^\mu_\Part(J_4)$ is transferral-free and incomplete.
\end{lemma}

\proof
We introduce constants $\mu',\beta,\beta'$ with $\mu \ll \mu' \ll \beta \ll \beta' \ll 1$. Suppose for a contradiction we have a partition $\Part$ of $V$ into parts $V_1,V_2,V_3,V_4$ of size at least $n/4 - \mu n$ such that $L=L^\mu_\Part(J_4)$ is incomplete and transferral-free.  Note that all parts have size at most $n/4 + 3\mu n$. Recall that $\delta(J) \ge (n,n-1,(3/4-c)n,(1/4-3c)n)$ by (\ref{eq:4complexdeg}). Now we need the following claim.

\begin{claim} \label{claim:4split} $ $
\begin{enumerate}[(i)]
\item If $\ib \in \Z^4$ with $|\ib| = 4$ and $d_\ib(J) \geq \beta$ then $d_{\ib}(J) \geq 1-\beta$.
\item If $\ib \in \Z^4$ with $|\ib| = 3$ and $d_\ib(J) \geq 2\beta$ then $d_{\ib}(J) \geq 1-\beta$.
\end{enumerate}
\end{claim}

To prove the claim, first consider any $\ib \in \Z^4$ with $|\ib| =4$ and $d_\ib(J) \geq \beta$. Note that $\ib \in L$. By Proposition~\ref{lymetc}(iv) there is some $\ib' = \ib - \unit_j \in \pl\ib$ such that $d_{\ib'}(J) \ge d_\ib(J)^{3/4} + O(1/n)$. Let $B$ be the set of edges $e' \in J_{\ib'}$ that lie in at least $\mu'n$ edges $e \in J_4$ with $\ib(e) \neq \ib$. Then there is some $j' \ne j$ for which the number of edges with index vector $\ib - \unit_j + \unit_{j'}$ is at least $|B|\mu'n/12$. This must be less than $\mu n^4$, otherwise we have $\ib - \unit_j + \unit_{j'} \in L$, so $\unit_j - \unit_{j'} \in L$, contradicting the fact that $L$ is transferral-free. It follows that $d_{\ib'}(B) \le \mu'$ (say), so $d_{\ib'}(J \sm B) \ge d_\ib(J)^{3/4} - 2\mu'$. Next we show that if $e' \in J \sm B$ then $e' \cup \{x\} \in J_\ib$ for all but $2\mu'n$ vertices $x \in V_j$. This holds because there are at least $\delta_3(J) \geq (1/4 - 3c)n$ vertices $x$ such that $e' \cup \{x\} \in J_4$, all but at most $\mu'n$ of these lie in $V_j$ (as $e' \notin B$), and all parts have size at most $n/4 + 3\mu n$. It follows that 
$$d_\ib(J) \ge (1-10\mu')d_{\ib'}(J \sm B) \ge (1-10\mu')(d_\ib(J)^{3/4} - 2\mu').$$ Since $d_\ib(J) \geq \beta$ this implies $d_{\ib}(J) \geq 1-\beta$, so we have proved (i). For (ii), observe that by Proposition~\ref{lymetc}(ii) there must be some $\ib'$ such that $|\ib'| = 4$, $\ib \leq \ib'$ and $d_{\ib'}(J) \geq (3/4-c)2\beta \ge \beta$. Then by (i) we have $d_\ib(J) \geq 1-\beta$, so by Proposition~\ref{lymetc}(i) we have $d_{\ib'}(J) \geq d_{\ib}(J) \geq 1-\beta$. This completes the proof of Claim~\ref{claim:4split}. \medskip

Returning to the proof of the lemma, we show next that $d_\ib(J) \leq 2\beta$ whenever $\ib$ is a permutation of $(3,0,0,0)$. For suppose this fails, say for $\ib=(3,0,0,0)$. Then Claim~\ref{claim:4split} gives $d_{(3,0,0,0)}(J) \ge 1-\beta$. Now Proposition~\ref{lymetc}(iii) gives $d_{(4,0,0,0)}(J) \geq 1 - 4\beta$, so $(4,0,0,0) \in L$. Next, Proposition~\ref{lymetc}(ii) for $\ib' = (2,0,0,0)$ gives $(3/4-c)n \le (n/4 + 3\mu n) \sum_{j \in [4]}  d_{\ib'+\unit_j}(J)$, so without loss of generality $d_{(2,1,0,0)}(J) \geq 2\beta$. Then Claim~\ref{claim:4split} gives $d_{(2,1,0,0)}(J) \ge 1-\beta$. Now Proposition~\ref{lymetc}(iii) gives $d_{(3,1,0,0)}(J) \geq 1 - 4\beta$, so $(3,1,0,0) \in L$. But then $(4,0,0,0) - (3,1,0,0) = (1,-1,0,0) \in L$ contradicts the fact that $L$ is transferral-free. Thus $d_\ib(J) \leq 2\beta$ whenever $\ib$ is a permutation of $(3,0,0,0)$.

Now returning to the above inequality 
$$(3/4-c)n \le (n/4 + 3\mu n) \sum_{j \in [4]} d_{\ib'+\unit_j}(J),$$ where $\ib'$ is any permutation of $(2,0,0,0)$, we see that for any $\ib''$ which is a permutation of $(2,1,0,0)$ we have $d_{\ib''}(J) \geq 2\beta$, and therefore $d_{\ib''}(J) \geq (1-\beta)$ by Claim~\ref{claim:4split}.
Also, applying Proposition~\ref{lymetc}(ii) to $(1,1,0,0)$, we see that  without loss of generality $d_{(1,1,1,0)}(J) \geq 2\beta$. Then Claim~\ref{claim:4split} gives $d_{(1,1,1,0)}(J) \geq 1-\beta$. Now each of $d_{(2,2,0,0)}(J)$ and $d_{(2,1,1,0)}(J)$ is at least $1-4\beta$ by Proposition~\ref{lymetc}(iii). It follows that $L$ contains $(2,2,0,0)$, $(2,1,1,0)$ and $(2,2,0,0) - (2,1,1,0) = (0,1,-1,0)$, again contradicting the fact that $L$ is transferral-free. Thus there is no such partition $\Part$.
\endproof

\section{The main case of Theorem~\ref{TETRAPACK}.} \label{tetraproof}

As outlined earlier, Theorem~\ref{newmain} will imply the desired result, except for $3$-graphs that are close to the extremal configuration. We start with a lemma that handles such $3$-graphs. The proof requires a bound on the minimum vertex degree that guarantees a perfect matching in a $4$-graph. We just quote the following (slight weakening of a) result of Daykin and H\"aggkvist~\cite{DH} that suffices for our purposes. (The tight bound was recently obtained by Khan \cite{Kh}.)

\medskip
\begin{theo} \label{vertexmatching}
Let $H$ be a $4$-graph on $n$ vertices, where $4 \mid n$, such that every vertex belongs to at least $\frac{3}{4}\binom{n}{3}$ edges. Then $H$ contains a perfect matching.
\end{theo}

The following lemma will be used for $3$-graphs that are close to the extremal configuration. Note that we will be able to satisfy the condition that $|V_1|$ is even by applying Lemma~\ref{obtainoddtetrahedron}.

\begin{lemma} \label{packinggivenstructure}
Suppose that $1/n \ll \gamma \ll c \ll 1$. Let $V$ be a set of $4n$ vertices partitioned into $V_1$, $V_2$, $V_3$, $V_4$ each of size at least $n - \gamma n$. Suppose $J$ is a $4$-complex on $V$ such that
\begin{enumerate}[(i)]
\item $d(J_4[V_j]) \geq 1-\gamma$ for each $j \in \{2,3,4\}$,
\item for each $\ib \in (2,1,1,0), (2,1,0,1), (2,0,1,1)$ we have $d_\ib(J) \geq 1-\gamma$,
\item for every vertex $v \in V$ at least $c n^4$ edges of $J_4$ contain $v$ and have an even number of vertices in $V_1$, and
\item $|V_1|$ is even.
\end{enumerate}
Then $J_4$ contains a perfect matching.
\end{lemma}

\proof
We say that an edge of $J_4$ is \emph{even} if it has an even number of vertices in $V_1$, and \emph{odd} otherwise.
For each $j \in \{2,3,4\}$, we say that a vertex $v \in V_j$ is \emph{good} if there are at least $(1-2\sqrt{\gamma}) \binom{n}{3}$ edges of $J_4[V_j]$ which contain $v$. By (i) at most $\sqrt{\gamma} n$ vertices of each of $V_2, V_3$ and $V_4$ are bad; here we note that each part has size at most $n + 3\gamma n$. We say that a pair $u, v \in V_1$ is \emph{good} if for each $\ib \in \{(2,0,1,1), (2,1,0,1), (2,1,1,0)\}$ there are at least $(1-2\sqrt{\gamma})n^2$ edges $e \in J_4$ with $u,v \in e$ and $\ib(e) = \ib$. For each such $\ib$, by (ii) at most $\sqrt{\gamma}\binom{n}{2}$ pairs $u,v \in V_1$ do not have this property, and so at most $3 \sqrt{\gamma} \binom{n}{2}$ pairs $u,v \in V_1$ are bad. We say that a vertex $u \in V_1$ is \emph{good} if it lies in at least $(1 - 2\gamma^{1/4})n$ good pairs $u, v \in V_1$. Then at most $3\gamma^{1/4} n$ vertices of $V_1$ are bad. In total, the number of bad vertices is at most $\gamma^{1/5} n$, say.

Now let $E$ be a maximal matching in $J$ such that every edge in $E$ has an even number of vertices in $V_1$ and contains a bad vertex. We claim that $E$ covers all bad vertices. For suppose some bad vertex $v$ is not covered by $E$. Since each edge in $E$ contains a bad vertex, at most $4\gamma^{1/5} n$ vertices are covered by $E$, so at most $4\gamma^{1/5} n^4$ edges of $J_4$ contain a vertex covered by $E$. Then by (iii) we may choose an edge $e \in J_4$ which contains $v$, has an even number of vertices in $V_1$ and doesn't contain any vertex covered by $E$, contradicting maximality of $E$. So $E$ must cover every bad vertex of $J$. We will include $E$ in our final matching and so we delete the vertices it covers. We also delete the vertices covered by another matching $E'$ of at most $16$ edges, disjoint from $E$, so as to leave parts $V'_1$, $V'_2$, $V'_3$, $V'_4$ such that $8$ divides $|V'_1|$ and $4$ divides each of $|V'_2|$, $|V'_3|$, $|V'_4|$. The edges in $E'$ will have index $(2,1,1,0)$, $(2,1,0,1)$ or $(2,0,1,1)$. Note that by (ii) we can greedily choose a matching disjoint from $E$ containing $16$ edges of each of these indices, so we only need to decide how many of these we want to include in $E'$. First we arrange that there are an even number of remaining vertices in each part. Note that $|V|$ and $|V_1|$ are even, and the number of vertices remaining in $V_1$ is still even, as we only used edges with an even number of vertices in $V_1$. Then an even number of parts have an odd number of vertices remaining (after the deletion of $E$). If there are two such parts, say $V_2$ and $V_4$, then we remedy this by including an edge of index $(2,1,0,1)$ in $E'$, so we can assume there are no such parts. Now there are an even number of parts in which the number of vertices remaining is not divisible by $4$. If $V_1$ is one of these parts we include an edge of each index $(2,1,1,0)$, $(2,1,0,1)$, $(2,0,1,1)$ in $E'$; then each part has an even number of remaining vertices and the number in $V_1$ is divisible by $4$. There may still be two parts where the number of vertices remaining is not divisible by $4$. If so, say they are $V_3$ and $V_4$, we include $2$ edges of index $(2,0,1,1)$ in $E'$. Thus we arrange that the number of vertices remaining in each part is divisible by $4$. Finally, if the number of remaining vertices in $V_1$ is not divisible by $8$ then we include $2$ more edges of each index in $E'$. Thus we obtain $E'$ with the desired properties. We delete the vertices of $E \cup E'$ and let $V'_1$, $V'_2$, $V'_3$, $V'_4$ denote the remaining parts.

Next we choose a perfect matching $M_1$ of good pairs in $V'_1$. This is possible because $|V'_1|$ is even, and any vertex in $V'_1$ is good, in that it was incident to at least $(1 - 2\gamma^{1/4})n$ good pairs in $V_1$, so is still incident to at least $n/2$ good pairs in $V'_1$. Now we greedily construct a matching in $J_4$, by considering each pair in $M_1$ in turn and choosing an edge of index $(2,1,1,0)$ containing that pair. There are at most $n/2$ pairs in $M_1$, so when we consider any pair $uv$ in $M_1$, at least $n/3$ vertices in each of $V'_2$ and $V'_3$ are still available, in that they have not already been selected when we considered some previous pair in $M_1$. Since $uv$ is good, at most $3\sqrt{\gamma}n^2 < (n/3)^2$ of these pairs do not form an edge with $uv$, so we can form an edge as required. Thus we construct a matching $E''$ that covers $V'_1$. We delete the vertices of $E''$ and let $V''_2$, $V''_3$, $V''_4$ denote the remaining parts. Note that the number of vertices used by $E''$ in each of $V'_2$ and $V'_3$ is $|V'_1|/2$, which is divisible by $4$. Since $4$ divides each of $|V'_2|$, $|V'_3|$, $|V'_4|$, it also divides each of  $|V''_2|$, $|V''_3|$, $|V''_4|$. Furthermore, for any $j \in \{2,3,4\}$ and $x \in V_j$, since $x$ is good, the number of edges of $J_4[V''_j]$ containing $x$ is at least $\binom{|V_j''|}{3} - 3\sqrt{\gamma}\binom{n}{3} > \frac{3}{4}\binom{|V_j''|}{3}$, so $J_4[V''_j]$ contains a perfect matching by Theorem~\ref{vertexmatching}.  Combining these matchings with $E$, $E'$ and $E''$ we obtain a perfect matching in~$J_4$.
\endproof
 
We can now give the proof of Theorem~\ref{TETRAPACK}, as outlined at the start of this chapter. Suppose that $G$ is a $3$-graph on $n$ vertices, where $n$ is sufficiently large and divisible by $4$. For now we assume  $\delta(G) \ge 3n/4-1$, postponing the case that $8 \mid n$ and $\delta(G) = 3n/4-2$ to the final section. We introduce constants $\alpha$, $\beta$ with $1/n \ll \alpha \ll \beta \ll 1$. We will show that $G$ has a perfect tetrahedron packing; equivalently, that $J_4$ has a perfect matching, where $J=J_4(G)$ is the clique $4$-complex of $G$. By (\ref{eq:4complexdeg}) we have $\delta(J) \ge (n,n-1,3n/4-2,n/4-3)$, so we can apply Theorem~\ref{newmain}. Supposing that $J_4$ does not contain a perfect matching, we conclude that there is a space barrier or a divisibility barrier. As in the proof of Theorem~\ref{almosttetra} there are no space barriers by Propositions~\ref{nospacebar} and~\ref{largesetstetra}, so we may choose a minimal divisibility barrier, i.e.~a partition $\Qart$ of $V(G)$ into parts of size at most $n/4 - \mu n$ such that $L_\Qart^\mu(J_4)$ is transferral-free and incomplete. We have excluded all but one possibility for $\Qart$. Indeed, $\Qart$ cannot have more than two parts by Lemmas~\ref{3split} and~\ref{4split}, so must have two parts. By Lemmas~\ref{2split} and~\ref{noeventetrahedron}, we deduce that $\Qart$ partitions $V(G)$ into parts $U$ and $V$ of size at least $n/4 - \alpha n$ such that at most $\alpha n^4$ edges of $J_4$ have an odd number of vertices in each part.

By Lemma~\ref{oddtetrahedraextremalstruct}, there is a partition $\Part$ of $V(G)$ into parts $V_1$, $V_2$, $V_3$ and $V_4$ each of size at least $n/4 - \beta n$, so that $V_1 = U$ (without loss of generality) and  $d_\ib(G) \geq 1-\beta$ for each good $\ib$ (recall that good index vectors are those that appear in the extremal example). To apply Lemma~\ref{packinggivenstructure} we also need to arrange that $|V_1|$ is even, and that every vertex is in many even edges of $J_4$, where we say that an edge of $J_4$ is odd or even according to the parity of its intersection with $V_1$. We say that a vertex is \emph{good} if it lies in fewer than $n^3/400$ odd edges of $J_4$ (and \emph{bad} otherwise). Then at most $1600 \alpha n$ vertices are bad, since $J_4$ has at most $\alpha n^4$ odd edges. We let $V'_1$ consist of all good vertices in $V_1$ and all bad vertices in $V_2 \cup V_3 \cup V_4$, let $V'_2$ consist of all good vertices in $V_2$ and all bad vertices in $V_1$, and let $V'_3$ and $V'_4$ consist of all good vertices in $V_3$ and $V_4$ respectively. Note that any vertex that has a different index with respect to the partition $\Part'$ of $V(G)$ into the $V'_j$ than with respect to $\Part$ is bad, so we have $|V'_1| = n/4 \pm 2\beta n$. So by Lemma~\ref{obtainoddtetrahedron}, there exists a copy of $K^3_4$ in $G$ (i.e. an edge of $J_4$) with an odd number of vertices in $V'_1$. If $|V'_1|$ is odd, we choose one such edge of $J_4$, and delete the vertices of this edge from $V'_1$, $V'_2$, $V'_3$, $V'_4$, and $J$. Thus we can arrange that $|V'_1|$ is even.

Having changed the index of at most $1600\alpha n$ vertices of $J$ and deleted at most four vertices, we have $|V'_i| \geq n/4 - 2\beta n$ for each $i \in [4]$, and $d_\ib(G) \geq 1-2\beta$ for each good $\ib$ with respect to the partition into $V'_1$, $V'_2$, $V'_3$, $V_4'$. It remains to show that any undeleted vertex $v$ belongs to many edges with an even number of vertices in $V'_1$. To avoid confusion, we now use the terms $X$-odd and $X$-even, where $X$ is $V_1$ or $V'_1$, to describe edges according to the parity of their intersection with $X$. First suppose that $v$ is good. Then any edge of $J_4$ containing $v$ that is $V'_1$-odd but not $V_1$-odd must contain a bad vertex, so there are at most $1600 \alpha n^3$ such edges. Since $v$ belongs to at most $n^3/400$ $V_1$-odd edges of $J_4$, it belongs to at most $n^3/400 + 1600 \alpha n^3$ $V'_1$-odd edges of $J_4$. Also, $v$ belongs to at least $n^3/200$ edges of $J_4$, by the lower bound on $\delta(J)$, so $v$ belongs to at least $n^3/500$ $V'_1$-even edges. Now suppose that $v$ is bad. Then any $V_1$-odd edge of $J_4$ containing $v$ and no other bad or deleted vertices is $V'_1$-even. The number of such edges is at least $n^3/400 - 1600 \alpha n^3 - 4n^2 \geq n^3/500$. This shows that every vertex is contained in at least $n^3/500$ $V'_1$-even edges of $J_4$.

Finally, Lemma~\ref{packinggivenstructure} implies that $J_4$ has a perfect matching.

\section{The case when $8$ divides $n$.}
It remains to consider the case when $8 \mid n$ and $\delta(G) = 3n/4-2$. We apply the same proof as in the previous case, noting that we only needed $\delta(G) \ge 3n/4-1$ if $|V'_1|$ was odd, when we used it to find a $V'_1$-odd tetrahedron. Examining the proof of Proposition~\ref{obtainoddtetrahedron}, we see that a $V'_1$-odd tetrahedron exists under the assumption $\delta(G) = 3n/4-2$, unless we have $n/4 - 2 \le |V'_1| \le n/4$. Furthermore, we only need a $V'_1$-odd tetrahedron when $|V'_1|$ is odd, so we only need to consider the exceptional case that $|V'_1| = n/4-1$ and there is no $V'_1$-odd tetrahedron.

Let $c_1$ satisfy $\beta \ll c_1 \ll 1$. Recall that $d_\ib(G) \ge 1 - 2\beta$ for all good $\ib$. Note that for each bad triple $e$ there is a part $V'_j$ such that for every $x \in V'_j \sm e$, $e \cup \{x\}$ is $V'_1$-odd, and every triple in $e \cup \{x\}$ apart from $e$ is good. (We describe triples as good or bad according to their index vectors.) If $e$ is an edge, then some such triple is not an edge. Since there is no $V'_1$-odd tetrahedron, there are at most $\beta n^3$ bad edges by Proposition~\ref{lymetc}(iii). For any $2/n < c < 1/4$, we say that a pair is {\em $c$-bad} if it is contained in at most $(3/4-2c)n$ good edges. Since $\delta(G) = 3n/4-2$, any $c$-bad pair is contained in at least $cn$ bad edges. If $\beta \ll c$ then there are at most $cn^2$ $c$-bad pairs. Since there is no $V'_1$-odd tetrahedron, any bad edge of $G$ contains at least one $1/30$-bad pair.

Without loss of generality we have $|V'_4| \ge n/4+1$. Thus for each $v_2 \in V'_2$ and $v_3 \in V'_3$ there is some $v_4 \in V'_4$ with $v_2v_3v_4 \in G$. Each such edge is bad, so contains a $1/30$-bad pair. The bad pair can only be $v_2v_3$ for at most $c_1 n^2$ such triples. Thus without loss of generality there is a vertex $v \in V_4$ such that at least $n/10$ pairs $v_3v$ with $v_3 \in V'_3$ are $1/30$-bad. For each such pair there are at least $n/30$ vertices $v_2 \in V'_2$ such that $v_2v_3v$ is a bad edge. For each such bad edge and each $v_1 \in V'_1$, some triple in $v_1v_2v_3v$ is not an edge. Since $d_{(1,1,1,0)}(G) \ge 1-2\beta$, there are at most $c_1 n$ vertices $v_1 \in V'_1$ such that $v_1v_2v_3$ is not an edge for at least $c_1 n^2$ pairs $v_2v_3$. Thus for all but at most $c_1 n$ vertices $v_1 \in V'_1$ there are at least $n^2/300 - c_1 n^2$ pairs $v_2v_3$ such that one of $v_1v_2v$ or $v_1v_3v$ is not an edge. Then $v_1v$ is $1/200$-bad, so there are at least $n/200$ vertices $v_4 \in V'_4$ for which $v_1vv_4$ is a (bad) edge.

Now for each $v_1 \in V'_1$, $v_4 \in V'_4$ such that $v_1vv_4$ is an edge and $x \in V'_2 \cup V'_3$, at least one of $v_1xv_4$, $v_1xv$, $xvv_4$ is not an edge. Since $d_{(1,1,0,1)}(G)$ and $d_{(1,1,0,1)}(G)$ are at least $1-2\beta$, there are at most $c_1 n$ vertices $v_1 \in V'_1$ such that $v_1xv_4$ is not an edge for at least $c_1 n^2$ pairs $xv_4$. For $v_1 \in V'_1$, let $N_4(v_1)$ be the set of vertices  $v_4 \in V'_4$ such that $v_1vv_4$ is an edge, and $v_1xv_4$ is an edge for all but at most $600c_1 n$ vertices $x \in V'_2 \cup V'_3$. Let $V''_1$ be the set of vertices $v_1 \in V'_1$ such that $|N_4(v_1)| \ge n/300$. Then $|V'_1 \sm V''_1| \le 2c_1n$, otherwise there are at least $c_1 n$ vertices $v_1 \in V'_1 \sm V''_1$ such that there are at least $n/200 - |N_4(v_1)| \ge n/600$ vertices $v_4 \in V'_4$ such that $v_1vv_4$ is an edge and $v_4 \notin N_4(v_1)$, so $v_1xv_4$ is not an edge for at least $600c_1n \cdot n/600 = c_1 n^2$ pairs $xv_4$, which is a contradiction.

For each $v_1 \in V''_1$, $v_4 \in N_4(v_1)$ and all but at most $600c_1 n$ vertices $x \in V'_2 \cup V'_3$, $v_1xv_4$ is an edge, so at least one of $v_1xv$ and $xvv_4$ is not an edge. Since $d(v_1v)+d(vv_4) \ge 2\delta(G) = 3n/2-4$, and $|V'_j| \ge n/4 - 2\beta n$ for each $j \in [4]$, precisely one of $v_1xv$ and $xvv_4$ is not an edge for all but at most $2000c_1 n$ vertices $x \in V'_2 \cup V'_3$. Also, $v_1xv$ and $xvv_4$ are edges for all but at most $2000c_1 n$ vertices $x \in V'_1 \cup V'_4$. Now there can be at most $c_1 n$ vertices $v_1 \in V''_1$ such that $v_1v_2v$ is an edge for at least $3000c_1 n$ vertices $v_2 \in V'_2$. Otherwise, since we have at least $n/300$ choices of $v_4 \in N_4(v_1)$, so $v_1vv_4$ is an edge, and $1000c_1 n$ choices of $v_2$ such that $v_2vv_4$ is an edge, there must be at least $c_1^2 n^3$ good triples $v_1v_2v_4$ that are not edges, contradicting $d_{(1,1,0,1)}(G) \ge 1-2\beta$. Similarly, there are at most $c_1 n$ vertices $v_1 \in V''_1$ such that $v_1v_3v$ is an edge for at least $3000c_1 n$ vertices $v_3 \in V'_3$. But then we can choose $v_1 \in V''_1$ such that $v_1xv$ is an edge for at most $6000c_1 n$ vertices $x \in V'_2 \cup V'_3$, which contradicts $\delta(G) = 3n/4-2$.

Thus the exceptional case considered here cannot occur. The rest of the proof goes through as in the previous section,
so this proves Theorem~\ref{TETRAPACK}. \qed

\section{Strong stability for perfect matchings.}
 
We conclude this chapter by applying a similar argument to that of Lemma~\ref{packinggivenstructure} to prove Theorem~\ref{pmstrongstab}.
For this we need another theorem of Daykin and H\"aggkvist~\cite{DH}, which gives a bound on the vertex degree needed to guarantee a perfect matching in a $k$-partite $k$-graph.

\begin{theo}[\cite{DH}] \label{partitevertexdegree}
Suppose that $G$ is a $k$-partite $k$-graph whose vertex classes each have size $n$, in which every vertex lies in at least $\frac{k-1}{k}n^{k-1}$ edges. Then $G$ contains a perfect matching.
\end{theo}


The proof of Theorem~\ref{pmstrongstab} begins with Theorem~\ref{pmstab}, from which we obtain a set $V_1$ whose intersection with almost all edges of $G$ has the same parity. This gives us a great deal of structural information on $G$, namely that $G$ contains almost every possible edge whose intersection with $V_1$ has the `correct' parity. Given a single edge of the opposite parity, we can then delete this edge (if necessary) to `correct' the parity of $|V_1|$, whereupon the high density of edges of `correct' parity guarantees the existence of a perfect matching in $H$ by Theorem~\ref{partitevertexdegree}. We now give the details of the proof, but first restate the theorem in a slightly different form.

\medskip \noindent 
{\bf Theorem~\ref{pmstrongstab}.}
\emph{Suppose that $1/n \ll c \ll 1/k$, $k \ge 3$ and $k \mid n$, and let $G$ be a $k$-graph on $n$ vertices with $\delta(G) \ge (1/2-c)n$. Then $G$ does not contain a perfect matching if and only if there is a partition of $V(G)$ into parts $V_1,V_2$ of size at least $\delta(G)$ and $a \in \{0,1\}$ so that $|V_1| \neq an/k$ mod $2$ and $|e \cap V_1|=a$ mod $2$ for all edges $e$ of $G$.}

\proof
First observe that if such a partition exists there can be no perfect matching in~$G$. Indeed, the $n/k$ edges in such a matching would each have $|e \cap V_1|=a$ mod $2$, giving $|V_1| = an/k$ mod $2$, a contradiction. So assume that $G$ contains no perfect matching; to complete the proof it suffices to show that $G$ admits a partition as described.

Introduce new constants $b, b'$ with $1/n \ll b' \ll b \ll c \ll 1/k$. By Theorem~\ref{pmstab}, there is a partition of $V(G)$ into parts $V_1,V_2$ of size at least $\delta(G)$ and $a \in \{0,1\}$ so that all but at most $b' n^k$ edges $e \in G$ have $|e \cap V_1|=a$ mod $2$. Say that an edge $e$ is {\em good} if $|e \cap V_1|=a$ mod $2$, and $\emph{bad}$ otherwise. We begin by moving any vertex of $V_1$ which lies in fewer than $n^{k-1}/6(k-1)!$ good edges to $V_2$, and moving any vertex of $V_2$ which lies in fewer than $n^{k-1}/6(k-1)!$ good edges to $V_1$. Having made these moves, we update our definitions of `good' and `bad' edges to the new partition (so some edges will have changed from good to bad, and vice versa). We claim that the modified partition has parts of size at least $\delta(G) - bn$, that at most $bn^k$ edges are now bad, and any vertex now lies in at least $n^{k-1}/7(k-1)!$ good edges of $G$. To see this, first note that the minimum degree of $G$ implies that every vertex of $G$ lies in at least $n^{k-1}/3(k-1)!$ edges of $G$. So any vertex of $G$ which we moved originally lay in at least $n^{k-1}/6(k-1)!$ bad edges of $G$; since there were at most $b'n^k$ bad edges in total we conclude that at most $6k!b'n \leq bn$ vertices were moved, giving the bound on the new part sizes. Any edge which is now bad was either one of the at most $b'n^k$ edges which were previously bad, or one of the at most $6k!b'n^k \leq bn^k/2$ edges which contain a vertex we moved, so there are at most $bn^k$ bad edges after the vertex movements. If a vertex was not moved, then it previously lay in at least $n^{k-1}/6(k-1)!$ good edges of $G$; each of these edges is now good unless it is one of the at most $6k!b'n^{k-1}$ edges which also contain another moved vertex. Similarly, if a vertex was moved, then it previously lay in at least $n^{k-1}/6(k-1)!$ bad edges of $G$; each of these edges is now good unless it is one of the at most $6k!b'n^{k-1}$ edges which also contain another moved vertex. In either case we conclude that the vertex now lies in at least $n^{k-1}/7(k-1)!$ good edges of $G$. Having fixed our new partition, we say that a vertex is \emph{bad} if it is contained in at least $k\sqrt{b}n^{k-1}$ bad edges. In particular the set $B$ of bad vertices has size $|B| \leq \sqrt{b}n$. 

Suppose first that $G$ contains a bad edge $e^* = \{u_1, \dots, u_k\}$. Then we may choose vertex-disjoint good edges $e_1, \dots, e_k$ such that $u_i \in e_i$ and the edges $e_i$ contain no bad vertices except possibly the vertices $u_i$. Having done this, greedily choose a matching $E$ in $G$ of size at most $|B|$ which covers all bad vertices and is vertex-disjoint from $e_1, \dots, e_k$ and $e^*$. Now consider the matchings $E_1 := E \cup \{e_1, \dots, e_k\}$ and $E_2 := E \cup \{e^*\}$. We will choose either $E^* = E_1$ or $E^* = E_2$, so that, writing $V'_1$ and $V'_2$ for the vertices remaining in $V_1$ and $V_2$ respectively after deleting the vertices covered by $E^*$, and $n' := |V'_1 \cup V'_2|$, we have that $|V'_1|$ has the same parity as $an'/k$. To see that this is possible, note that $an'/k = an/k - a|V(E^*)|/k$, and $a|V(E_1)|/k - a|V(E_2)|/k = (k-1)a$ mod $2$. On the other hand, $|V(E_1) \cap V_1| - |V(E_2) \cap V_1| = (k-1)a+1$ mod $2$. Since these two differences have different parities, we may choose $E^*$ as required.

Now suppose instead that $G$ has no bad edges. In this case, we take $E^*$ to be empty, and assume that $|V'_1|$ has the same parity as $an'/k$. In either case we will obtain a contradiction to this choice of $E^*$, proving that in fact $G$ must have no bad edges and the size of $V'_1 = V_1$ has different parity to $an/k$, and therefore that $\Part'$ is the partition we wished to find.

Fix such a matching $E^*$, and note that $k \mid n'$ since $n' = n - |V(E)|$. We shall find a perfect matching in $G[V_1' \cup V_2']$, that is, covering the vertices which are not covered by $E^*$. 
Let $I = \{0 \leq i \leq k : i = a \mod 2\}$. We next choose numbers $n_i \ge cn$ for $i \in I$, such that $\sum_{i \in I} n_i = n'/k$ and $\sum_{i \in I} in_i = |V'_1|$. To see that this is possible we use a variational argument: we start with all $n_i$ equal to $cn$ except for $n_0$ or $n_1$ (according as $0 \in I$ or $1 \in I$), which is chosen so that $\sum_i n_i = n'/k$, so $\sum_{i \in I} i n_i$ initially is at most $n'/k + k^2cn \leq |V'_1|$, and has the same parity as $|V'_1|$. We repeatedly decrease some $n_i \geq cn+1$ by $1$ and increase $n_{i+2}$ by $1$; this increases $\sum_{i \in I} i n_i$ by 2, so we eventually achieve $\sum_i in_i = |V'_1|$, as required. Now we choose partitions $\Part_1$ of $V'_1$ and $\Part_2$ of $V'_2$ uniformly at random from those such that $\Part_1$ has $n_i$ parts of size $i$ and $\Part_2$ has $n_i$ parts of size $k-i$ for each $i$. For each $i \in I$, let $X^i_1,\dots,X^i_{n_i}$ be the parts in $\Part_1$ of size $i$, and $Z^i_1,\dots,Z^i_{n_i}$ be the parts in $\Part_2$ of size $k-i$. Then we may form a partition of $V(G)$ into parts $Y^i_j$ for $i \in I$ and $j \in [k]$ uniformly at random, where for each $i$ and $j$ the part $Y^i_j$ contains one vertex from each $X^i_\alpha$ if $j \le i$, and one vertex from each $Z^i_\alpha$ if $i+1 \le j$. So each $Y^i_j$ has size $n_i \geq cn$. We consider for each $i \in I$ an auxiliary $k$-partite $k$-graph $H^i$ on vertex classes $Y^i_1, \dots, Y^i_k$ whose edges are those $k$-partite $k$-tuples that are edges of $G$. So the $H^i$ are vertex-disjoint, and to find a perfect matching in $G$, it suffices to show that with high probability each $H^i$ has a perfect matching. 

Fix some $i \in I$. Since we covered all bad vertices by $E^*$, all vertices in $H^i$ belong to at most $k\sqrt{b}n^{k-1}$ bad edges of $G$. Now fix $y_1 \in Y^i_1$. Say that a $(k-2)$-tuple $(y_j)_{j=2}^{k-1}$ with $y_j \in Y^i_j$ for $j \in \{2, \dots, k-1\}$ is {\em $y_1$-bad} if $y_1 \dots y_{k-1}$ belongs  to more than $kb^{1/4}n$ bad edges. Then there are at most $b^{1/4}n^{k-2}$ $y_1$-bad $(k-2)$-tuples. Suppose that $i < k$, so $Y^i_k \sub V_2$; then a $k$-tuple $y_1 \dots y_{k-1}x$ is good if $x \in V_2$ or bad if $x \in V_1$. Since $\delta(G) \ge (1/2-c)n$ and $|V_2| \le (1/2+c + b)n$, if $(y_i)_{i=2}^{k-1}$ is $y_1$-good then $y_1 \dots y_{k-1}x$ is an edge for all but at most $3cn$ vertices $x \in V_2$. The number of these vertices $x$ that lie in $Y_k^i$ is hypergeometric with mean at least $(1-3c)|Y_k^i|$, where $|Y_k^i|=n_i \ge cn$. So by the Chernoff bound (Lemma~\ref{chernoff}), with probability $1 - o(1)$ there are at least $(1-4c)n_i$ such vertices $x$ for any choice of $y_1 \in Y^i_1$ and $y_1$-good $(y_i)_{i=2}^{k-1}$. Thus $y_1$ is contained in at least $$(n_i^{k-2} - b^{1/4}n^{k-2})(1-4c)n_i > (1-5c)n_i^{k-1}$$ edges of $H^i$. By symmetry, with positive probability the same bound holds for all vertices in any $H^i$ (for $i=k$ we have $Y_k \sub V_1$, in which case we proceed similarly with the roles of $V_1$ and $V_2$ switched). Then each $H^i$ has a perfect matching by Theorem~\ref{partitevertexdegree}, contradicting our assumption that $G$ does not have a perfect matching. So we must have the case that no edges of $G$ are bad and $V_1$ has different parity to $an/k$, as required.\endproof

\chapter{The general theory} \label{sec:general}

In this final chapter, we give a more general result, which epitomises the geometric nature of the theory, in that it almost entirely dispenses with degree assumptions. We cannot hope to avoid such assumptions entirely when using methods based on hypergraph regularity, as sparse hypergraphs may have empty reduced systems. Our degree assumptions are as weak as possible within this context, in that the proportional degrees can be $o(1)$ as $n \to \infty$. The hypotheses of our result are framed in terms of the reduced system provided by hypergraph regularity, so it takes a while to set up the statement, and it is not as clean as that of our main theorems (which is why we have deferred the statement until now). However, the extra generality provided by the geometric context will be important for future applications, even in the context of minimum degree thresholds for hypergraph packing problems, where the tight minimum degree may not imply the minimum degree sequence required by our main theorems. (We intend to return to this point in a future paper.) In an attempt to avoid too much generality, we will restrict attention to the non-partite setting here.

First we describe the setting for our theorem. We start with the hypergraph regularity decomposition.\medskip

\nib{Regularity setting.}
Let $J$ be a $k$-complex on $n$ vertices, where $k \mid n$.
Let $\Qart$ be a balanced partition of $V(J)$ into $h$ parts,
and $J'$ be the $k$-complex of $\Qart$-partite edges in $J$.
Let $P$ be an $a$-bounded $\eps$-regular vertex-equitable $\Qart$-partition $(k-1)$-complex on~$V(J)$,
with clusters $V_1,\dots,V_{m_1}$ of size $n_1$.
Let $G$ be a $\Qart$-partite $k$-graph on~$V(J)$ that is $\nu$-close to~$J'_k$ and perfectly $\eps$-regular with respect to~$P$.
Let $Z = G \bigtriangleup J'_k$. \medskip

Next we describe the setting for the reduced system, which names the constructions that were used in the proof of Theorem~\ref{perfectmatchingF}. While we do not specify the source of the subsystem $(R_0,M_0)$, it is perhaps helpful to think of it as being randomly chosen so as to inherit the properties of $(R,M)$, as in the proof of Lemma~\ref{RANDOMEDGESELECTION}.\medskip
 
\nib{Reduced system setting.}
Let $R^1 := R^{J'Z}_{P\Qart}(\nu,\cb)$ on $[m_1]$. Also let $(R,M)$ be a matched $k$-system on $[m]$, where $R$ is the restriction of $R^1$ to $[m]$. Given any $R'$ and $M'$ which can be formed from $R$ and $M$ respectively by the deletion of the vertices of at most $\psi m$ edges of $M$, and any vertices $u,v \in V(R')$, we will choose $M_0 \sub M'$ such that $V_0 := V(M_0)$ satisfies $m_0 := |V_0| \leq m_0^*$ and $M_0$ includes the edges of $M'$ containing $u,v$. Let $R_0 = R[V_0]$, $X=X(R_0,M_0) = \{\chi(e)-\chi(e'): e \in R_0,e' \in M_0\} \sub \R^{m_0}$
and $\Pi = \{\xb \in \Z^{m_0} : \xb \cdot \1 = 0 \}$. \medskip

Now we can state our general theorem (in the above setting). 

\begin{theo} \label{general}
Suppose that $k \geq 3$ and $1/n \ll \eps \ll 1/a \ll \nu, 1/h \ll c_k \ll \dots \ll c_1 \ll \psi \ll \delta, 1/m_0^*, 1/\ell \ll \alpha, \gamma, 1/k$.
Suppose that $m \geq (1-\psi) m_1$, and the following conditions hold.
\begin{enumerate}[(i)]
\item For any $R', M', u$ and $v$ given as above, we may choose $M_0$ so that $(R_0, M_0)$ satisfies $B(\0, \delta) \cap \Pi \sub CH(X)$.
\item $\delta^+(D_\ell(R,M)) \ge \gamma m$.
\item $L_\Part(R_k)$ is complete for any partition $\Part$ of $V(R)$ into parts of size at least $(\gamma-\alpha)m$.
\item Every vertex is contained in at least $\gamma n^{k-1}$ edges of $J$.
\end{enumerate}
Then $J_k$ contains a perfect matching.
\end{theo}

\proof
We follow the proof of Theorem~\ref{perfectmatchingF}, outlining the modifications. Introduce constants with
\begin{align*}
& 1/n \ll \eps \ll d^* \ll d_a \ll 1/a \ll \nu, 1/h \ll d, c \ll c_k \ll \dots \ll c_1\\ 
&  \ll \psi \ll 1/C' \ll 1/B, 1/C \ll \delta, 1/m_0^*, 1/\ell \ll \alpha, \gamma,1/k.
\end{align*}
Our set-up already provides $\Qart$, $J'$, $P$, $G$, $Z$, $R$, $M$, as in that proof.
(Note that we are working in the non-partite setting, so there is no $\Part$, and $F$ consists of the unique function $f: [k] \to [1]$.
We also do not need to consider two reduced systems.) There is no need for an analogue of Claim~\ref{matchingclaim2} as we already have $M$. 
The proof of Claim~\ref{claim_partitionJ} is the same, except that for (iii), instead of
the minimum $F$-degree assumption, we apply Lemma~\ref{randomsplit} to the edges of $J$ containing $v$,
using the assumption that every vertex is contained in at least $\gamma n^{k-1}$ edges of $J$.
To prove Claim~\ref{completetransdigraphF}, instead of Lemma~\ref{main-reduction-partiteF}, we first apply
Lemma~\ref{ballgivesgentransferrals} to see that $(R_0,M_0)$ is $(B,C)$-irreducible for any $u$ and $v$. It follows that $(R',M')$ is $(B,C)$-irreducible, then we apply Lemma~\ref{complete_digraph} to see that $D_{C'}(R',M')$ is complete.
The remainder of the proof only uses these claims, so goes through as before.
\endproof

Each of the conditions of Theorem~\ref{general} is necessary for the proof strategy used in Theorem~\ref{perfectmatchingF}, so we cannot strengthen this result further using this argument. 

We conclude with a remark comparing our techniques to the absorbing method, which has been successfully applied to many hypergraph matching problems. The idea of this method is to randomly select a small matching $M_0$ in a $k$-graph $G$, and show that with high probability it has the property that it can `absorb' any small set of vertices $V_0$, in that there is a matching covering $V(M_0) \cup V_0$. Given such a matching $M_0$, to find a perfect matching in $G$, it suffices to find an almost-perfect matching in $G \sm V(M_0)$, which is often a much simpler problem.

For example, the essence of Lo and Markstr\"om's independent proof of Theorem~\ref{partitehajnalszem} was to show that, in any $k$-partite graph $G$ which satisfies the degree conditions of this theorem, a randomly-chosen collection $M_0$ of $\eps n$ disjoint copies of $K_k$ in $G$ is absorbing. In this setting this means that for any set $V_0 \subseteq V(G)$ whose size is divisible by $k$, there exists a perfect $K_k$-packing in $G[V_0 \cup V(M_0)]$. They then showed that $G \sm V(M_0)$ must contain an $K_k$-packing $M_1$ covering almost all the vertices; the absorbing property of $M_1$ then implies that $G \sm V(M_1)$ contains a perfect $K_k$-packing, which combined with $M_1$ yields a perfect $K_k$-packing in $G$. See~\cite{LM2} for further details. 

The advantage of the absorbing method is that it avoids the use of the hypergraph blow-up lemma (and often avoids hypergraph regularity altogether), and so leads to shorter proofs, when it works. However, the existence of an absorbing matching seems to be a fortuitous circumstance in each application of the method, and it is not clear how one could hope to find it in general problems. By contrast, our theory explains `why' there is a perfect matching, by analysing the only possible obstructions (space barriers and divisibility barriers). This has the additional advantage of giving structural information, so we can apply the stability method to obtain exact results.

\medskip

\textbf{Acknowledgements.}
We would like to thank Allan Lo for pointing out an inconsistency in the proof of Theorem~\ref{partitehajnalszem} in an earlier version, and also an anonymous referee for offering helpful suggestions which have led to an improvement in the presentation of this paper.

\backmatter
\bibliographystyle{amsalpha}



\end{document}